\documentclass[12pt]{article}
\pdfoutput=1 

\usepackage[final]{epsfig}
\usepackage{ulem}
\usepackage{ltxtable}

\usepackage[hang]{subfigure}

\def\proof{\noindent {\bf Proof. }}

\def\R{I\!\!R}

\newcommand{\B}{\mathcal{B}} 
\newcommand{\Ho}{\mathcal{H}}

\newcommand{\Ghat}{\widehat{G}}

\newcommand{\G}{\widetilde{G}}

\newcommand{\poly}{\mathcal{P}}
\newcommand{\Phat}{\widehat{\mathcal{P}}}
\newcommand{\Qhat}{\widehat{\mathcal{Q}}}

\newcommand{\PEhat}{\widehat{\mathcal{P}^{E}}}
\newcommand{\PBhat}{\widehat{\mathcal{P}^{B}}}

\newcommand{\T}{\mathbf T}

\newcommand{\tower}{\mathcal{T}}

\newcommand{\cylinder}{\mathcal{C}}

\newcommand{\mc}{\mathcal}

\newsavebox{\cuadrito}
\sbox{\cuadrito}{\framebox[7pt]{ }}

\newcommand{\qed}{\makebox[8pt]{}\hfill {\usebox{\cuadrito}}\medskip}

\newtheorem{theorem}{THEOREM}[section]

\newtheorem{definition}{Definition}[section]
\newtheorem{lemma}[theorem]{LEMMA}
\newtheorem{proposition}[theorem]{PROPOSITION}

\newtheorem{conjecture}{CONJECTURE}[section]
\newtheorem{corollary}[theorem]{COROLLARY}

\newtheorem{algorithm}[theorem]{ALGORITHM}

\newcounter{mycomplaints}
\def\complaint#1{\refstepcounter{mycomplaints}%
\ifhmode%
\unskip%
{\dimen1=\baselineskip \divide\dimen1 by 2 %
\raise\dimen1\llap{\tiny -\themycomplaints-}}\fi%
\marginpar{\tiny [\themycomplaints]: #1}}%

\newcounter{example}
  \newenvironment{example}{\refstepcounter{example}
    \subsubsection{Example
     }}{$\qed$ \\}

\usepackage{url}

\begin{document}

\title{Isostatic Block and Hole Frameworks}
\author{
{Wendy Finbow-Singh \thanks{Department of Mathematics, St. Mary's University, Wendy.Finbow-Singh@sma.ca
Supported in part by a grant from NSERC, and by York University}} 
\\
{ Walter Whiteley
\thanks{Department of Mathematics and Statistics, York University,    
whiteley@mathstat.yorku.ca.  Supported in part by a grant
from  NSERC (Canada).  }}
}
\maketitle
\begin{abstract} A longstanding problem in rigidity theory is to characterize the graphs which are minimally generically rigid in $3$-space.
The results of Cauchy, Dehn, and Alexandrov give one important class: the triangulated convex spheres, but there is an ongoing desire for further classes.
We provide such a class, along with methods to verify generic rigidity that can be extended to other classes.  
These methods are based on a controlled sequence of vertex splits, 
a graph theoretic operation known to take a minimally generically rigid framework to a new minimally generically 
rigid framework with one more vertex.  

One motivation for this is to have well-understood frameworks  which can be used to explore Mathematical Allostery - frameworks in which adding bars at one site, 
causes changes in rigidity at a distant site.  This is an initial step in exploring the possibility of mechanical models for an important behaviour in proteins.



\end{abstract}

\section{Introduction }

The combinatorial characterization of general graphs that can be realized 
in 3-space as isostatic
(rigid and independent) bar and joint frameworks is a major unsolved problem in
rigidity theory \cite{graver,taywhiteley,wchapter,whandbook}.  In the absence of a general characterization,
it is significant to investigate additional classes of graphs and confirm the rigidity
and independence 
of almost all realizations of such graph in 3-space (generic rigidity), something that follows from
the existence of even one isostatic realization.

Historically, from the work of Cauchy and Dehn \cite{cauchy,dehn},
we know that arbitrary convex triangulated spheres are isostatic,
 and therefore any realization of these 3-connected planar graphs $G=(V,E)$,
 with $|E|=3|V|-6$, at generic positions of the vertices
 is also isostatic \cite{gluck}.  More generally, the results of Alexandrov \cite{alexandrov}
 show that for any convex polyhedron, one can subdivide the natural edges
 with additional  vertices, and then triangulate the `faces' on all their vertices,
 to create an isostatic framework~\cite{infp1}.  From one point of view, the steps from this
 example becomes a guide for the steps in some of the proofs we explore in this paper. 

\begin{figure}[ht]
 \begin{center}
 \psfig{file=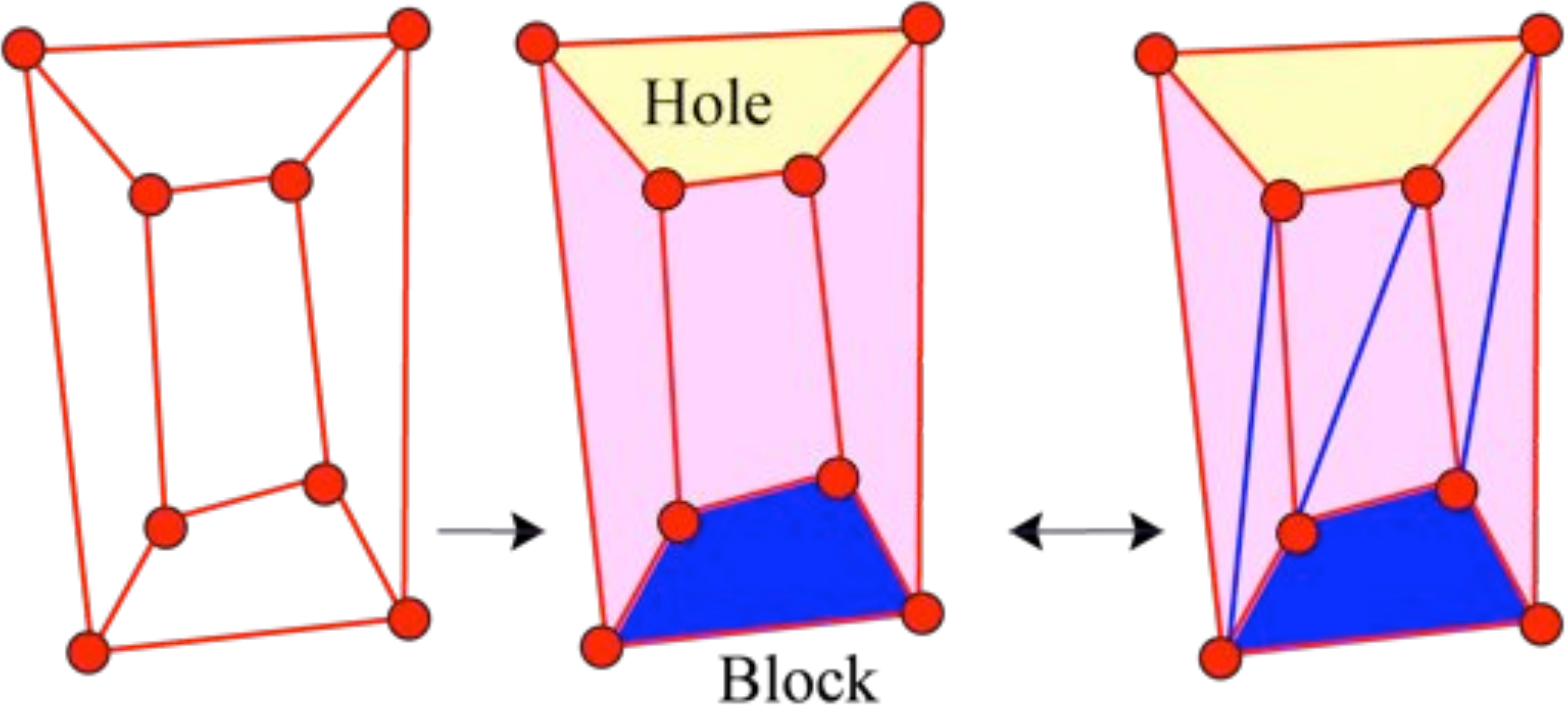, scale=.08}
  \caption{Creating a block and hole polyhedron.  Moving from a 3-connected
 planar graph, (a), through selecting faces as blocks and holes, (b), to triangulating the
the surface discs, (c).\label{PolyBlock}}
 \end{center}
 \end{figure}

A second model for the results in this paper is Theorem 4.1 of Whiteley, \cite{infp2}, which
verified a generic version of a conjecture of Kuiper  \cite{kuiper}.   
If we start with a triangulated sphere 
and remove one edge (creating a quadrilateral {\it hole}) and insert
a new edge somewhere
else, connecting the two vertices of triangles which share an edge (creating a tetrahedral {\it block})
then the resulting graph is isostatic for almost all realizations provided
the hole and the block are connected by four vertex disjoint paths from the four vertices
of the hole to the four vertices of the block (Figure~\ref{PolyBlock}).

In this paper, we are interested in the behaviour of general triangulated spheres modified with some holes and some isostatic blocks.  
The holes and blocks are created in such a way as to balance the count; 
that is, the number of edges removed from the holes will equal the number of edges added to form the blocks, so $|E| = 3|V| - 6$.

For example, when  building a geodesic dome, you cut off one part of the triangulated convex sphere and attach this base to the ground.  The ground is a large block.  If the polygon on the ground has $k$ vertices, you can now remove up to $k-3$ other edges to make windows, doors, or other holes in the remaining dome.  Which edges can you remove?  We will provide methods which can be applied to verify the rigidity of a proposed set of added holes, at least for generic positions of the vertices. 

\begin{figure}[h] 
    \begin{center}
    \psfig{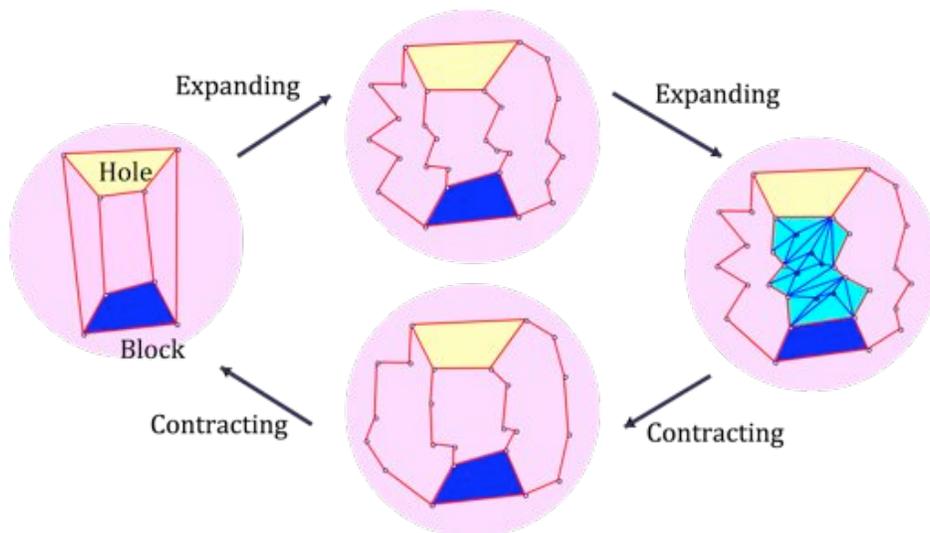}
    \caption{Expanding a block and hole polyhedron from the base (a)
 through expanding the edges between surface discs (b) to adding additional vertices
 and triangulating the surface discs (c). Carefully selected contractions are used to find a sequence of vertex splits which grow the base $\Phat$ to the expanded polyhedron $\PEhat$. \label{PolyExpand}}
 \end{center}
 \end{figure}

We cast our solution in a pattern reminiscent of some steps of Alexandrov's Theorem
(Figure~\ref{PolyExpand}):

\begin{enumerate}
 \item[(i)]  pick a 3-connected planar graph $G$ - the graph of the polyhedron  $\mc{P}$;
 
  \item[(ii)]  select certain faces to become holes (no added interior edges in the face) and
  other faces to become blocks (inserting isostatic subframeworks on the vertices of the face) to form
  a graph $\widehat{\mc{P}}$, leaving the remaining faces as {\it triangulated surfaces faces}, creating a 
  {\it base polyhedron};
  
  \item[(iii)]  expand the graph by inserting new vertices along edges which do not
  belong to a hole or a block (edges which separate surface faces);
  
 \item[(iv)]  further expand the vertices of this polyhedron by inserting vertices inside any surface faces
  and provide an arbitrary planar triangulation of these surface `faces' (including all
  added vertices of the face), creating a graph $G^{*}$ for the expanded polyhedron $\PEhat$;
  
 \item[(v)] verify that the expanded polyhedron  $\PEhat$ can be reached from the base polyhedron through a sequence of vertex splits by carefully selecting the reverse sequence of operations of contracting edges;
 
  \item[(vi)] verify that the base polyhedron $\Phat$, with blocks and holes and each of the possible 
 triangulations of the surface faces, is isostatic;
 
 \item[(vii)] conclude that the expanded polyhedron $\PEhat$ is isostatic for almost all (generic) realizations, and arbitrary retriangulations of the expanded faces.
 \end{enumerate}

Our basic tool for connecting the base polyhedron with the expanded polyhedron is
vertex splitting \cite{vertexsplit}.  This is a technique developed to expand
isostatic frameworks by inserting a new vertex and two triangles in place of two adjacent
edges, creating a larger isostatic framework.  It was first extensively applied to verify the generic rigidity of 
 fully triangulated surfaces \cite{fogelsanger, vertexsplit}.  

Knowing that vertex splitting preserves first-order rigidity and independence as we expand, 
the central problem is to show that a given final expanded structure $G^{*}=\PEhat$ comes from some 
triangulation of the base polyhedron $\widehat{P}$ by a 
sequence of vertex splits.  This paper will present some new techniques to confirm
the existence of such a construction sequence, based on the local topology of how the triangulated
discs are embedded in the larger framework, as well as how two such discs and 
their shared boundary path are embedded in the larger framework. 

Given a triangulated disc,  we verify the existence of the 
desired vertex splits by working back via the reverse of vertex splitting - edge contraction,
for three cases:
\begin{enumerate}
 \item  inductively contract edges so that no interior edges remain in a triangulated disc
 (Section~3.1)
 \item  inductively contract edges so that no vertices remain interior to triangulated discs
 (Section~3.2);
 \item  inductively contract the path vertices interior to the original edges between triangulated
 discs, so that we have a triangulation of the original surface faces of the base polyhedron 
 (Section~3.3) .
\end{enumerate}

The conclusion will be that the expanded polyhedron is constructed from the base by a sequence of vertex splits.  
Therefore, if the 
base polyhedron is generically rigid (independent, isostatic), then the final polyhedron will be generically  rigid.

 This leaves the second challenge of creating an appropriate isostatic base polyhedron (Section 5), and ways to `find'
 an appropriate base polyhedron in a larger framework which has been formed within a large
 triangulated sphere by selecting holes and blocks in an appropriate balance: a
{\it block and hole polyhedron}.   For example,
 for our new proof of Kuiper's Conjecture, we need 4-connectivity (in a vertex sense)
 to ensure that we 
 could find the base `quadrilateral tower' and complete the proof (\S5.2). 

 In Section~\ref{sec:base} we also provide a general conjecture for
 when such a block and hole polyhedron is isostatic.  This conjecture would imply that this
 class of frameworks can be tested for being isostatic by a simple polynomial time counting
 algorithm and a local condition, in contrast to the general problem of frameworks in 3-space.

 For the last five years, the class of frameworks we are investigating in \S5 has been applied in modified form as a
 `toy template' for a rigidity approach to allostery (transmission of shape and rigidity across large
 biomolecules \cite{allostery}).  With this modeling of `transmission of degrees of freedom'
 across frameworks in mind, we purse some extensions of this work (\S6).  For example, 
 if we started with two pentagons, as a block and a hole, we will see that making both of
 these into holes creates a residual framework with four `shared' degrees of freedom.
 In a hypothetical docking of a ligand on one site, under the assumption of 5-connectivity
 between the sites, we will confirm we can remove  0, 1, 2, 3 or 4
 degrees of freedom from both sites.

For many purposes, the blocks and holes  are simply `reserved faces` which will not be altered during the contractions and vertex splits, without paying any attention to what is happening inside, behind or even between those faces.  In Section~\ref{sec:cycle} we step back to observe that the methods and results of the Sections 4,5,6 contain a general splitting process which can be applied to an arbitrary generically rigid graph $G$:  {\it Given any cycle, this can be split into two cycles, and a general triangulated cylinder inserted, creating a larger generically rigid graph}.   Then for many purposes the blocks and holes have simply been reserving selected faces which will not be altered, without paying any attention to what is happening inside or `behind' those faces.

Throughout this paper, the reader might notice that the conditions on blocks and on
holes are equivalent, suggesting that if a structure  $G^{*}$ is generically isostatic,
with a set of blocks $\B$ and holes $\Ho$, then the swapped framework $\widetilde G$
will also be isostatic with blocks $\Ho^{*}$ and holes $\B^{*}$.
In a companion paper \cite{frw}, we verified this swapping
in a strong geometric form:  {\it If a block and
hole framework $(G^{*},p)$ is isostatic, then the swapped framework $\widetilde (G,p)$ is also
isostatic, for the same positions of the polyhedral vertices.}  When it helps organizing the pieces here, we will recall this swapping within this the paper.

\section{Basic Vocabularly and Concepts}\label{sec:basic}

A {\it graph} $G = (V, E)$ is an ordered pair of sets $V$ and $E$, where $V = \{1, 2, \ldots, v\}$ is a set of {\it vertices}, 
and $E$ is a set of unordered pairs of distinct vertices called {\it edges}.  

A {\it path} is a sequence $(x_1, x_2,  \ldots, x_n)$, $1\leq n$,  of distinct vertices of a graph such that 
consecutive vertices in the sequence are edges of the graph.  We define the {\it length of a path} to be the number of edges in the path.  
A {\it cycle} is a path $(x_1, x_2,  \ldots, x_n)$, $1< n$ such that $x_1 = x_n$.
Two vertices, $x$ and $y$, of a graph $G$ are said to be {\it connected} if there is a path from $x$ to $y$ in the graph.  If all pairs of vertices are connected, the graph $G$ is {\it connected}.  Otherwise the graph $G$ is  {\it disconnected}.   
A graph $G$ is called {\it $k$-connected} if there does not exist a set of $k-1$ vertices whose removal disconnects the graph. For example, a simple cycle is $2$-connected.  

A {\it planar graph} $G$, is a graph that can be drawn in the plane (or sphere) such that no pair of edges of $G$ intersect except at their endpoints. 
When drawn on a sphere, a 3-connected planar graph creates a unique set of disjoint regions with no interior vertices, each surrounded by a boundary polygon of vertices and edges \cite{whitney}.  We call these regions, with their boundary polygons the {\it faces} of the {\it topological graph} $\G$. Since these $3$-connected planar graphs can also be realized as polyhedra (by Steinitz's Theorem), we will describe them as polyhedra.  For example, the graph $K_4$ is visualized as a spherical tetrahedron with four triangular faces.  

Given such a topological graph, we designate some faces as the set of holes $\Ho$, some other faces to be a set of blocks 
$\B$, and the remaining faces will be the set  $T$ (which  will become the triangulated discs).  Notice that in this partition into faces induced by a $3$-connected planar graph, two  faces  $F_1$ and $F_2$ intersect in: (i) the empty set; or (ii) one vertex; or (iii) in a single edge.   A topological graph with such a partition of the faces is denoted by both $\Ghat$ and $\poly$.

A {\it triangulated face} of a topological graph is a face which has been triangulated on its boundary polygon.  When we triangulate all faces in $T$, and insert an isostatic subgraph into each block (see below), we have created the {\it block and hole polyhedron $\Phat$}.

Two such block and hole polyedra $\Phat, \Qhat$  are called {\it topologically equivalent} if they have the same blocks $\B$ and holes $\Ho$, and the same face polygons $T$, with possibly different triangulations of the faces in  $T$.

\subsection{Generic rigidity}
Our main theorems will confirm the `generic rigidity' of certain graphs based on block and hole polyhedra.  We also want to insert a `generically rigid' graph for each of the blocks in our polyhedron.  
We build up the definition of generic rigidity starting with infinitesimal rigidity of frameworks $(G,p)$ at a specific geometric configuration $p$, and 
then shift to `generic' configurations \cite{graver,taywhiteley,wchapter,whandbook}.  We will present all the definitions and results  for $\R^3$.

A {\it framework} $\mathcal{F} = (G, p)$ in $\R^3$ is a graph $G$, together with configuration $p$ (an embedding of the vertices), 
$p: V \rightarrow \R^3$.  An {\it infinitesimal motion} of the framework 
$\mathcal{F} = (G = (V, E), p)$, is a function $\nu: V \to \R^3$ such that 
for every edge $\{x, y\} \in E$
$$(p(x) - p(y)){\bf \cdot}(\nu(x) - \nu(y)) = 0.  $$
Assuming that the configuration $p$ spans $\R^3$ \cite{taywhiteley,wchapter}, an  {\it infinitesimal rigid motion}  is a function $\nu^*: V \to \R^3$ such that 
for every pair of vertices $x$, and $y$ in $V$
$$(p(x) - p(y)){\bf \cdot}(\nu^*(x) - \nu^*(y)) = 0.  $$
Other infinitesimal motions which are not rigid motions are called {\it infinitesimal flexes}. 
The framework $\mathcal{F}$ is {\it infinitesimally rigid} if every infinitesimal motion of the framework is an infinitesimal rigid motion. 

The framework $(G,p)$ is {\it independent} if removing any one edge increases the space of infinitesimal motions. 
A framework $(G,p)$ is {\it isostatic} if it is infinitesimally rigid and independent.  
In other words, the framework is called isostatic if it is minimally infinitesimally rigid.

The rigidity properties can be represented in the linear algebra of the homogeneous equations above.  In this language the rigid frameworks form a basis for the row space of the matrix for the complete graph, the independent frameworks are independent sets of rows, and the isostatic frameworks are bases. 
An isostatic framework in $\R^3$, with  $|V| > 3$, will have  $|E|=3|V| - 6$ edges.  For example, all convex triangulated spheres are isostatic \cite{dehn}.  \\

 Throughout this paper we will be focusing on graphs whose frameworks are isostatic for `almost all' configurations.  
We will be working with the combinatorics of the graph $G$, not with any special properties of the configuration $p$. 
A configuration $p$ is  {\it generic} if  the coordinates are not the solution to any algebraic equations.  The generic configurations are an open 
dense subset of all configurations in  $\R^{3|V|}$ - so we sometimes say `almost all configurations' (those found with probability one if we make random 
choices of the coordinates). 
A graph $G$ is  {\it generically rigid} if every framework with generic vertices is rigid.  The theory guarantees that a graph $G$ is generically rigid if and only if some framework $(G,p)$  is infinitesimally rigid \cite{graver,wchapter}.  

Similarly, a graph $G$ is  {\it generically independent} if the graph is independent in some framework $(G,p)$.   A graph $G$ is  {\it generically isostatic} if it is generically rigid and generically independent.  Our goal is to prove that a range of graphs $G$ are generically isostatic.  

We have imaged the blocks as faces which are themselves generically rigid, when restricted to the vertices of their boundary cycle.  The next theorem says we can substitute one generically rigid block for another one, so that which generically rigid subframework insert does not matter.  

\begin{theorem} [General isostatic Substitution Principle \cite{frw}]
Given a generic framework $\mc{F}$ with an isostatic subframework $\mc{F}'$, then replacing $\mc{F}'$ with another
isostatic subframework $\mc{F}''$ gives a new framework $\mc{F}^*$ on the same vertices which has an isomorphic space of infinitesimally motions. In particular, if $\mc{F}$ is isostatic (independent, rigid) then  $\mc{F}^*$  is isostatic (independent, rigid, resp). 
\end{theorem}

With this in hand we will assume that some edges  have been added to the boundary cycle of each block, so that the subframework on these vertices is generically isostatic.  This was a step in our definition of the block and hole polyhedron $\Phat$ above.  It will not matter which generically isostatic subgraph we add for a block, as long as it uses the edges of the boundary cycle of the block.   

\subsection{Expansions of polyhedra}  
Given a block and hole polyhedron $\Phat$, with blocks, holes, and triangulations in the faces of $T$, we can expand it to a larger block and hole polyhedron $\PEhat$, as informally described in the introduction (Figure~{\ref{PolyExpand}}).   There are several steps to creating the larger graph:
\begin{itemize}
\item[(i)] remove the edges interior to the triangulation of boundary polygons of $T$, leaving the boundary polygons; 

\item[(ii)] subdivide edges which separate two faces in $T$, creating expanded boundary polygons for these faces;

\item[(iii)] add vertices interior to  individual faces of $T$
 
\item[(iv)] select new triangulations for each of the faces  $F\in T$, using the new boundary polygon and all the new  vertices interior to the face.
\end{itemize}
We call the new polyhedron $\PEhat$, and abuse notation to continue to describe the set of modified triangulated faces as $T$.  We call a face $F\in T$ together with its interior vertices and the triangulation a {\it triangulated disc} $D$.  

Our goal is to prove results which apply for all choices the new triangulations of the extended faces in $T$. We note there is still a shared `face' structure of $\Phat$ and $\PEhat$.   Any two such expansions $\PEhat$ and ${\PEhat}{'}$ have the same blocks and holes.  They have {\it topologically equivalent faces $T$} - in the sense that if we remove all but the boundary edges of each triangulated disc and if we identify the expanded paths as single abstract edges we have the shared underling topology, including blocks and holes of $\Phat$.   In particular, two faces either: (i) do not intersect; (ii)
intersect in a single vertex; or (iii) intersect in a single path. We say that  $\Phat$,  $\PEhat$ and ${\PEhat}{'}$ are all {\it topologically equivalent}. 

For later reference, an edge is said to be in the {\it interior} of a triangulated  disc $D$ if neither of it's vertices lies on the polygonal boundary of the disc.
Other edges in the triangulation of $D\in T$, that is edges with one or both vertices on the boundary, are not interior.  Edges on the intersection of two discs, including the faces of a hole or a block are called {\it boundary edges}. 


\subsection{Contracting edges}  
The operation of contracting an edge of a graph will be central to Sections 3 and 4 of this paper.  
Informally, edge contraction is an operation on a graph which merges the vertices of an edge, say $\{x, y\}$, into a new vertex $z$ of the graph.  
The edge being contracted is removed from the graph, and any vertices adjacent $x$ or $y$ before the merge become adjacent to the new vertex afterwards.   In this paper, contraction will only be applied to edges of discs in $T$, not to edges in blocks, or the boundaries of blocks or holes. 

We state this formally as follows:
Let $G = (V, E)$ be a graph such that $\{\{x, y\}, \{x, u_1\}, \{x, u_2\}, \ldots, \{x, u_m\}$, $ \{y, v_1\}, \{y, v_2\}, \ldots, \{y, v_n\} \} \subset E$. 
The graph resulting from {\it contracting the edge} $\{x, y\}$ to a single vertex $z$ is the graph $G' = (V', E')$, where
$V' = V \backslash \{x, y\} \cup \{z\}$, and $E' = E \backslash \{ \{x, y\}, \{x, u_1\}, \{x, u_2\}, \ldots, \{x, u_m\}$, $\{y, v_1\}, \{y, v_2\}, 
\ldots, \{y, v_n\} \} $$
\cup \{ \{z, u_1\}, \{z, u_2\}, \ldots, \{z, u_m\}$, $\{z, v_1\}$, $\{z, v_2\}, \ldots $, $\{z, v_n\} \}$.  


When contraction is applied to an edge of $\Phat$, it is possible that the topology will be changed.  
A {\it non-facial triangle}, or {\it non-surface triangle} $\T$ of a polyhedron, $\Phat$, 
is a $3$-cycle whose vertices and edges are vertices and edges of 
$\Phat$ (not inside any block), however, the triangle itself is not a face of $\Phat$, or included as a triangle in the given triangulations of discs $D\in T$.      
Consider the corresponding planar graph, $G$, of $\Phat$.  $\T$ will be a separating cycle in $\Phat$; that is, 
 $\T$ contains at least one vertex in it's interior (with respect to the embedding in the plane) and one vertex in it's exterior.  
A non-facial triangle, $\T$, of $\Phat$ can be thought of as a triangular `waist' of $\Phat$.
We call the edges of any non-facial triangles of $\Phat$ {\it short edges}, and all other edges of $\Phat$ {\it long edges}. 

We may contract a long edge of $\Phat$ which is inside a disc, or is on a boundary path of length greater than one, without changing the topology of $\Phat$.  
If we contract an edge of boundary path of length one separating two discs, this keeps the spherical topology, but changes the underlying face structure, creating a new polyhedron $\Phat{'}$ which is not topologically equivalent, as we have defined this. 

If a short edge in a polyhedron $\Phat$ is contracted, this will leave  a single edge separating the larger graph with two smaller polyhedra.  This would definitely alter the topology.  Throughout Section~3 we will be working to avoid contractions on short edges and preserve the topology of the polyhedron.

\section{Inductions to Contract Discs and Paths}

In this section we are working with a block and hole polyhedron $\Phat$ which can be viewed as an expansion $\PEhat$ of a base polyhedron $\PBhat$, as described in Sections~1,2.
In this setting, the intersection of two surface discs $D_{1}, D_{2}$ in $\Phat$ corresponds to an edge of the base polyhedra (perhaps currently an expanded path), and the 
intersection of three or more surface discs of the polyhedron corresponds to a vertex of the base polyhedra. 
Throughout this section, we impose conditions that apply when we are working with the types of triangulated discs which arise in this context. However, the inductive steps have wider applications, so we will be explicit about the local conditions on a disc, or on a pair of discs in the following lemmas.  

Our goal is to start with a disc, or a set of discs and contract edges to produce a smaller polyhedron, aiming towards a simpler base polyhedron. We do this in sequence:
\begin{enumerate}
\item  contracting to remove all edges between interior vertices of a disc (Lemma~\ref{shrink_interior_edge_lemma}); 
\item  contracting to remove all remaining vertices inside a disc (hubs of wheels)  (Lemma~\ref{spoke contraction});
\item contracting the lengths of the path between two surface discs to a single edge (Lemma~\ref{pathshrink}). 
\end{enumerate}
When done over an appropriate set of discs, this will reverse the expansion process  back to the original
polyhedron - based only on the combinatorics within the discs and between the discs.  This contraction sequence will also yield a reverse sequence of vertex splits (see \S4).  For later reference, we will call this sequence of inductively contracting the edges of the triangulated sphere a {\it Contraction Sequence}.

\subsection {Contracting interior edges of a triangulated disc}

Throughout this section, we assume $D$ to be a triangulated disc in a graph $G$.  The key properties of $D$ used are:
\begin{itemize}
\item[(i)] all edges in the large graph $G$ connecting to interior vertices are in the triangulation; 
\item[(ii)] the interior of the disc is fully triangulated, that is, every edge non on the boundary of $D$ is part of two triangles, and 
\item[(iii)] every triangle formed by three vertices of the disc, at least one interior to the disc, is either part of the triangulation or separates the vertices of the disc. 
\end {itemize}
We first show that if there is an edge joining two interior vertices of $D$, then we can
find a long edge in the interior of $D$ to contract. 

\begin{figure}[h]
\begin{center}
  \subfigure[] { \psfig{file=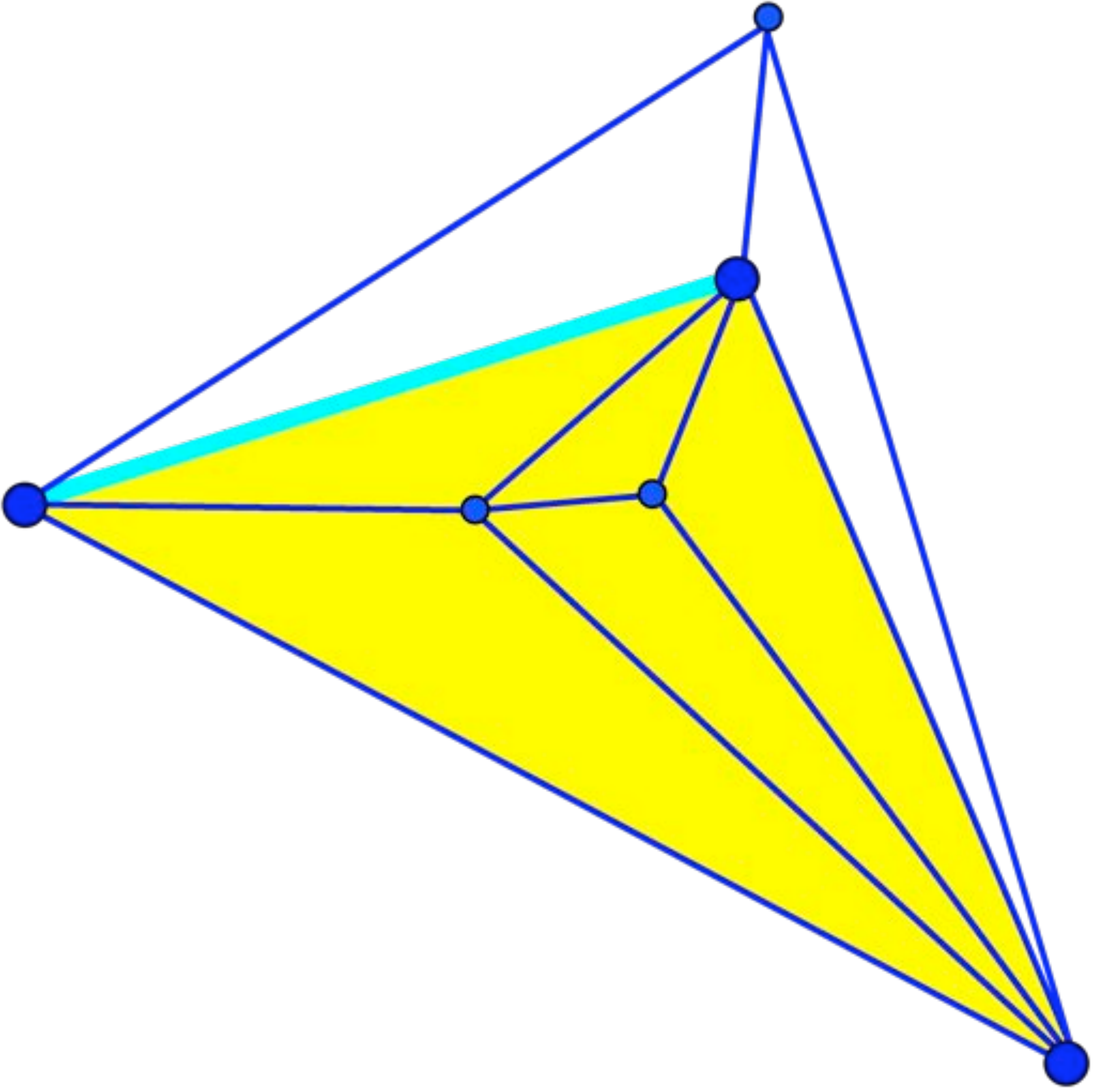, scale=.035}}\quad\quad
  \subfigure[] { \psfig{file=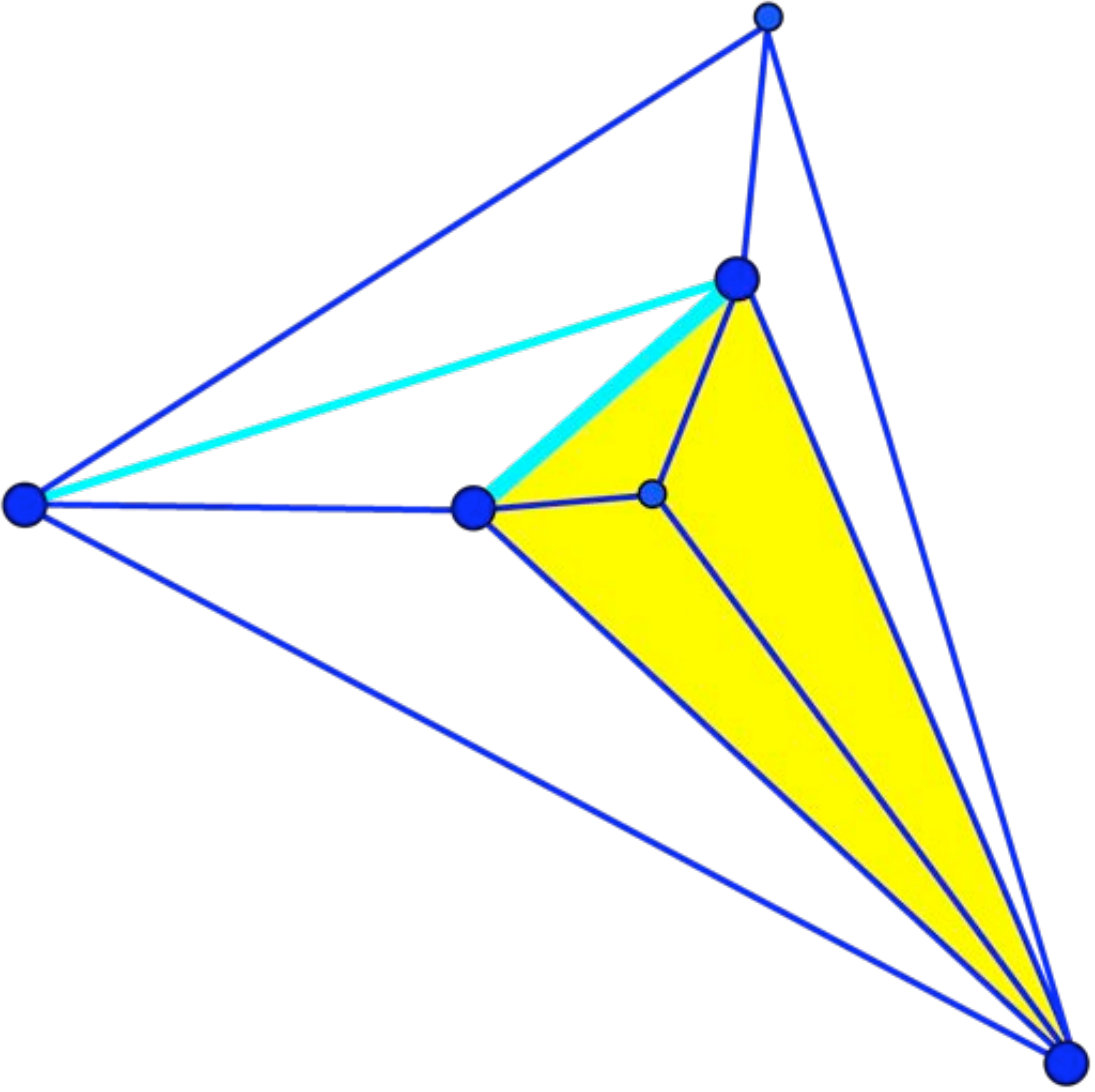, scale=.035}}\quad
  \subfigure[]{ \psfig{file=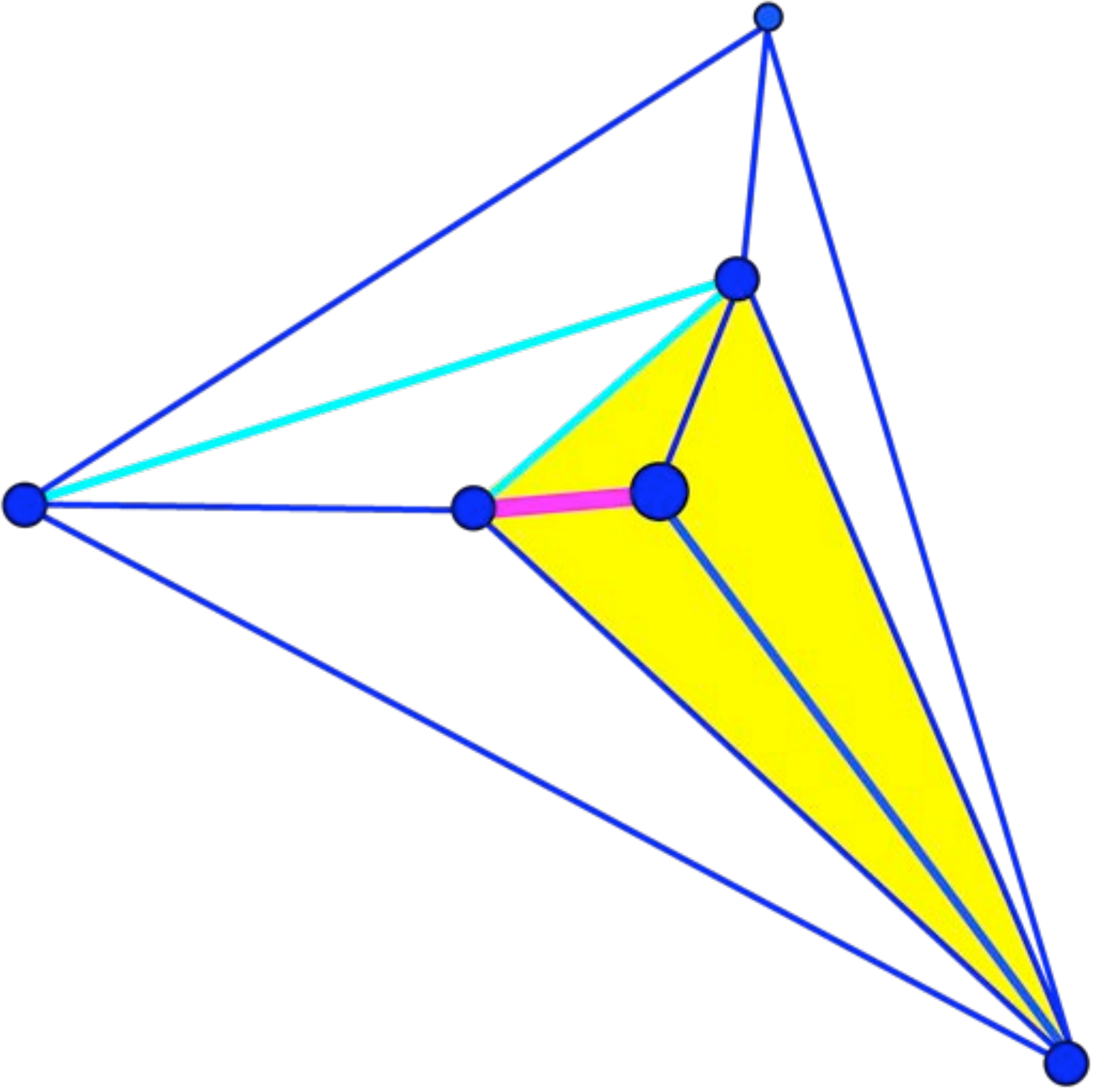, scale=.035}}
\caption{Given an edge joining two interior vertices,(a), either the edge identifies a short 
(non-facial) triangle (b), with a new interior edge, or, in the end, we have identified
a long interior edge (c).\label{Interior}}
\end{center}
\end{figure}
\begin{lemma} [{\bf Step 1}]  Let $D$ be a triangulated disc, and 
 $E$ be an edge joining two vertices interior to $D$.
Then there is an edge in the interior of $D$ that is long.
\label{shrink_interior_edge_lemma}
\end{lemma}

\proof If $E$ is long, we are done.
Suppose $E$ is short and therefore an edge of a non-facial triangle, $\T$ which separates the disc. 

Therefore, $\T$ is a separating non-facial triangle in $D$, with edges and vertices in its interior.  
Also, the boundary triangle $T$ contains at most one vertex, say $v$, of the boundary polygon of $D$.  
Therefore, the interior of $\T$ does not contain any vertices of the boundary polygon of $D$.  
  
We prove, by induction on the number of vertices interior to $\T$, that there is a long edge inside $\T$.
If there is only one vertex interior to $\T$, then that vertex is 3-valent.  We may contract any one
of the three edges. 

Next, assume that whenever there are $n$ or fewer vertices in $\T$, there is an interior edge we can contract, and suppose there are $n+1$
vertices in the interior of the triangle $\T$.  Then there is an interior edge inside of $\T$.  If that edge is long, we are
done.  Otherwise, it is short and we have a smaller non-facial triangle $\T'$ which contains fewer than $n+1$
vertices, since it is inside $\T$.  By the induction hypothesis, 
there is a long edge in the interior of $\T'$ to contract.  \qed 


\subsection{Contracting the  Hub Vertices  within a disc}

\begin{figure}[h]
 \begin{center}
   \subfigure[] { \psfig{file=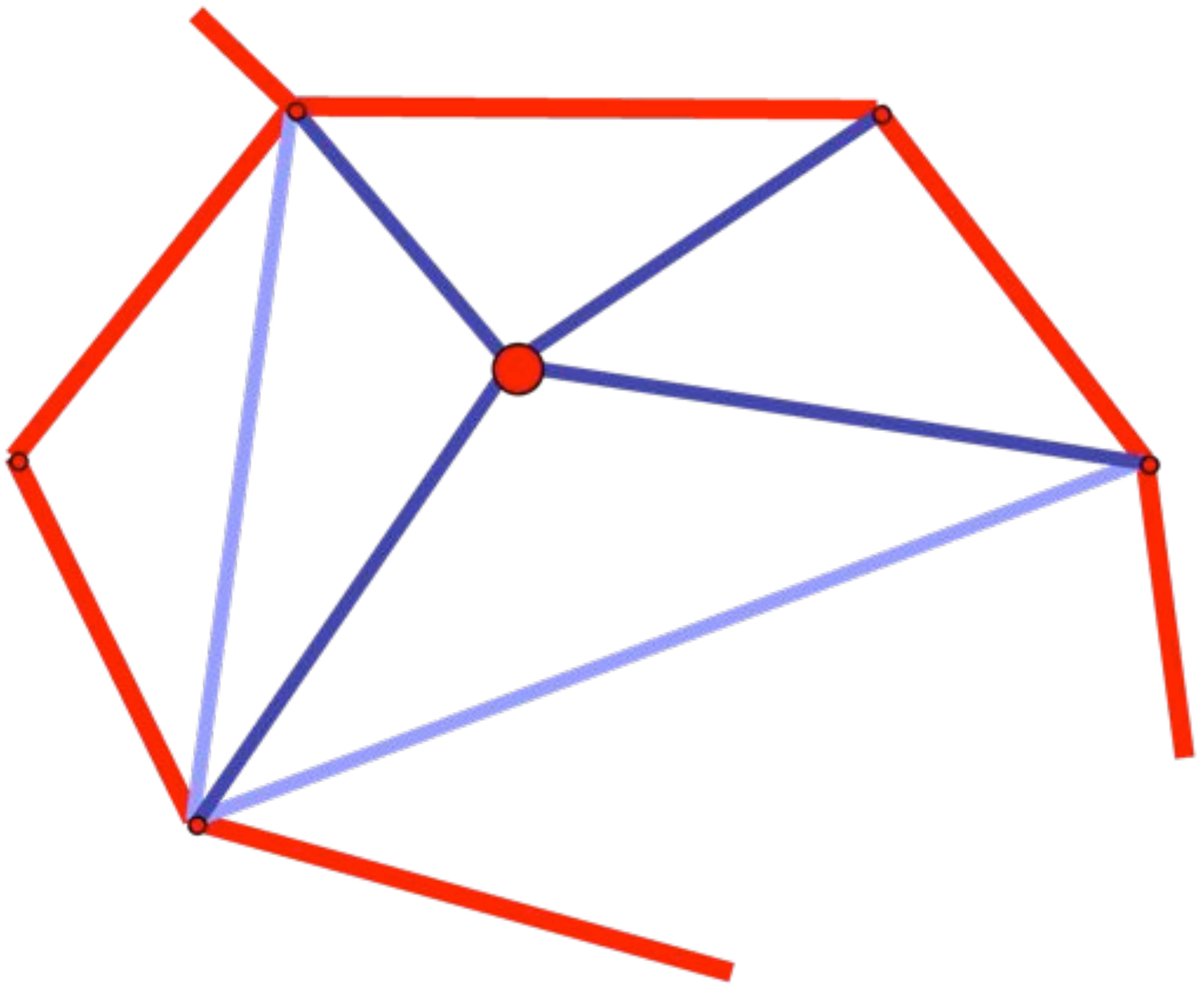, scale=.035}}\quad\quad
  \subfigure[] { \psfig{file=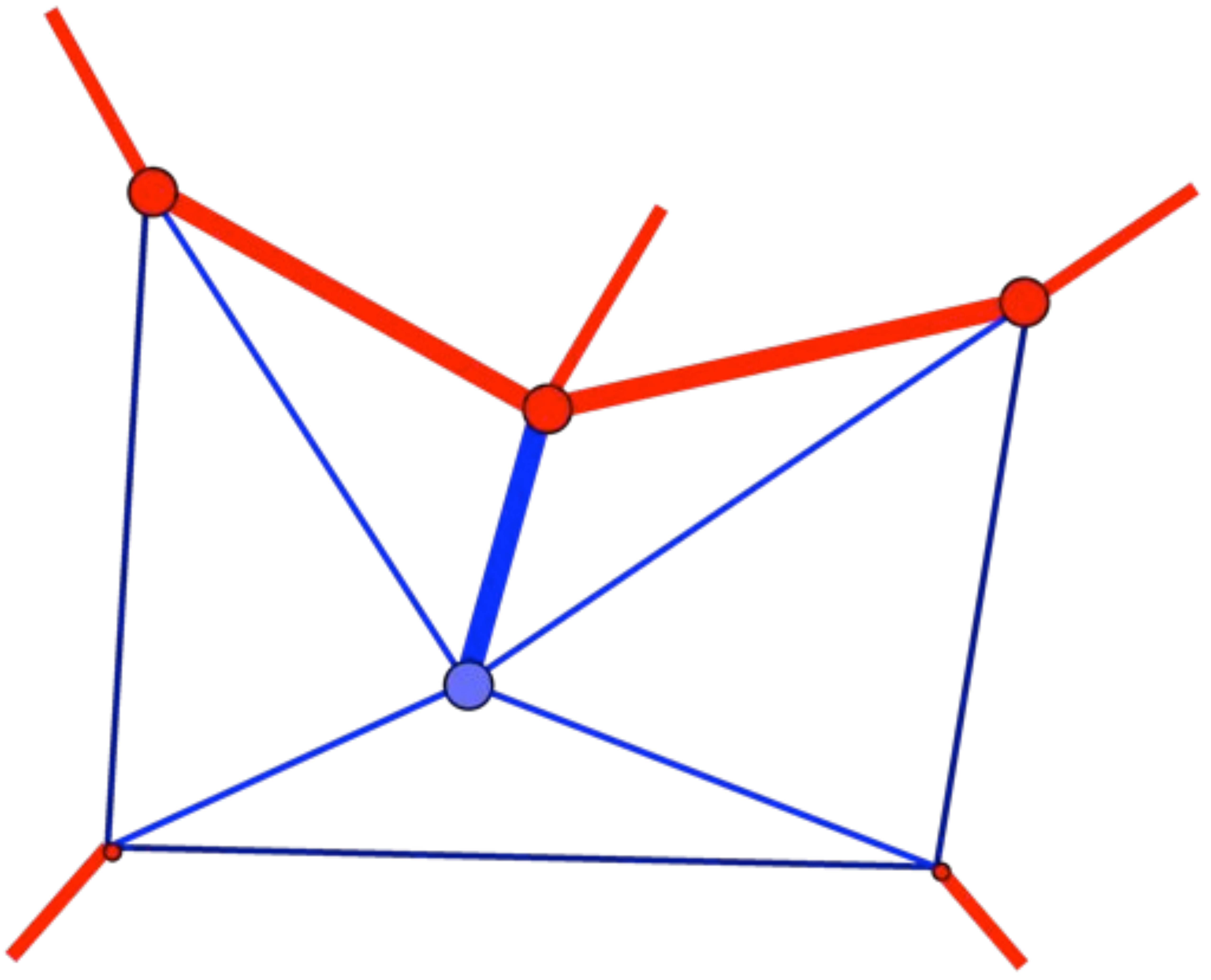, scale=.035}}\quad
  \subfigure[]{ \psfig{file=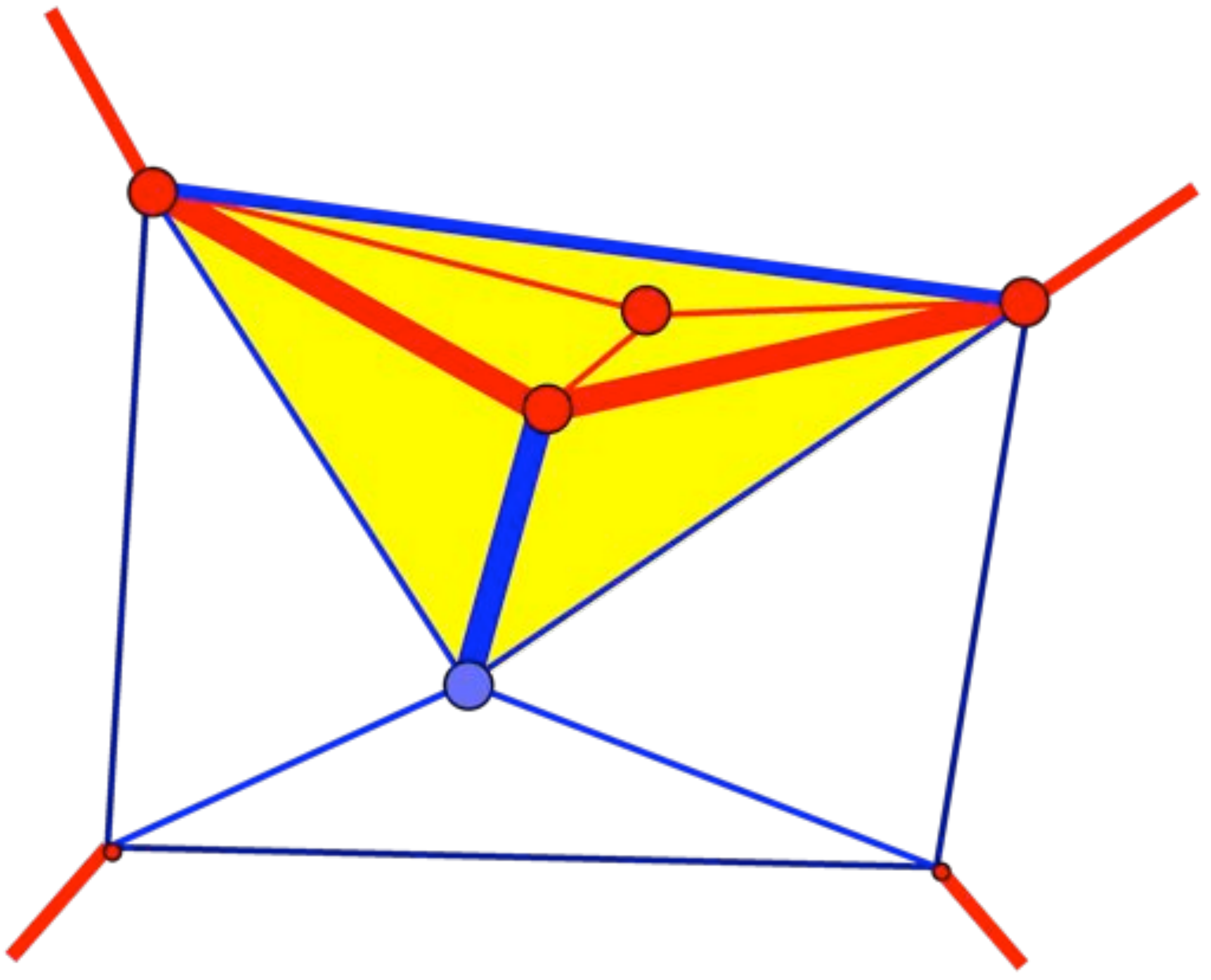, scale=.035}}
 \caption{A wheel with a central hub (a).  Given two non-adjacent rim vertices on different paths (b), there is no edge 
 connecting them provided the boundary of the disc is clean. (c) gives a boundary which is not clean.}  \label{fig:unclean}
 \end{center} 
 \end{figure}

As an immediate corollary of Lemma \ref{shrink_interior_edge_lemma},
we can contract all the edges between interior vertices in a
triangulated disc until we are left with no edges between interior vertices, thus reducing the number of vertices in the interior of the 
triangulated disc.  The only remaining interior vertices will be adjacent only to vertices on the boundary of the disc - which we call a `wheel' Figure~\ref{fig:unclean}(a).  

A {\it wheel} is a graph with vertex set $V = \{h, r_1, r_2, \ldots, r_n\}$.  Setting $r_{n+1} = r_1$, the set of edges
of the wheel is the set $E = \{ \{r_i, r_{i+1}\}|i = 1, 2, \ldots, n\} \cup \{ \{r_i, h\}|i = 1, \ldots, n\} $.  We call
the vertex $h$ the {\it hub} of the wheel, and the vertices $r_i, i = 1, 2, \ldots, n$ the {\it rim vertices} of the
wheel.  The edges of the form $\{ \{r_i, r_{i+1}\}|i = 1, 2, \ldots, n\}$ are the {\it rim edges} of the wheel, and the
edges of the form $\{ \{r_i, h\}|i = 1, \ldots, n\}$ are called {\it spokes}. 

From the previous discussion, the remaining interior vertices of a disc $D$ will be hub vertices of wheels, with all their rim vertices on the boundary polygon of $D$.  
We now consider conditions under which we can shrink a spoke and reduce the number of hubs of wheels.  

If the disc has a {\it clean boundary}, that is, all edges between boundary vertices are interior to the triangulation of the disc, 
then we may shrink hubs to the rim.  See Figure~\ref{fig:unclean}(b)(c).
It is possible to assume less than this strong property of the boundary, with a more delicate set of assumptions, but that is not needed for our proofs here. 

\begin{lemma}  Let  $D$ be  triangulated disc with a clean boundary and no edges joining interior vertices, and $h$ be a hub vertex in the triangulated disc.   
Then each spoke to a rim vertex $r_i$ is long. 
\label{spoke contraction}
\end{lemma}

\proof Take a spoke $h,r_{i} $.  Assume this is short.  
Since there are no other interior vertices adjacent to $h$, any non-facial triangle would be of the form $h,r_{i},r_{j} $.  
This would require an edge $r_{i},r_{j}$.  
Since the boundary polygon is clean, this edge must be on the boundary or interior to the disc.   
So the non-facial triangle must contain interior vertices of the disc; i.e. the hub of some other wheel in the disc!   
However, given a triangulation of $D$ this would require some edge between two hub vertices, contradicting the assumption that there were no such interior edges.  
\qed

Note that shrinking some spoke may generate a new edge between vertices in the boundary, but these edges are also interior to the triangulated disc; 
so the boundary remains clean, and we can continue to shrink spokes and reduce the number of vertices interior to the disc.

\subsection{Contracting Path Edges}

The third step in the contraction sequence is to shrink the path along which two triangulated discs have been connected down to length 1. Later, we will consider  shrinking such paths to length $0$.   

We will consider two discs $D_{1},D_{2}$, each with a clean boundary.  
These discs are {\it well-connected} if there is a path $P=x_{1}...x_{k}$ which is the complete intersection of the two discs, 
and all vertices of the path, except the end points, are adjacent only to other vertices in the two discs.  

\begin{figure}[h]
 \begin{center}
   \subfigure[] { \psfig{file=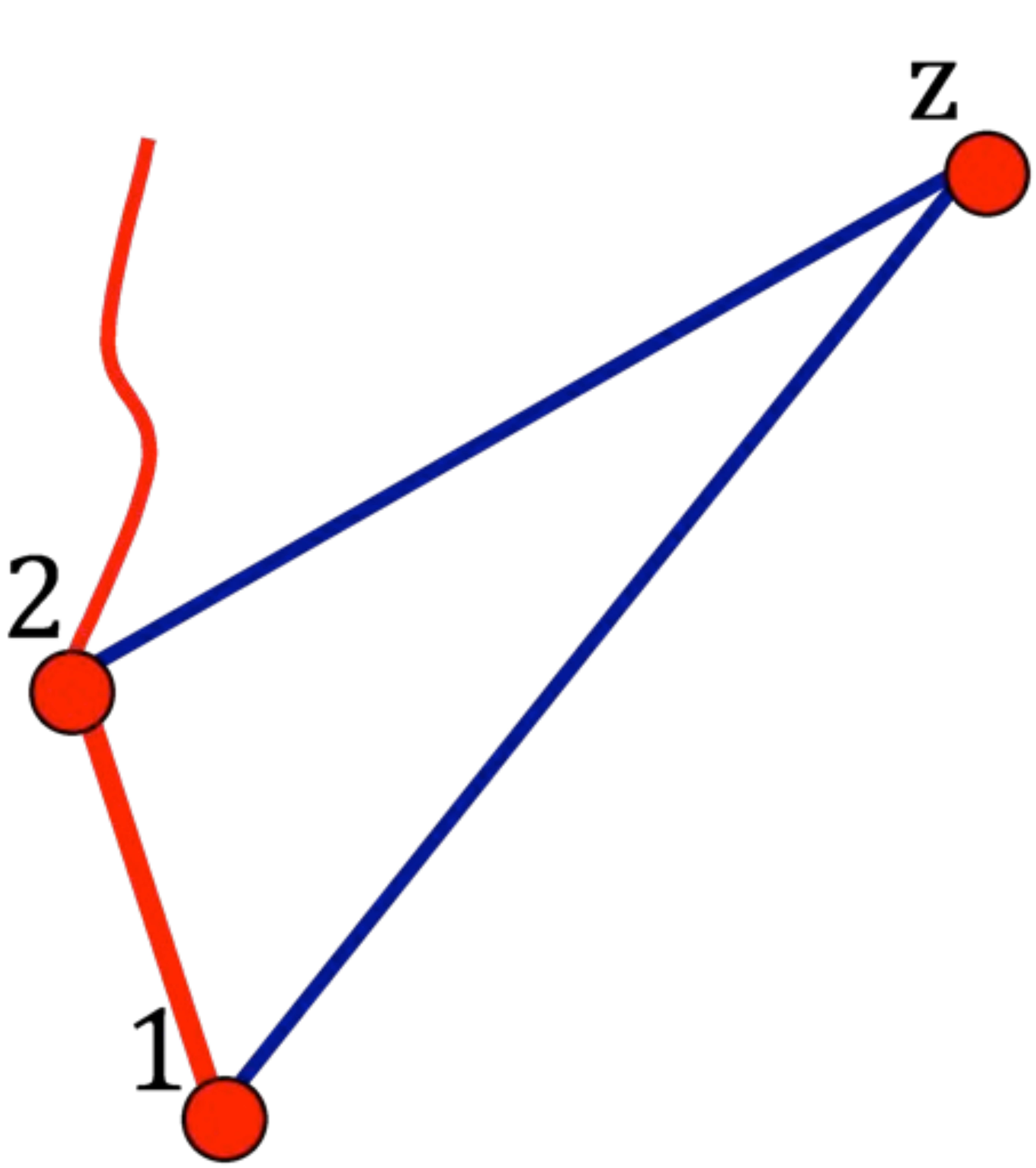, scale=.05}}\quad\quad
  \subfigure[] { \psfig{file=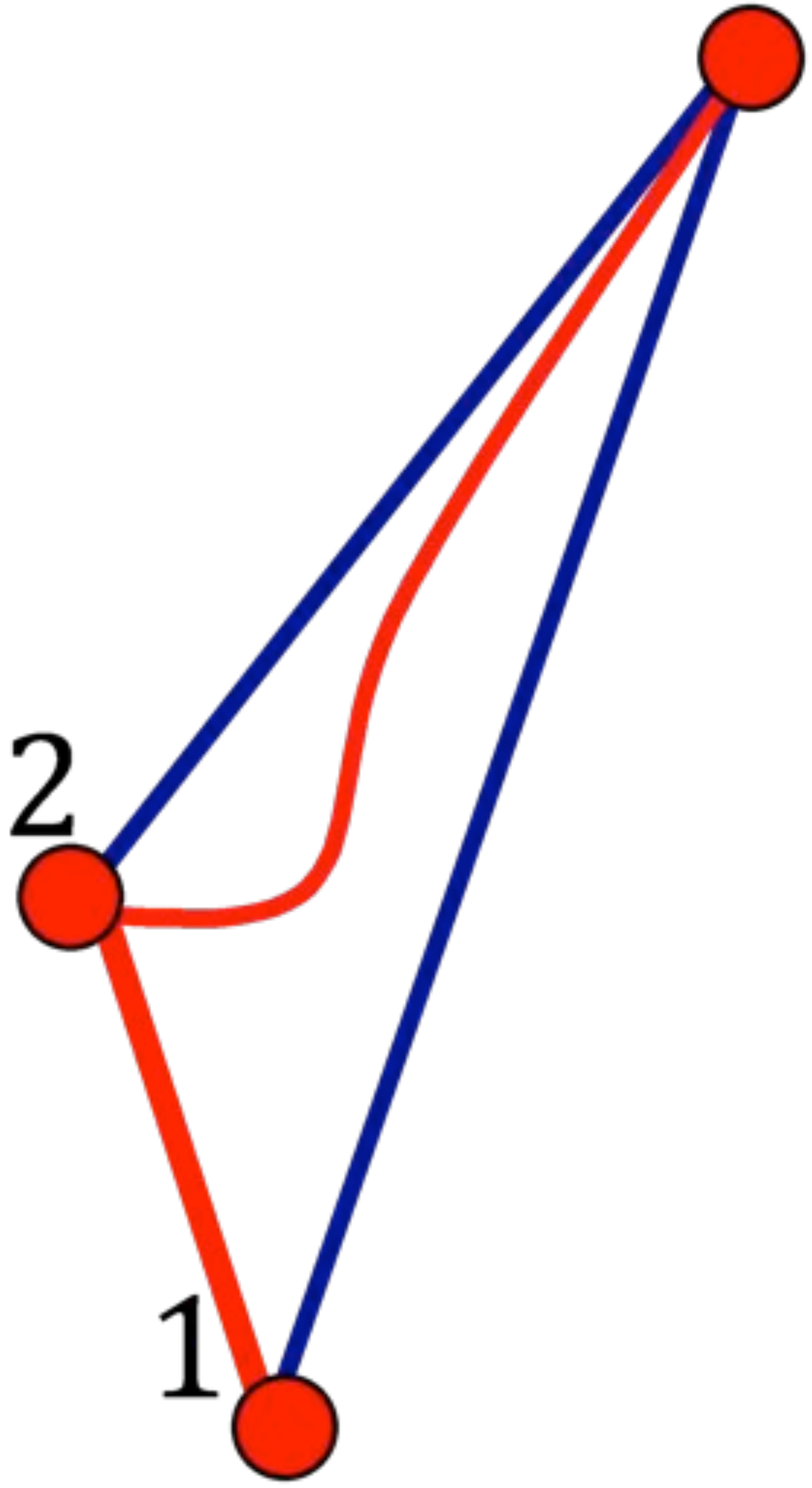, scale=.05}}\quad
 \caption{On a path edge of well-connected discs $D_{1},D_{2}$, a short triangle leads to a contradiction(b). }  \label{fig:pathcontract}
 \end{center} 
 \end{figure}
\begin{lemma} 
If $D_1$ and $D_2$ contain no interior vertices, have clean boundaries, and are well-connected along the path $P$, then each edge of $P$ is a long edge.  
\label{pathshrink}
\end{lemma}

\proof  Let  $P=x_{1}, \ldots, x_{k}$.  We claim that $\{x_1, x_2\}$ is long.  

If not, $\{x_1, x_2,z\}$ is in a non-facial triangle.  
If $z$ is not on the path ($z$ is only on one of the discs), then the edges $x_{1},z$ and $x_{2},z $ must be in the same disc.  
That means that the non-facial triangle is within a single disc - and therefore there are interior vertices in this disc - which is a contradiction.  

Otherwise,  $z$ must be within the path - with edges  $x_{1},z$ and $x_{2},z $.  
If these edges are in the same disc, then the argument above yields a contradiction.  
Alternatively, they are in distinct discs, and without loss of geneality, suppose  $x_{2},z$ is in $D_{2}$.  
That would mean that $x_{2},z $ joins two vertices of the boundary of $D_{1}$ but is not in the triangulation of $D_{1}$.  
This is also a contradiction of the boundary of $D_{1}$ being clean.  
 
We conclude that every edge interior to the path is long.  
\qed 

Shrinking such a long path edge does not add any interior vertices to either disc. Nor does it generate new edges outside the two discs, which could make either boundary not clean.  
Therefore, Lemma~\ref{pathshrink} will permit us to shrink all paths  to length $1$.  
Note that shrinking a path down to length $1$ could not change the fact that the boundary of any disc (including these discs) is clean. It will also not change the topology of the underlying polyhedron,   

\subsection{Cleaning the boundary of a triangulated disc. }
If we are given a triangulated disc, surrounded by other discs (some triangulated, some not) then the boundary may not be `clean'.  
Under what conditions can we clean the boundary, or guarantee the boundary is clean?   

The first step is to ensure a {\it clear path} separating two triangulated discs.  We give a negative definition: if two vertices on the shared boundary of two 
discs are adjacent within one disc (and not adjacent on the path), then the boundary path of the other disc is {\it not clear}.   
We will reroute the boundary between the discs to make it shorter and clear for both discs (Lemma~\ref{pathclearing}).   
We call these local moves, involving just two discs and their shared boundary {\it clearing the path}.  

In order not to swallow up a disc in the following process, we call two discs {\it well-attached} if they are well-connected and, 
in addition, there are boundary edges of each disc which are not part of the path (avoiding Figure~\ref{fig:unclean}(c)).  

We consider two discs, $D_{1}$ and $D_{2}$ which are well-connected along   
$P=x_{1}...x_{k}$, $k>1$, and assume $\{x_{1},x_{k}\}$ is not an edge of the larger graph. 
Then $D_{1}$ and $D_{2}$ are {\it clearly attached} if there is no 
edge in the graph between pairs of vertices $x_{i},x_{j}$ which are not adjacent on the path. Otherwise it is not clearly attached. 

\begin{lemma} 
If $D_1$ and $D_2$ are well-attached along the path $P=x_{1}...x_{k}$ with $k>2$ but not clearly attached, 
then we can reroute the path through the combined triangulation on the disc, using only vertices of the path, to 
$P'=x_{1}\ldots x_{i'} \ldots x_{k}$, with the same endpoints, so that the modified discs remain well-attached and the attachment is   
clean.
\label{pathclearing}
\end{lemma}

\proof  Assume there is an edge $\{x_{i},x_{j}\}$ in the triangulation of one disc, say $D_{2}$, between path vertices $x_i$ and $x_j$ which are not adjacent on the path.  
Then the subpath $P^{i,j}$ of $P$ joining $x_{i}$ to $,x_{j}$, along with the edge $x_{i},x_{j}$ will, topologically, surround a smaller disc 
within $D_{2}$.  We shift the edges, vertices (including those in $P^{i,j}$), and triangles into $D_{1}$, and reroute 
the boundary between the discs to $P'=x_{1}\ldots x_{i},x_{j}\ldots x_{k}$.

This does not add edges outside of the discs, and keeps the discs well-attached with a shorter boundary path.  However, there will now be additional interior edges within $D_{1}$. 
Repeated application of this rerouting will end up with a clear path (perhaps of length $1$).  
\qed

Note that if we do only insist that the discs are well-connected, it is possible that the rerouting encloses all vertices, edges and faces of $D_{2}$ and we lose the disc entirely, creating a new polyhedron which is not topologically equivalent.  

If we repeat this clearing of the paths for all the boundaries of a disc $D_{1}$, we are close to having a clean boundary.  There is one other possible 
difficulty: edges which run from one vertex on the boundary, across to another vertex on the boundary, but not through a well-attached disc.  
For example, if a third disc $D_{3}$ had a boundary intersection which was not a single path, then we could still have such problematic edges.  See Figure~\ref{fig:unclean}(c) where the intersection of two discs is a pair of vertices with no shared path. 
Also, if there are edges in the larger graph which are not restricted to specific discs, the same problem could arise.

In the applications we intend, there is a set of discs (not necessarily all triangulated) which will surround the triangulated discs, with a few added 
edges which are guaranteed not to join vertices on the boundary of a single disc.   This will be enough to permit us to clean up the boundary 
of a triangulated disc, and apply the full effect of the previous lemmas; shrinking back all vertices within triangulated discs and shrinking 
all paths between discs to length $1$.  

A triangulated disc $D$ is {\it well-surrounded} by a cycle of discs $C=D_{1} \ldots D_{m}$ if:
\begin{enumerate}
\item [(i)] each disc $D_{i}$ of the cycle intersects $D$ in either one vertex or a connected path of edges;
\item  [(ii)]each disc $D_{i}$ intersecting $D$ in more than one edge is also a well-attached triangulated disc;
\item  [(iii)]all edges between vertices of the boundary of $D$ are in one of these discs.  
\end{enumerate}

\begin{lemma} 
If $D$ is well-surrounded by a cycle of discs $C=D_{1} \ldots D_{m}$, then we can clear paths shared by $D$ and triangulated discs 
$D_{i}$ so that $D^{*}$ is well-surrounded and has a clean boundary.  
\label{boundarycleaning}
\end{lemma}

\proof  Along the boundary, we can clear each path which is shared with a well-attached disc, by Lemma~\ref{pathclearing}.   
	By property (iii) any edge joining two vertices of the boundary of $D$ has to run through one of these surrounding discs.  
But this is either within an adjacent triangulated disc (contradicting the step of clearing the shared path) or involves a disc 
meeting the boundary at vertices not along a path, contradicting (i). 
\qed

\section{The Main Theorem}\label{main}

We are now ready to prove the first main result - that an expanded polyhedron $\PEhat$ expanded from an isostatic block and hole polyhedron $\Phat$ is also isostatic.  The key tool is the graph theoretic inverse of the edge contractions of  Section 3, called vertex splitting. 

\subsection{Vertex Splitting} \label{sec:vertexsplit}

We recall a graph theory technique called vertex splitting, which will later be used to add a vertex to an isostatic structure 
in such a way as to preserve it's infinitesimal rigidity and independence \cite{vertexsplit}.   

\begin{figure}[h]
 \begin{center}
  \subfigure[] { \psfig{file=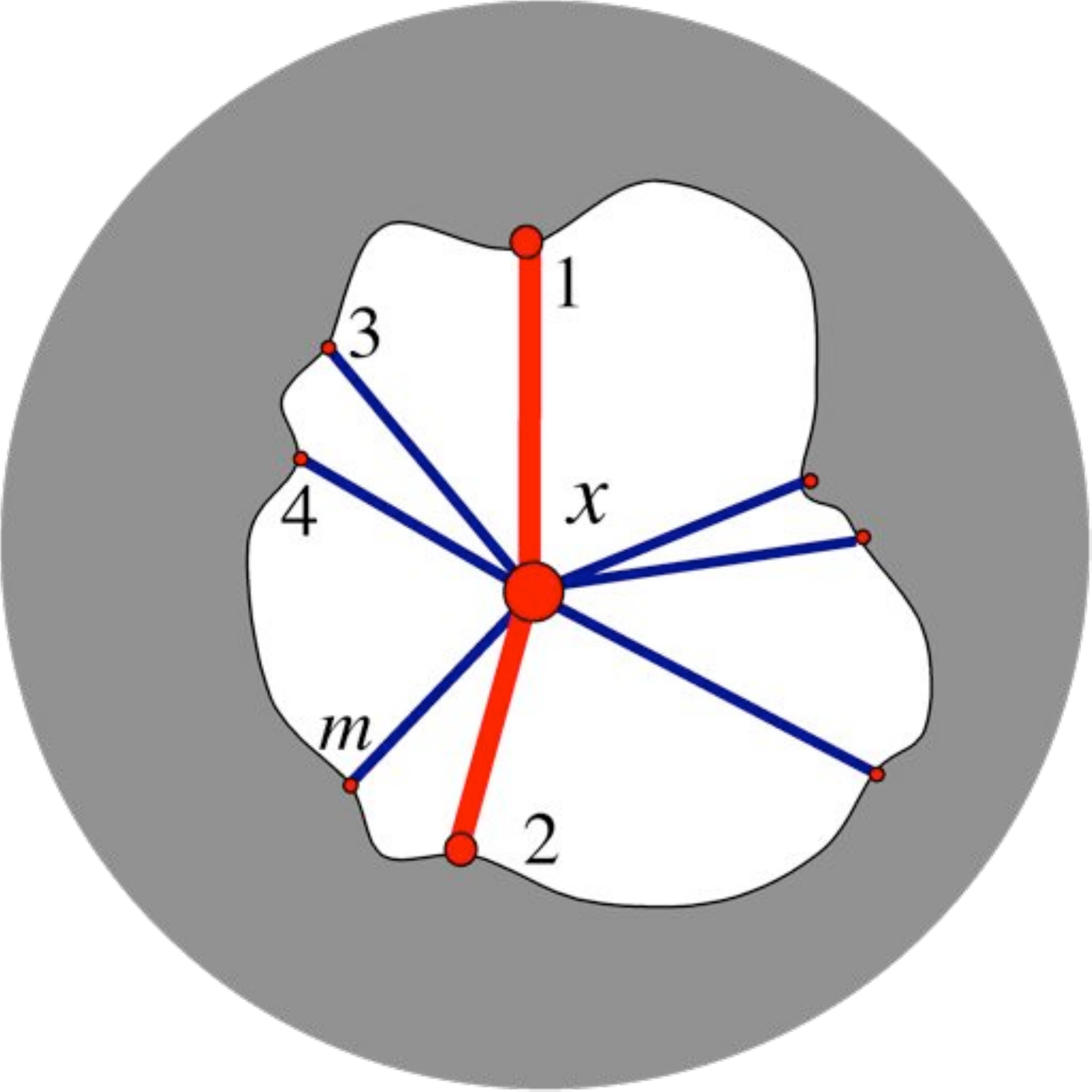, scale=.06}}\quad\quad\quad
   \subfigure[]{ \psfig{file=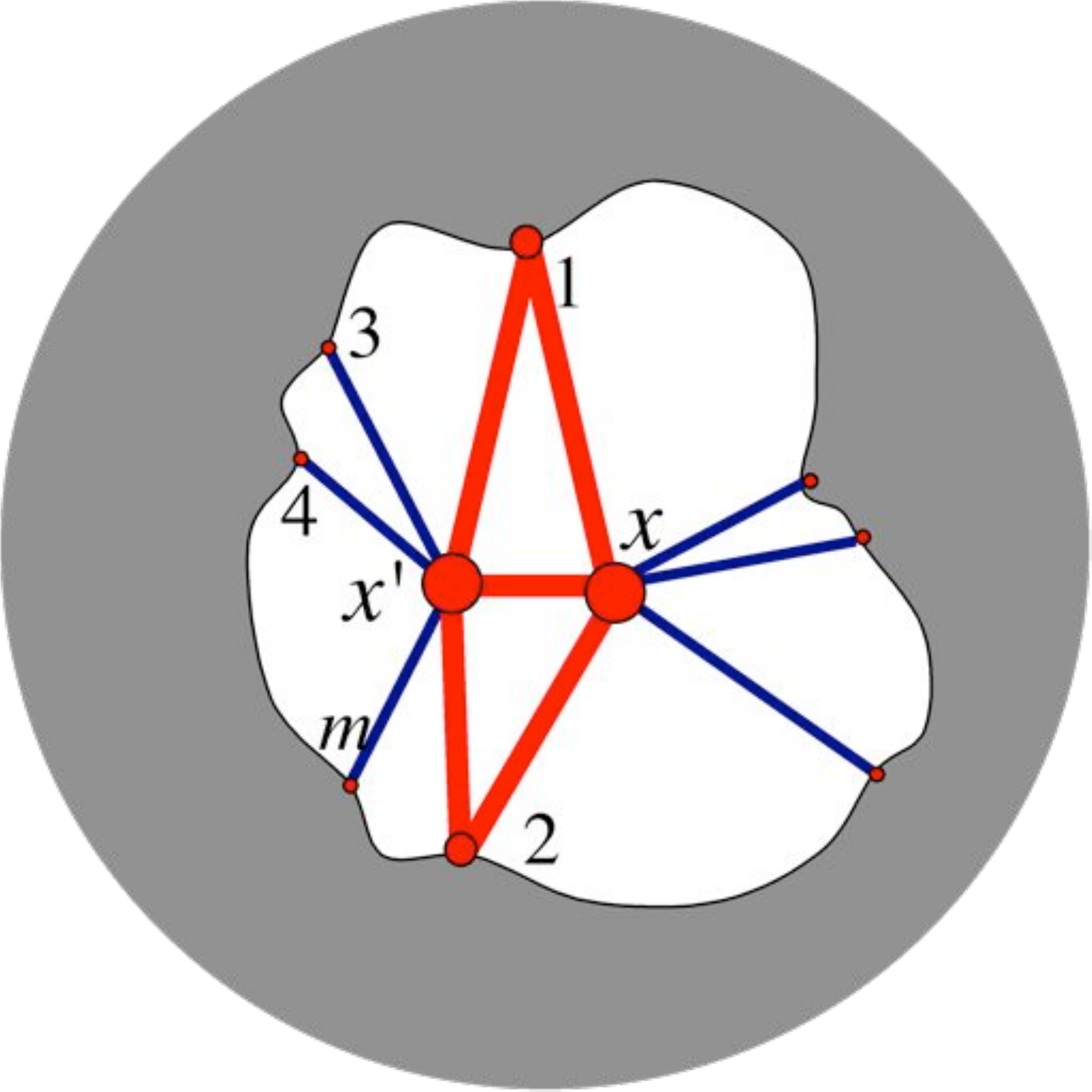, scale=.06}}
 \caption{A vertex split applies locally to a vertex and two attached edges (a).  This adds one vertex and a net of $3$ new edges.}
\label{fig:vertexsplit}
 \end{center} 
 \end{figure}
\begin{definition}
{\rm   Given a graph, $G = (V, E)$ with a vertex $x$ and 
subsets of edges $E_1 = \{(x, 1), (x, 2)\}$, and $E_2 = \{ (x, 3), (x, 4), \ldots, (x, m)\}$, a} vertex split 
{\rm of $x$ on the edges in $E_1$ is the modified graph $G' = (V', E')$, with new vertex $x'$, so $V' = V \cup \{x'\}$, and
an edge set $E' = (E \backslash E_2)  \cup \{(x, x'), (x', 1), (x', 2), \ldots, (x', m)\}$ .}
\end{definition}
We note that we are adding one new vertex and three new edges, preserving the count for being isostatic: $|E| = 3|V| - 6$.  The following result is central to our use of vertex splitting: 

\begin{theorem}[Vertex Splitting: Whiteley\cite{vertexsplit}]
Given a generically isostatic (independent, infinitesimally rigid) graph $G$ in $3$-space then the new graph $G'$ formed by a vertex split on two edges is generically isostatic (independent, infinitesimally rigid). \label{thm:vertex-split}
\end{theorem}

\begin{figure}[h]
 \begin{center}
  \subfigure[] { \psfig{file=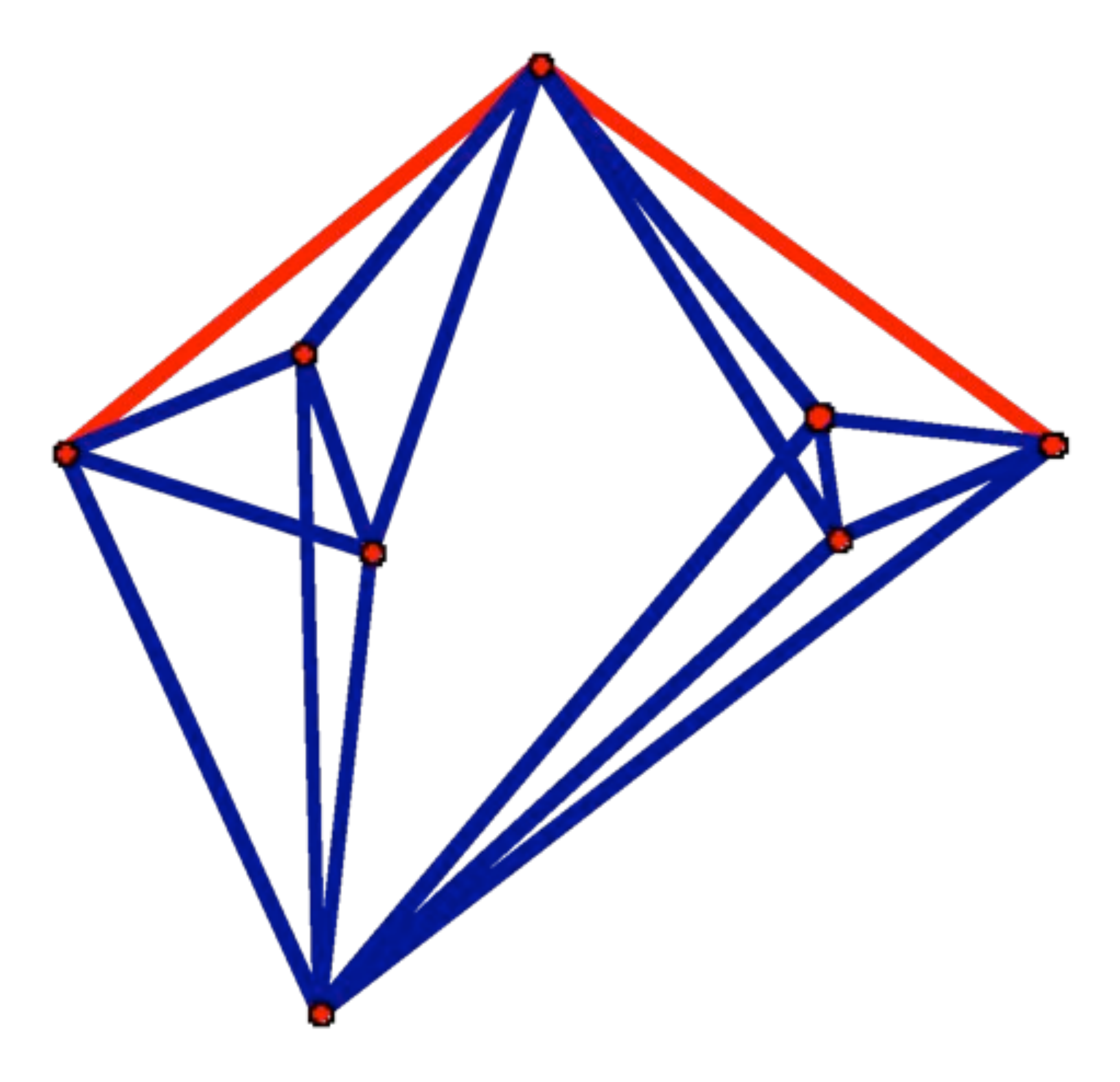, scale=.13}}\quad\quad\quad
  \subfigure[]{ \psfig{file=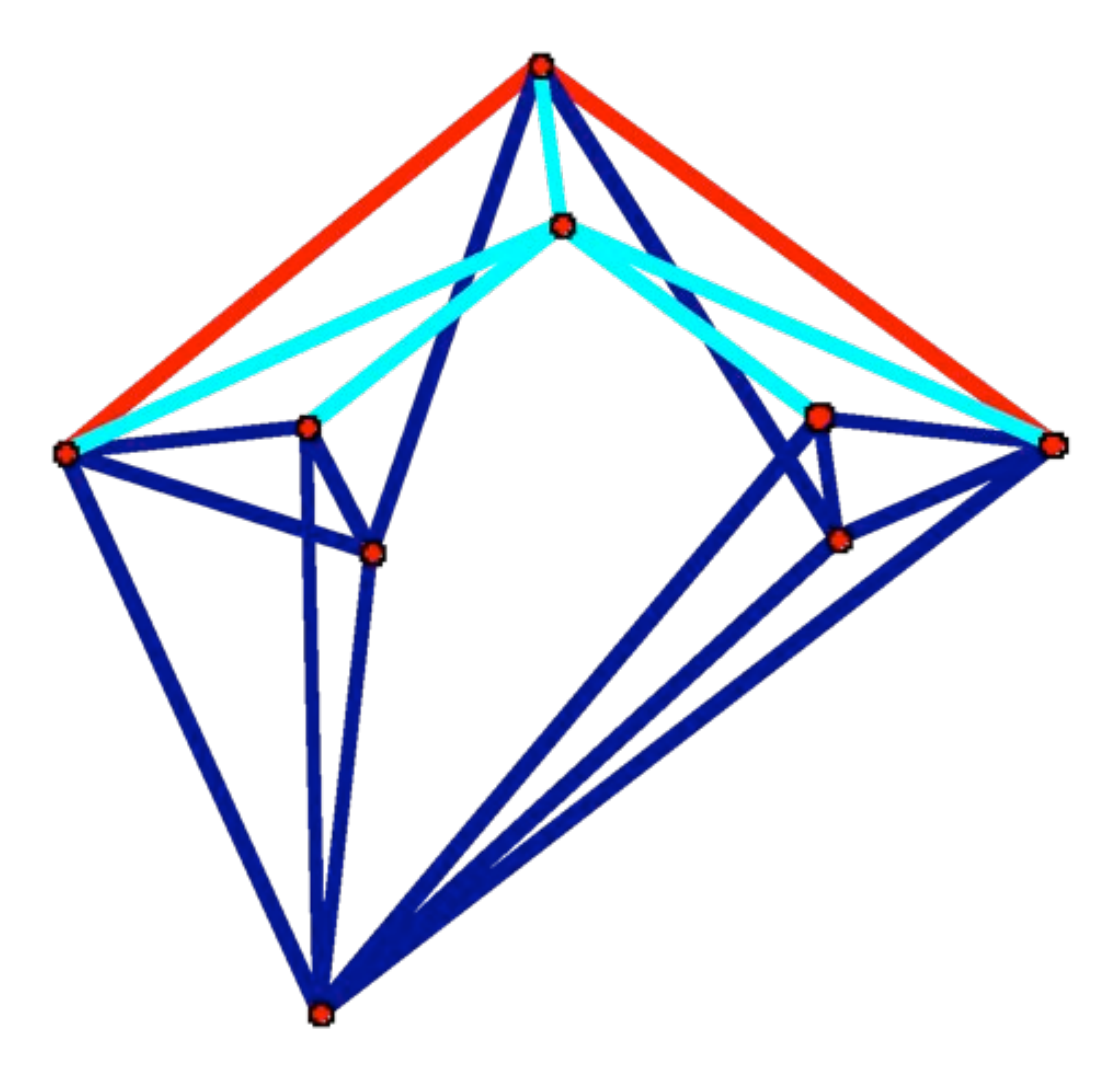, scale=.13}}
 \caption{The double banana (a) is generically flexible, but a vertex split takes it to a generically isostatic framework (b). }
\label{fig:bananasplit}
 \end{center} 
 \end{figure}

The converse is not true; a vertex split on two edges on a non-rigid framework may yield a rigid framework.  
For example, starting with the double banana in Figure \ref{fig:bananasplit}a, we can vertex split on the red vertex and two red edges to obtain the rigid framework in 
Figure \ref{fig:bananasplit}(b).
 
\begin{lemma}
A surface vertex splitting in a block and hole polyhedron $\Phat$ is the graph theory inverse operation to edge contraction on an edge in exactly two triangles.   

Given a polyhedral block and hole polyhedron $\Phat$, vertex splitting on vertices bounding or in surface discs, with the adjacent edges separated following the topology of the edge cycle at $x$ is the topological inverse of edge contraction on an edge in a surface disc or between two surface discs. 
 
 Such a vertex split generates topologically equivalent polyhedron $\Phat_{1}$
\end{lemma}
\proof
The proof follows from the definitions of vertex splitting and edge contraction.
\qed

We recall that short edges will not be contracted in this paper. If an edge is short, it is the edge of three triangles (including a non-facial triangle).  When such an edge is contracted in a generically isostatic graph, 
we lose one vertex and four edges so the graph becomes undercounted ($|E'| = 3|V'| - 7$) and therefore the graph becomes generically flexible.  


As a simple application of Lemma \ref{shrink_interior_edge_lemma} and Theorem~\ref{thm:vertex-split}, we verify Gluck's result that all triangulated spheres in 3-space are generically rigid, by showing that they are constructed from the double triangle by vertex splitting (see  \cite{vertexsplit} for a related proof).

\begin{corollary}[Gluck \cite{gluck}] All triangulated spheres are generically isostatic in 
$3$-space. 
\end{corollary}

\proof  We can separate the triangulated sphere $\poly$ into two surface discs
by just taking one triangle $\T$ as the first disc, and the rest as the second disc $D$. 
The proof is then by induction on the number $m$
of vertices of second disc $D$. If the triangulated
sphere has shrunk to a double triangle $\T,\T$, we are finished. 
If the number of vertices of the disc $D$ is greater
than 3, by Lemma~\ref{shrink_interior_edge_lemma} we can find a contractable edge. 
Since the double triangle is generically rigid in 3-space,
and the triangulated sphere can be obtained by vertex-splitting, the reverse of edge contraction, and hence the
triangulated sphere $\poly$ is rigid.  \qed 

\subsection{The Contraction Sequence}\label{sec:contraction}
In this section, we will assume that we are given a block and hole polyehdron $\Phat$  that is sub-divided into  blocks,  holes and  triangulated discs.  We also require that 
\begin{enumerate}
\item[1)]   each pair of discs is well-attached; 
\item[2)]   each triangulated disc is well-surrounded;  and 
\item[3)]   all vertices are in some disc.  
\end{enumerate} 
Since these are properties to be preserved throughout our basic contraction sequence, we name these as {\it  well-designed} block and hole polyhedra.  

In Algorithm~\ref{contractionSequence} below, we describe the sequence of steps. 
We verify that when each of the steps is applied to change the polyhedron, either 
\begin{itemize}
\item[(i)] the step reduces the number of vertices of the boundary paths, leaving the total number of vertices unchanges; or 
\item[(ii)] the step reduces the the number of vertices overall with no increase in the number of vertices on any boundary paths.  
\end{itemize} 
In Section 4.3, we confirm that this algorithm inductively contracts our expanded polyhedron down to our base polyhedra in steps reversed by vertex splitting.  

\begin{algorithm}[Contraction Sequence]
{\bf do}\\
\indent {\bf Step 1}: Clean the disc boundary of all triangulated discs.\\
\indent {\bf Step 2}: Contract all interior edges.\\
\indent {\bf Step 3}: Contract all contractible spokes.\\
\noindent {\bf while} (there is a disc that is not clean)\\
\noindent {\bf Step 4}: Contract path edges to length 1.
\label{contractionSequence}
\end{algorithm}

The first step in the Contraction Sequence is to clean the boundary paths of length greater than one, 
thereby reducing the lengths of the boundary paths between triangulated discs.
We note that if a triangulated disc is adjacent to a block or a hole, 
then the shared boundary path has length one and is already clear. \\

\noindent {\bf Step 1}: Clean the disc boundary of all triangulated discs. 

\begin{proposition} If application of Step $1$ to a well-designed block and hole polyhedron $\Phat$  yields the polyhedron $\Phat_1$, then:
\begin{itemize}
\item[(i)] each disc of $\Phat_1$ has a clean boundary;
\item[(ii)] if $\Phat_1\neq \Phat$, then the length of each cleaned boundary path between two triangulated discs in $\Phat_1$ will be less than the length of the corresponding boundary path in $\Phat$,  other  boundary paths are unchanged and the total number of vertices is unchanged; and 
\item[(iii)] $\Phat_1$ is topologically equivalent to $\Phat$ and therefore $\Phat_2$ is a well-designed block and hole polyhedron. \end{itemize}
\label{cleaning}
\end{proposition}

\proof Lemma~\ref{boundarycleaning} yields \ref{cleaning}~(i) and \ref{cleaning}~(iii).
 
If the boundary between two discs in $\Phat$ is clear, the length of that boundary path will not change. On the other hand, if the path is not clear 
the process of clearing the path reroutes part of the path through an edge, thereby shortening the path (Lemma~\ref{pathclearing}), and the new path between the two discs will be shorter than the original path. 

These steps keep the same disc structure in the topology, including the properties of discs being well attached and well-surrounded. 
 \qed 

Notice that application of Step 1 does not change the underlying graph - just the subdivision into triangulated discs. As such, it does not alter the generic rigidity or generic independence of the associated polyhedron. 

After completion of Step 1 on any polyhedron with an unclean boundary, we are left with a polyhedron $\Phat_1$ whose boundary paths have fewer vertices and edges than the polyhedron $\Phat$.  The total number of vertices is not changed. 
By cleaning the discs we have not changed the topology of our disc structure so the polyhedron $\Phat_1$
remains sub-divided into discs that are well-designed, and we note that the boundary of each triangulated disc is clean.  
Our next step will, for every triangulated disc, remove all edges interior to the disc.  \\

\noindent {\bf Step 2}: Contract all interior edges. 

\begin{proposition} If application of Step $2$ to the  well-designed block and hole polyhedron $\Phat_{1}$ output from Step $1$ yields the polyhedron $\Phat_2$, then:
\begin{itemize}
\item[(i)] the only vertices interior to triangulated discs of $\Phat_2$ are the hub vertices of wheels;
\item[(ii)]  $\Phat_1$ can be obtained from $\Phat_2$ by a construction sequence of vertex splits;
\item[(iii)] if $\Phat_2 \neq \Phat_{1}$, then the number of vertices has decreased and the number of path vertices is unchanged; and 
\item[(iv)]  $\Phat_2$ is topologically equivalent to $\Phat_{1}$ and therefore $\Phat_2$ is a well-designed block and hole polyhedron.
\end{itemize}
\label{only_hubs}
\end{proposition}

\proof By Lemma \ref{shrink_interior_edge_lemma},
we can inductively contract all interior edges of the
discs.  This does not alter the underlying topology of the discs, so $\Phat_{2}$ is also  a well-designed block and hole polyhedron, giving (i), (iii) and (iv).  

Since vertex splitting is the reverse of edge contraction, the Proposition follows.  \qed

After completion of Step 2 on any interior edges, we are left with a polyhedron $\Phat_2$ which has fewer vertices and edges than polyhedron $\Phat_1$.  
Our next step will remove some (prehaps all) remaining vertices interior to a disc in the polyhedron $\Phat_2$.  

Recall that a spoke is an edge from a hub vertex to a path vertex.  We say a spoke is {\it contractible} if 
the boundary of the disc containing the spoke is clean. Our next step is:\\

\noindent {\bf Step 3}: Contract all contractible spokes.

\begin{proposition} [{\bf Step 3}]
If application of Step $3$ to the polyhedron $\Phat_2$ output from Step $2$ yields the polyhedron $\Phat_3$, then
\begin{itemize}
\item [(i)] some discs of $\Phat_3$ do not contain any interior vertices; 
\item [(ii)] $\Phat_3$ can be obtained from $\Phat_2$ by a construction sequence of vertex splits;
\item[(iii)] if $\Phat_3\neq \Phat_{2}$, then the number of vertices has decreased, and the number of path vertices is unchanged;
\item[(iv)] $\Phat_3$ is topologically equivalent to $\Phat_{2}$ and all pairs of faces well-attached.
\end{itemize}   
\label{wheelremoval}
\end{proposition}

\proof  
 Within any given disc of $\Phat_2$, each spoke adjacent to a hub is long by Lemma~\ref{spoke contraction}, and we may contract one spoke for each hub.   This operation conserves the underlying topology, and this still comes from a well-designed block and hole polyhedron.  
 
\begin{figure}[ht]
 \begin{center}
  \subfigure[] { \psfig{file=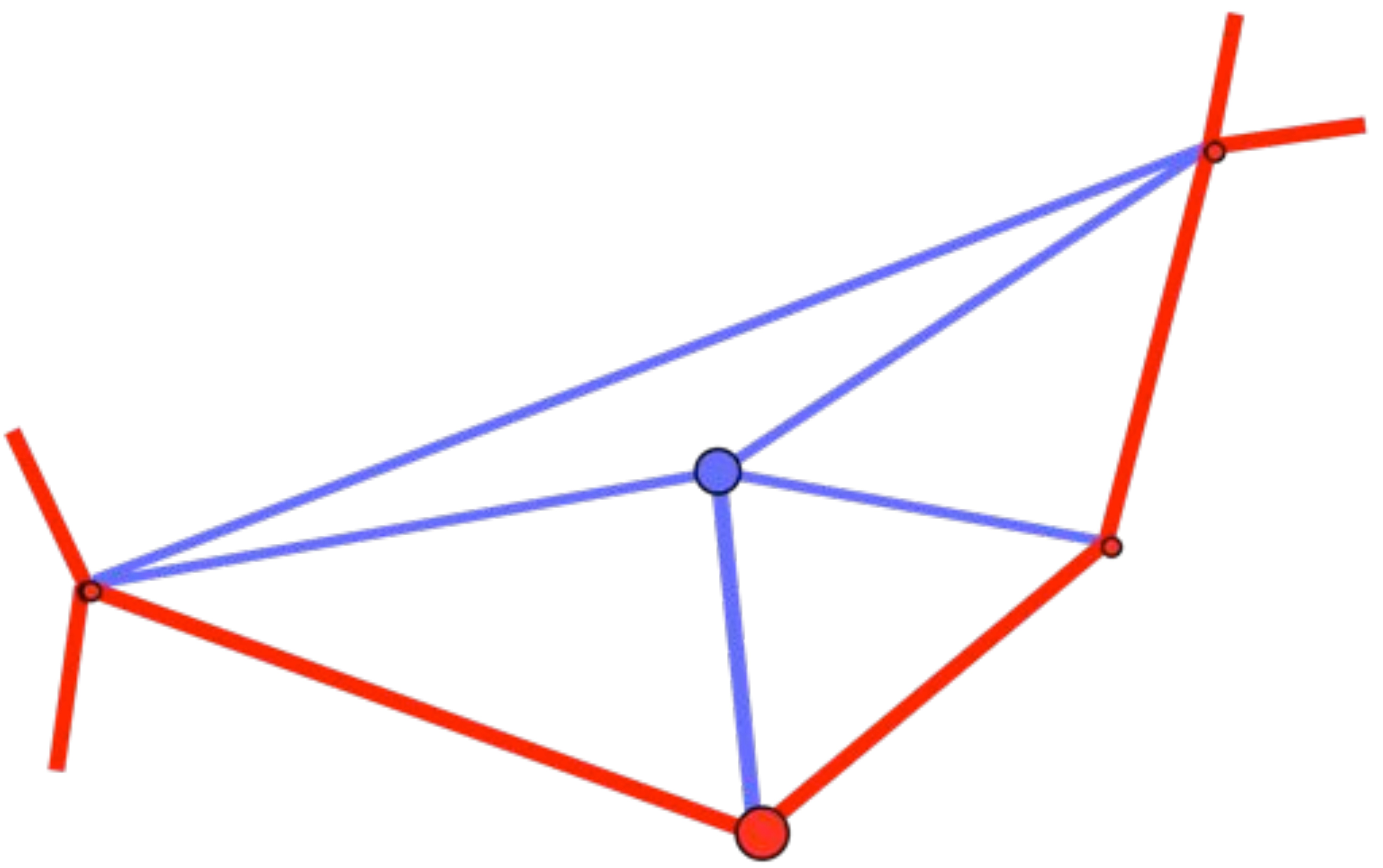, scale=.035}}\quad
  \subfigure[]{ \psfig{file=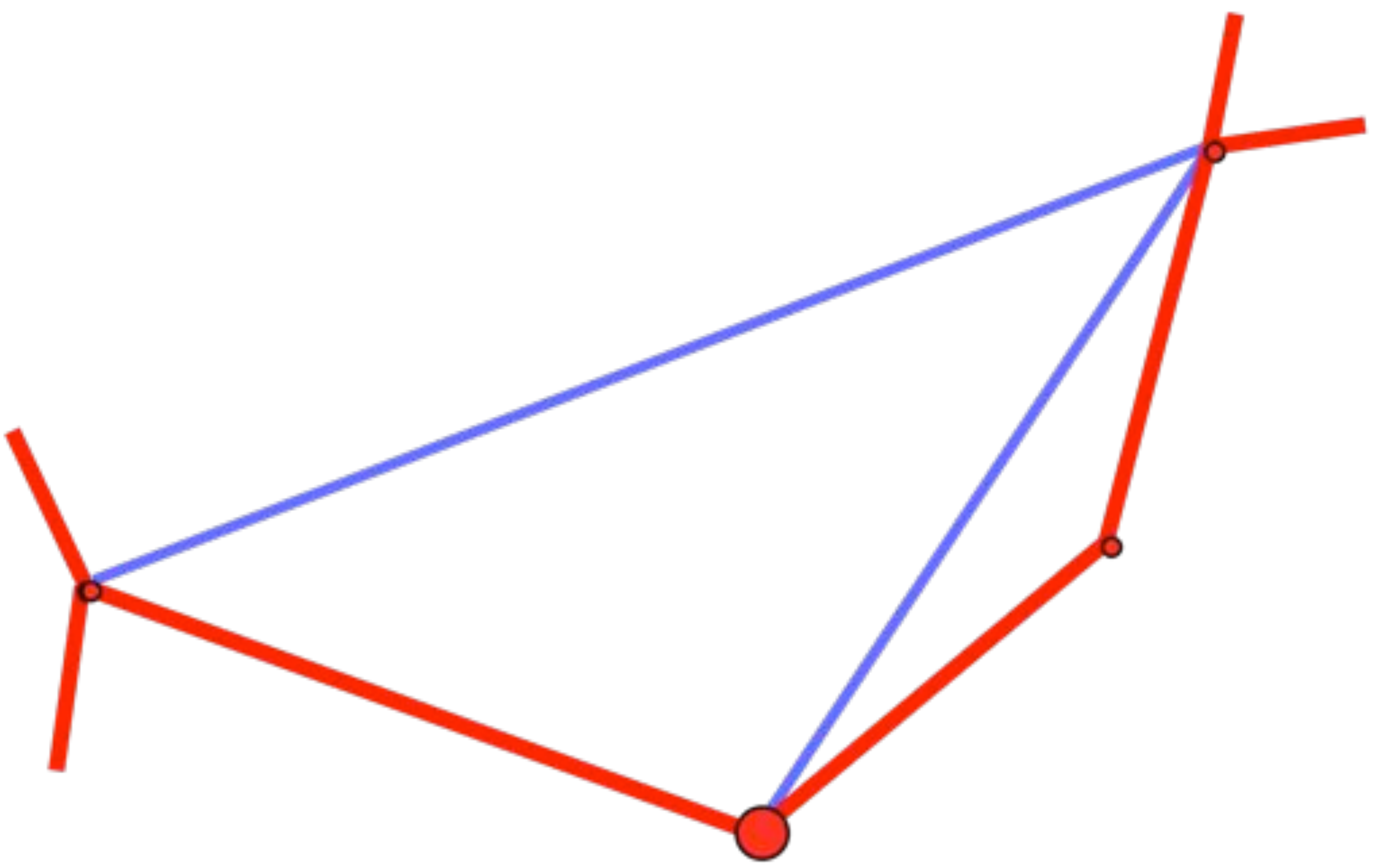, scale=.035}} \quad
  \subfigure[]{ \psfig{file=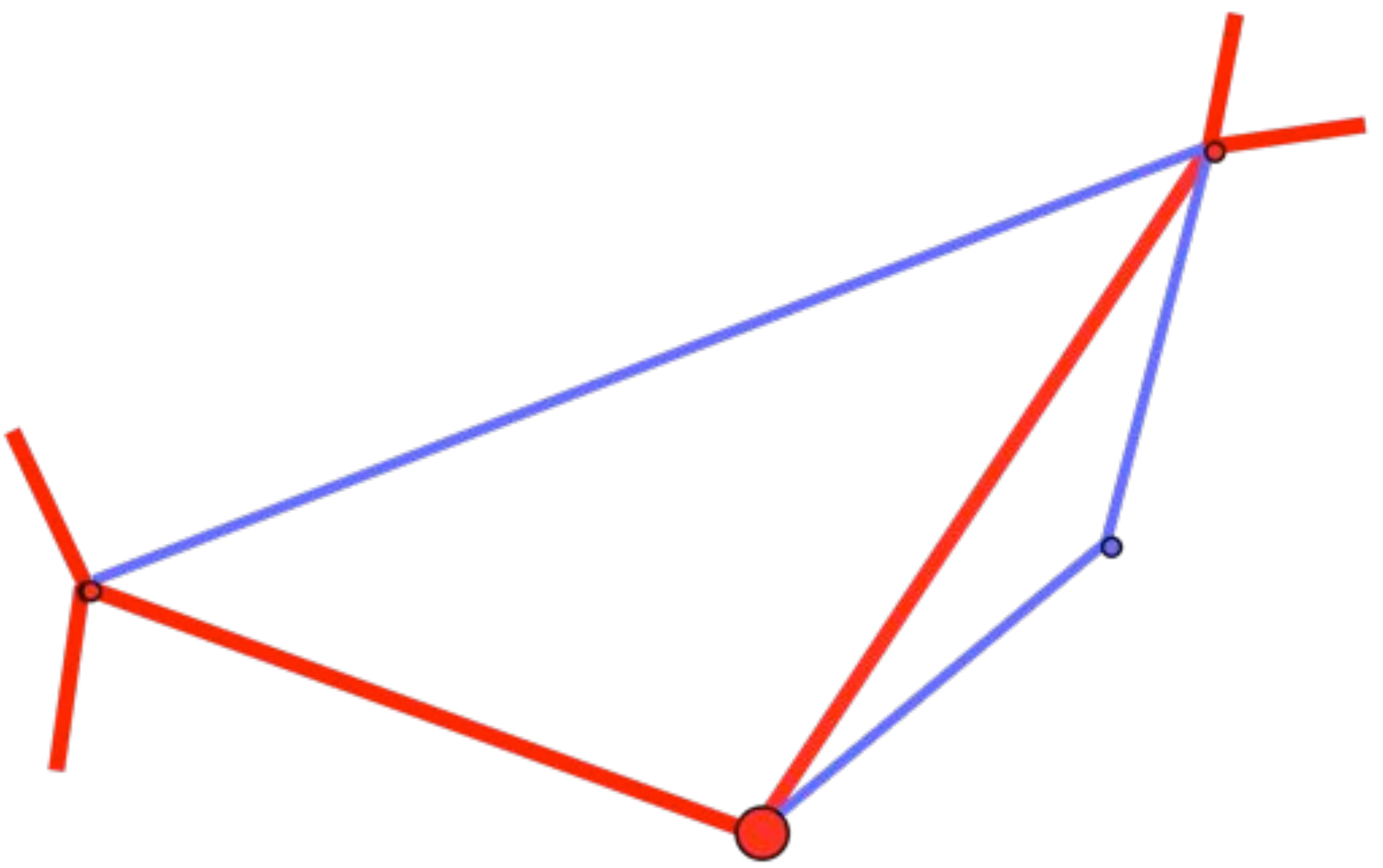, scale=.035}} 
  \caption{Given a spoke to a hub  (a), contracting may make an adjacent path unclear (b) (c).} \label{fig:Hubunclean}
 \end{center}
 \end{figure}

We note that in contracting a spoke in one disc, the boundary of an adjacent disc could become unclean (Figure~\ref{fig:Hubunclean}).  
In that case, we say the spokes in that disc are not contractible and we do not attempt to contract them in this step, 
however, we will contact spokes contained in other clean discs.  

Since vertex splitting is the reverse of edge contraction, the Proposition follows.
\qed

We begin by executing the first three steps in the order given.  Step 1 ensures that every disc is clean, and Step 2, contracting interior edges, does not 
cause any of the boundaries of the discs to become not clean.  However, Step 3, contracting spokes, has the potential to cause a boundary path to become
not clean (see Figure~\ref{fig:Hubunclean}).  

Prior to executing Step 4,  we have a conditional statement that has us cycle 
back to Step 1, in the case when some boundary path is no longer clear.   However, in this looping back, the number of vertices must have been reduced from the previous pass through Step 1, and the lengths of some paths have also been reduced while others remained constant. 
The loop eventually terminates because there are a finite number of vertices in the discs and on the paths.  The final polyhedron $\Phat_{3}$ still comes from a well-designed block and hole polyhedron, with the same blocks and holes.

In the end, perhaps after several cycles through Steps 1-3, 
$\Phat_3$ is  a well-designed polyhedron and  1) all paths are clear; and 2) there are no vertices interior to any discs in $\Phat_3$.  \\

\noindent {\bf Step 4}: Contract path edges to length 1. 

The next step is to shrink all paths between discs of $\Phat_3$ to length 1. 

\begin{proposition} [{\bf Step 4}]
If application of Step $4$ to the well-designed polyhedron $\Phat_3$ output from the final pass through Step $3$ yields the polyhedron $\Phat_4$, then 
\begin{itemize}
\item  [(i)] each pair of well-attached discs share a path of length 1 in $\Phat_4$; 
\item  [(ii)] $\Phat_3$ can be obtained from $\Phat_4$ by a construction sequence of vertex splits;
\item [(iii)] if $\Phat_4 \neq \Phat_{3}$, then the number of path vertices has decreased; and
\item[(iv)] $\Phat_4$ is topologically equivalent to $\Phat_{3}$ and therefore $\Phat_4$ is a well-designed block and hole polyhedron.
\end{itemize} 
\label{shrink_paths}
\end{proposition}

\proof
By Lemma~\ref{pathshrink}, each path edge is long and we may contract these to get a path of length 1. 
Contracting a path edge keeps the path clear because contracting doesn't create any new edges between pairs of vertices on the path. 
It does not change the underlying topology of the block and hole polyhedron. 

Since vertex splitting is the reverse of edge contraction, the Proposition follows.   
\qed

When we complete Step 4, all paths are length $1$.  The final polyhedron $\Phat_{4}$ comes from a well-designed block and hole polyhedron with the same blocks and holes as the original polyhedron $\Phat$, 
with all boundary paths of length $1$ and no vertices interior to the faces.  We call such a polyhedron {\it simplified}, and it comes from a base polyhedron 
of the topology we have expanded, though perhaps with a retriangulation of the surface discs. 


\begin{proposition} [{Contraction Sequence}]
Application of the Contraction Sequence (Algorithm~\ref{contractionSequence}) to the polyhedron $\PEhat$ of a well-designed block and hole polyhedron yields the simplified polyhedron $\Phat$ of a base 
well-designed block and hole polyhedron with an equivalent topology and all paths of length 1.  
 
$\PEhat$ can be obtained from $\Phat$ by a construction sequence of vertex splits. 
\end{proposition}

  \proof This follows from Propositions \ref{cleaning}, \ref{only_hubs}, \ref{wheelremoval}, and \ref{shrink_paths}.
\qed

\begin{corollary} The application of the Contraction Sequence and then the reverse vertex splits to a well-designed polyhedron $\PEhat$  yields the polyhedron
$\Phat^*$ which is essentially the same as the initial block and hole polyhedron, with some rerouting of boundary paths (reversing Step 1) and with some triangulation of the surface discs.
\end{corollary}
 
\proof 
Follows from assumptions and previous Propositions.  
\qed

\subsection{Main Result}\label{sec:main}

We apply the steps of Algorithm~\ref{contractionSequence}  to the expanded block and hole polyhedra $\PEhat$, to show that the expansion from a base polyhedron  can be accomplished by 
a sequence of vertex splits.  
  
Recall that an expanded polyhedron is created by starting with a base polyhedron $\poly$ in which we select certain faces of $P$ to  become holes and other faces of $P$ to become blocks (with isostatic graphs on each block, including the initial edges of the block face), forming the initial block and hole polyhedron $\Phat$.   
The remaining faces of $\Phat$ were called { surface discs}. 

Next, we expanded $\Phat$ by inserting vertices to split the edges of $\Phat$ which are not on the boundary of
either a block or a hole.   Finally we inserted additional vertices inside surface discs and triangulated these discs, on their boundaries and all their vertices.

This proposition confirms that all the conditions are met to allow us to apply Step 1 of the Contraction Sequence to clear the paths of in an expanded polyhedron $\PEhat$.

\begin{proposition} If the initial polyhedron $\Phat$ is a well-designed block and hole polyhedron, then $\PEhat$ is a well-designed block and hole polyhedron. 
\end{proposition}

\proof We note that the surfaces faces of a polyhedron are well-attached and well-surrounded, and the construction of the expanded polyhedron from 
a given polyhedron ensures that the topology of the two structures remain the same.   The proof now follows from the definitions and the prior results.
\qed

We are now ready to state the main theorems for expanded block and hole polyhedra.  The proofs are essentially contained in the preceding discussion. 

\begin{theorem} [{\bf Expansions as Vertex Splits}] An expanded block and hole polyhedron $\PEhat$, 
can be reduced by a sequence of contracting long edges to  the base block and hole polyhedron ${\PBhat}^{*}$
with some triangulation of the surface faces of the original $\PBhat$.  

Conversely,  the expanded block and hole polyhedron $\PEhat$ can be created from some triangulation of the
base block and hole polyhedron $\PBhat$ by a sequence of vertex splits.
\label{contractionsequence}
\end{theorem}

Since vertex splitting preserves first-order rigidity, and independence, we now have achieved our initial goal.

\begin{theorem} [{\bf Rigidity of Expanded Polyhedra}]
If a base polyhedron $\PBhat$ is generically isostatic (rigid, independent,
respectively) for
all triangulations of its surface faces,  
then every expanded polyhedron $\PEhat$ on this base is generically isostatic (rigid, independent,
respectively).\label{thm:ExpandedPolyhedron}
\end{theorem}

Theorem~\ref{thm:ExpandedPolyhedron} now sets the task of identifying base polyhedra which are generically isostatic for all triangulations of the surface faces.  We will illustrate that with an example below, and give some initial results in the next section.  

\subsection{Extensions} \label{sec:extensions}

\begin{figure}[h]
 \begin{center}
   \subfigure[] { \psfig{file=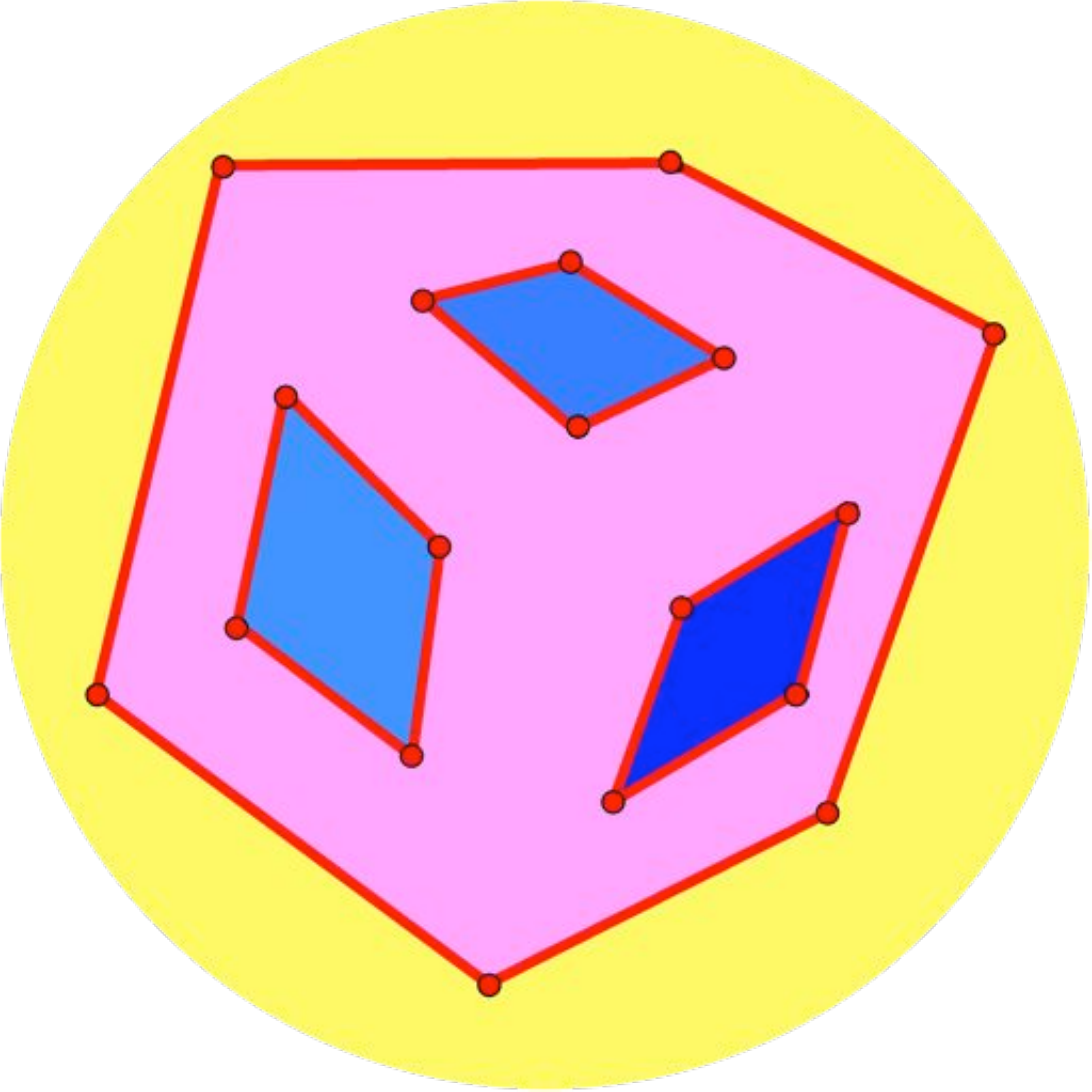, scale=.045}}\quad
 \subfigure[]{ \psfig{file=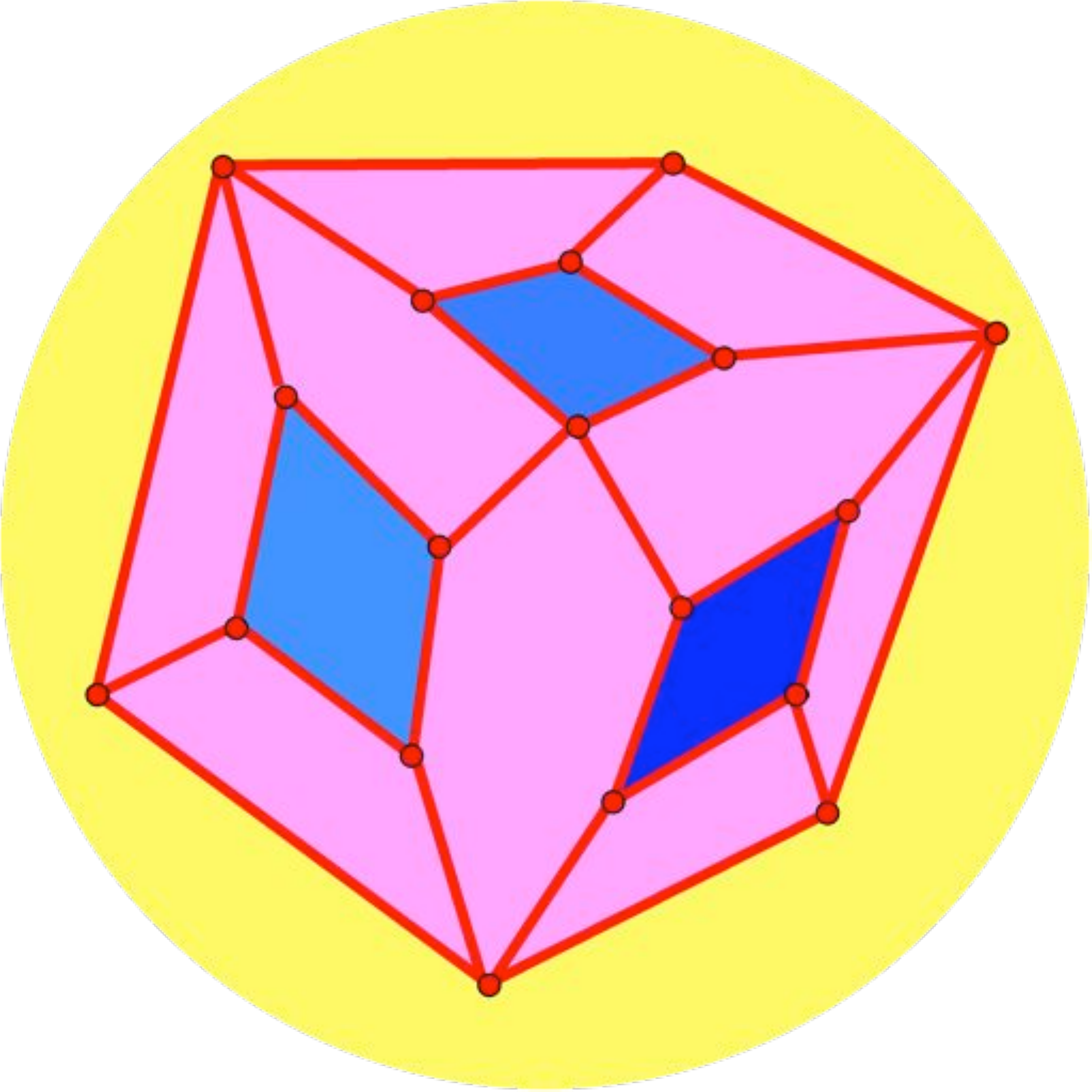, scale=.045}} \quad
   \subfigure[]{ \psfig{file=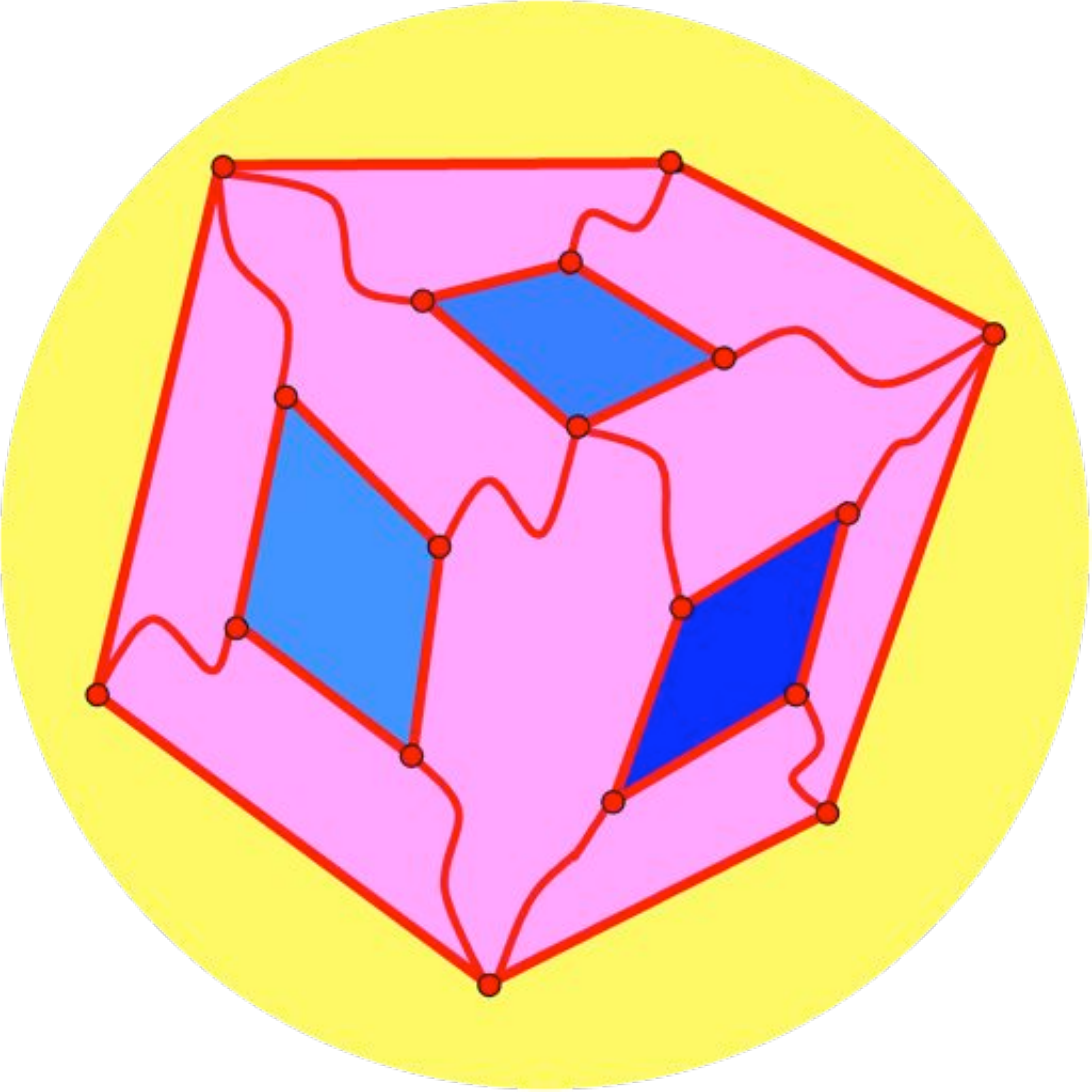, scale=.045}} 
  \caption{Given a spherical polyhedron with identified blocks and one exterior hexagon hole (a), we can locate a possible base polyhedron (b), and  search for paths that would lead us there, under contraction (c)}\label{fig:HexHole}
 \end{center}
 \end{figure}

Figure~\ref{fig:HexHole} illustrates another setting where we can apply these techniques.  The goal is to confirm that a selected spherical polyhedron $\poly$ with designated blocks and holes  and other faces triangulated (Figure~\ref{fig:HexHole}(a)), is indeed generically isostatic.  First we find a simple candidate base polyhedron with the same pattern of blocks and holes, which we can check is generically isostatic  (Figure~\ref{fig:HexHole}(b)).  We then  search for paths which contract the given polyhedron to the base, using the methods above.  When this is found (Figure~\ref{fig:HexHole}(c)),  our results on contractions demonstrate that the originally given block and hole polyhedron is generically isostatic.  

\begin{figure}[h]
 \begin{center}
   \subfigure[] { \psfig{file=HexHole4, scale=.05}}\quad\quad
 \subfigure[]{ \psfig{file=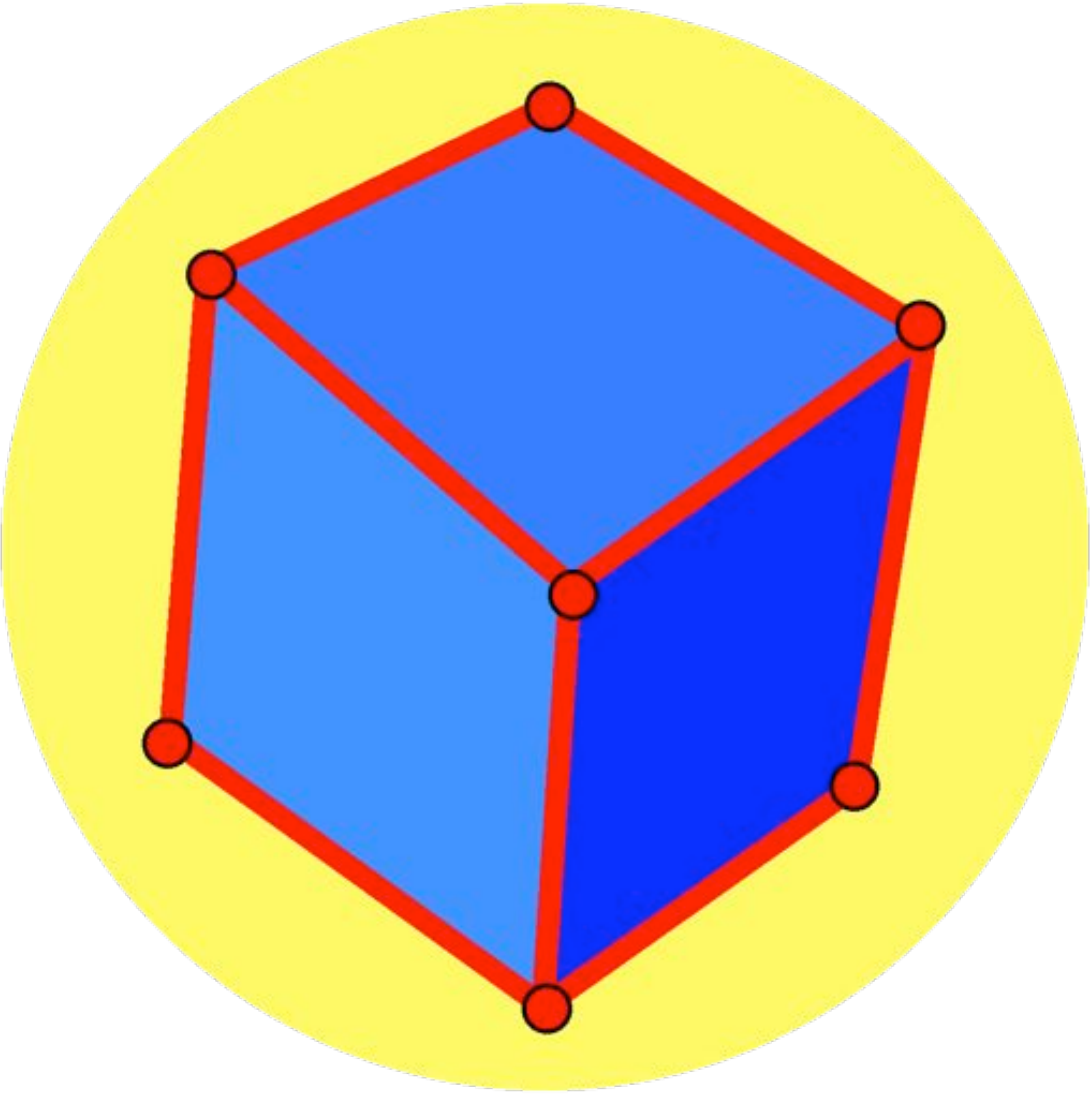, scale=.10}} \quad
  \caption{Given a spherical polyhedron with identified blocks and holes (a), and paths of length 1, we can locate a possible simpler base polyhedron with paths shrunk to 0 (b), which is easily identified as generically rigid.}\label{fig:HexHoleS}
 \end{center}
 \end{figure}

However, the polyhedron $\PBhat$ in Figure~\ref{fig:HexHole}(b) is still not very simple to check for generic rigidity. It is substantially simpler to consider the base as the further reduced polyhedron $\Phat$ in Figure~\ref{fig:HexHoleS}(b).  This modified base has one large hole (the exterior hexagon) and three blocks and is easily checked.   To do this next step, we need to consider when we can contract paths of length 1 between triangulated discs to paths of length 0, without disturbing the isostatic counts etc.. 

One answer to when we can safely contract paths of length 1 between triangulated discs is when we have a guide in the form of a known, isostatic base.  If the previous polyhedron also had the isostatic count $|E|=3|V|-6$, then the contraction only removed three edges  (to preserve the count when identifying the vertices). This contraction is then reversed by vertex splitting, which is all that we desire.   We record this as a further corollary.

\begin{corollary}[Simple Base Contraction] Let $\Phat$ be an original (expanded) polyhedron with identified blocks, holes, triangulated discs and boundaries which is well-designed.  Let $\PBhat$ be  base polyhedron with paths between pairs of triangulated discs contracted to length $1$ or $0$.  If  $\PBhat$ is generically isostatic for all re-triangulations of its triangulated discs, then 
\begin{itemize}
\item[(a)]  there is a sequence of vertex splits from a triangulation 
of  $\PBhat$ to $\Phat$.
\item[(b)] the initial block and hole polyehdron $\Phat$ is generically isostatic.
\end{itemize} \label{cor:simplebase} 
 \end{corollary} 

\proof  The proof follows from Theorem~\ref{thm:ExpandedPolyhedron} and the observations in the previous paragraph.
\qed 

For the initial example in Figure~\ref{fig:HexHole}(a) we have a complete verification of the generic isostaticity provided we can find the paths in Figure~\ref{fig:HexHole}(c) which can then be contracted to the generically isostatic base Figure~\ref{fig:HexHoleS}(b).  These are all reversed by vertex splitting, which preserves the generic isostaticity.

The Swapping Principle \cite{frw} applied to this example confirms that the same analysis applies to larger polyhedron with a hexagonal base and three windows.  More generally, this route can be used, for example, to analyze the grounded geodesic domes with desired windows and doors which were mentioned in the introduction.  \\

There is another corollary implicit in the algorithm. If we apply Steps 1-3 to any triangulated disc within {\it any} larger framework with a  boundary within the larger framework, the methods will apply provided that there are no edges outside the disc which connect two boundary vertices.   We call the boundary of such a disc {\it clear}.  We can show that the triangulated disc with a clear boundary can be obtained from 
a ``smaller" triangulated disc with the same boundary, and with no interior vertices.
  
\begin{corollary}[{Disc Clearing}] Let $D$ be a triangulated disc with a clear boundary.  Then there is a sequence of vertex splits from a triangulation 
$D^{*}$ of the same boundary, 
with no interior vertices, ending with $D$.  
\label{nointeriorvertices} 
 \end{corollary} 

If we know that one such $D^{*}$ makes the modified graph $G^{*}$ infinitesimally rigid, then some triangulations of the vertices of $D$ will make the modified graph $G^{*'}$  generically rigid.  However, this might not include the desired triangulation $D$. 
Note that, given a disc $D$ with a clean boundary within a larger graph $G$.  To effectively use these contractions without careful examination of the options, we need to know that {\it all possible} triangulations of $D^{*}$, embedded in the graph $G^{*}$ with only this disc modified, give generically rigid graphs.  We can then conclude from the algorithm that $G$ is generically rigid.

\section{Creating base polyhedra}\label{sec:base}

If our reduction has brought us back to a base block and hole polyhedron, with surface faces larger than a triangle, 
we may end with any of the possible triangulations of these faces.  One way to avoid checking all possible triangulations is to ensure that our base case has only blocks, 
holes and surface faces which are triangles.  

As in the previous subsection, one can  go to even smaller base polyhedra by shrinking paths of length $1$ between surface discs to paths of length $0$ (identifying vertices).  This changes the topology - and will require the added guidance from a known isostatic base.

There are some geometric alternatives in which non-triangular surface faces are realized as plane faces, with no three vertices collinear and with one of the triangulations known to be isostatic.  When such realizations are provided, the Isostatic Replacement Principle used in \cite{infp1} can be applied to substitute any one triangulation of the plane face  for any other triangulation with the same boundary edges.   We will not pursue this alternative here, but there may be situations where that approach is helpful. 

\subsection {Block and hole $n$-towers.} \label{sec:towers}
We demonstrate that one basic set of examples are generically isostatic.  
These examples will be explored further in Section~\ref{sec:allostery}. 

\begin{figure}[h]
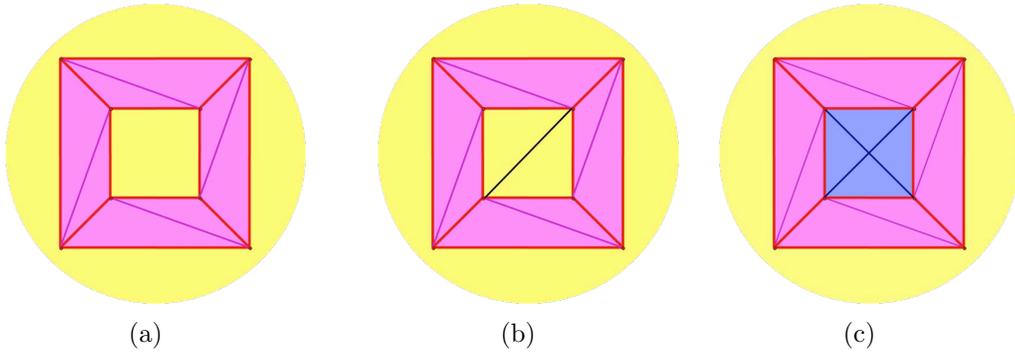

 \begin{center}
   \subfigure[] { \psfig{file=P44_0, scale=.027}}\quad\quad
  \subfigure[]{ \psfig{file=P44_1, scale=.027}}\quad
  \subfigure[]{ \psfig{file=P44_2, scale=.027}}
 \caption{Adding edges to a cylinder $\cylinder(4,4,4)$ to generically rigidify it. }
\label{fig:C444}
 \end{center} 
 \end{figure}

A {\it block and hole $n$-tower $\tower_{n}$} is a block and hole polyhedron with one block of size $n$ (thought of as the ground of the tower),  one hole of size $n$ 
(thought of as the open top of the tower), and a tube or cylinder of surface triangles between them (see Figure~\ref{fig:C444}).  
The block and hole polyhedron is {\it proper} if there are $n$ vertex-disjoint paths between the block and the hole.  
We make these $n$ paths, plus the boundaries of the block and hole into the boundary paths of the block and hole $n$-tower polyhedron.   It is simple to check that such an $n$-tower has $|E|=3|V|-6$, so it has the potential to be generically isostatic.  

\begin{figure}[h]
 \begin{center}
 \psfig{file=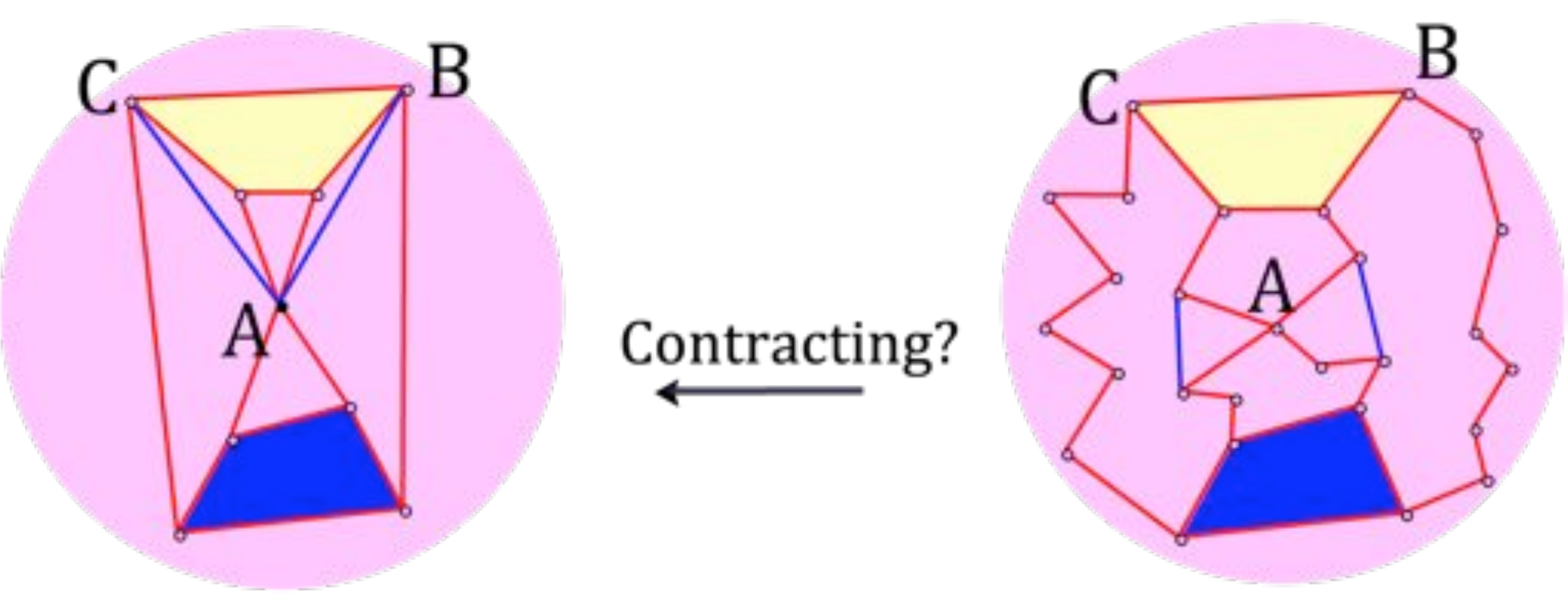, scale=.065}
 \caption{An expanded polyhedron (b) where the boundary paths do not $4$-connect the $4$-hole 
 and the $4$-block can contract to a base which is not $4$-connected (triangle ABC in (a)) . }
\label{fig:non4connected} 
\end{center} 
 \end{figure}

The existence of $n$ vertex-disjoint paths is equivalent to $n$-connectivity of an $n$-block and $n$-hole in the tower $\tower_{n}$ \cite{menger}.  
That is, there is no set of fewer than $n$ vertices which separates the block from the hole in the underlying graph.  
If $\tower_{n}$ was not $n$-connected, we can take the separating vertices, edges among these, along with the component which contains the block, and verify that this subframework has $|E'|> 3|V'|-6$ and has more than $2$ vertices.  
This subframework has a self-stress, so the overall framework cannot be infinitesimally rigid (or isostatic).  Thus we will focus on the proper $n$-towers.   Figure~\ref{fig:non4connected} shows such a non-example. 

\begin{figure}[h]
 \begin{center}
 \psfig{file=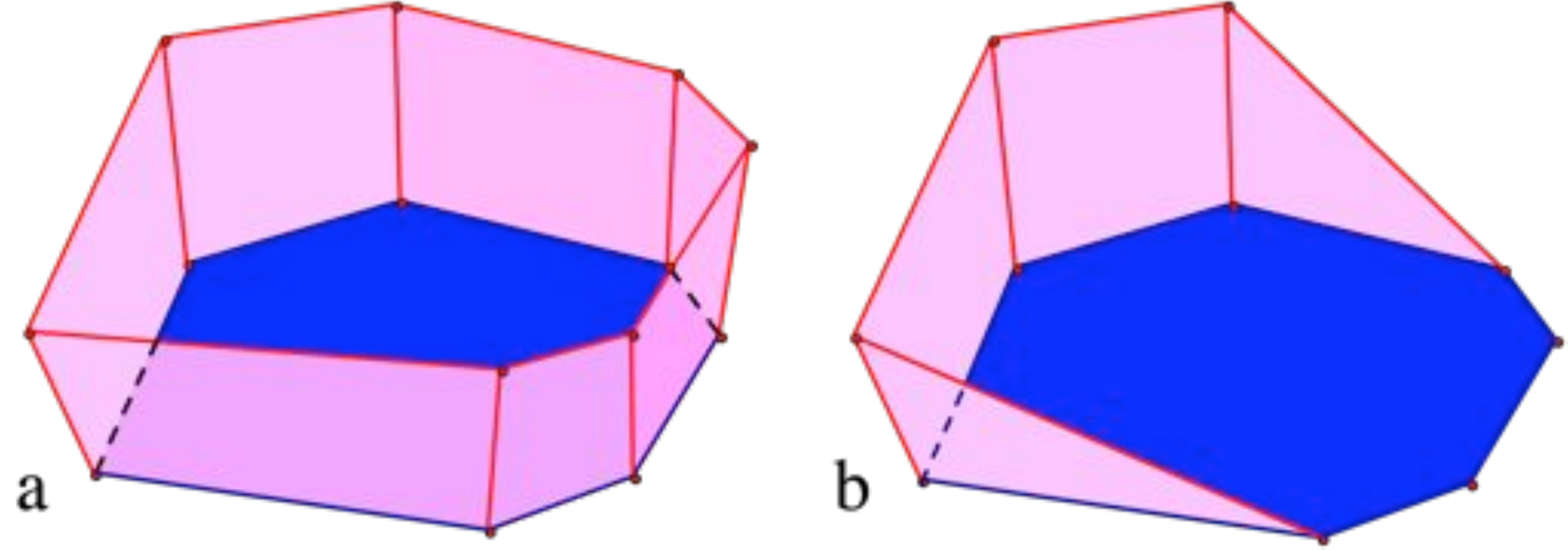, scale=.06}
 \caption{Sample cylinders with connecting vertex-disjoint paths of length $0$ and $1$. }
\label{fig:cylinderwall}
 \end{center} 
 \end{figure}

If we start with any proper $n$-tower, and apply the lemmas above, we can reduce down to a framework with all boundary edges between triangulated discs down to length $0$ or $1$. (Note that some or all of the the $n$ vertex-disjoint paths might just be  shared vertices on the block and the hole, so no further contracting is needed; see Figure~\ref{fig:cylinderwall} (b).)   

 Our desired base is single $n$-gon which is both blocked with additional edges to be a block, and is recorded as a second $n$-gon hole to make the base polyhedron.   The framework is the framework of the block, which is isostatic (both independent and rigid).   With this generically isostatic base, we apply Corollary~\ref{cor:simplebase} to shrink remaining paths of length $1$ to length $0$ and come to the base.  We conclude that every proper block and hole $n$-tower  $\tower_{n}$ is generically isostatic.  

\begin{theorem}[\bf Rigidity of $n$-towers]
Every proper $n$-tower is generically rigid.
\label{thm:$n$-tower}
\end{theorem}

\noindent {\bf Remark.}  At a geometric level, the infinitesimal rigidity of $(\tower_{n},p)$ is the same as the `swapped block and hole polyhedron' 
$(\widetilde{\tower_{n}},p)$ where the hole is turned to a block and the block is turned to a hole.  This an example of the general principles developed in \cite{frw}. 

\subsection{Necessary conditions for isostatic block and hole polyhedra.} \label{sec:necessary}
As we saw earlier, a triangulated spherical polyhedron is isostatic, 
with $|E|=3|V|-6$.  If we want a general  block and hole polyhedron to be isostatic, then the number of bars removed to create the 
holes should equal the number of edges inserted to make the blocks isostatic.   
If hole $H_{j}$ is a $h_{j}$-gon, then we remove $h_{j}-3$ edges.  If block $B_{i}$ is a $b_{i}$-gon, then
the number of edges to be added to the existing triangulation is $b_{j}-3$.  (Simply imagine a general postion
`convex' polygon in space with the top and bottom still triangulated to make a convex triangulated sphere. )  This observation gives the following neccessary condition. 

\begin {proposition} If a block and hole polyhedron $\poly$ is isostatic, then 
$$ \sum_{i\in \B} (b_{i} -3)= \sum_{j\in \Ho} (h_{j}-3)$$
\end{proposition}

Our experience with connectivity conditions for a single block and single hole demonstrates that this count, alone,
is not sufficient, if these are not sufficiently connected. (See Example~\ref{example5}) .  The counter-example involves some disconnecting polygon which, if treated as a hole or a block components
 would give one, and therefore both of the components an unbalanced count.  

\begin{definition}  A {\rm cut cycle} $\mathcal{C}$ is a sequence of $c$ vertices 
$(v_{1}, \ldots,  v_{i}, v_{i+1}, \ldots, v_{c})$, such that each adjacent pair 
$v_{i}, v_{i+1}$ is either an edge of the  base polyhedron $\Phat$, or if the pair is not such an edge,  then the two vertices share a triangulated disc.   
\end {definition} 

Note that if two vertices share a surface face, then their edge {\it might} be an edge of the triangulation of the face in $T$. 
Any claim we make about the rigidity for {\it all triangulations} needs to allow for this possibility, which this definition does.  

\begin{proposition} If a well-designed block and hole polyhedron $\Phat$ is generically isostatic for arbitrary retriangulations of the surface discs, then for any 
cut cycle $\mathcal{C}$ of size $c$, each component $\Phat_{C}{'}$ satisfies
$$ \sum_{i\in \B'} (b'_{i} -3) - \sum_{j\in \Ho'}( h'_{j}-3 )\leq c - 3. $$ \label{thm:cuts}
\end{proposition} 

\proof  If we have a component with
$\sum_{i\in \B'} (b'_{i} -3) - \sum_{j\in \Ho'}( h'_{j}-3 )\geq c - 3 $ then
this component is already overcounted, and there is a self-stress. \qed
 
The cut set condition 
$$ \sum_{i\in \B'} (b_{i} -3) - \sum_{j\in \Ho'}( h_{j}-3 )\leq c - 3 $$
is implicitly symmetric in blocks and holes and  is equivalent to  the condition
$$ \sum_{j\in \Ho''}( h_{j}-3 ) - \sum_{i\in \B''} (b_{i} -3)  \leq c - 3 $$
 on the other component because the overall count is:  $ \sum_{i\in \B' \cup \B''} (b_{i} -3)= 
\sum_{j\in \Ho' \cup \Ho''}( h_{j}-3 )$.  This equivalence of necessary conditions for generically isostatic block and hole polyhedra, 
if we swap blocks and holes is confirmed in \cite{frw}.   

The set of necessary conditions is not sufficient , as the following example confirms.  
 \begin{figure}[h]
 \begin{center}
   \subfigure[] { \psfig{file=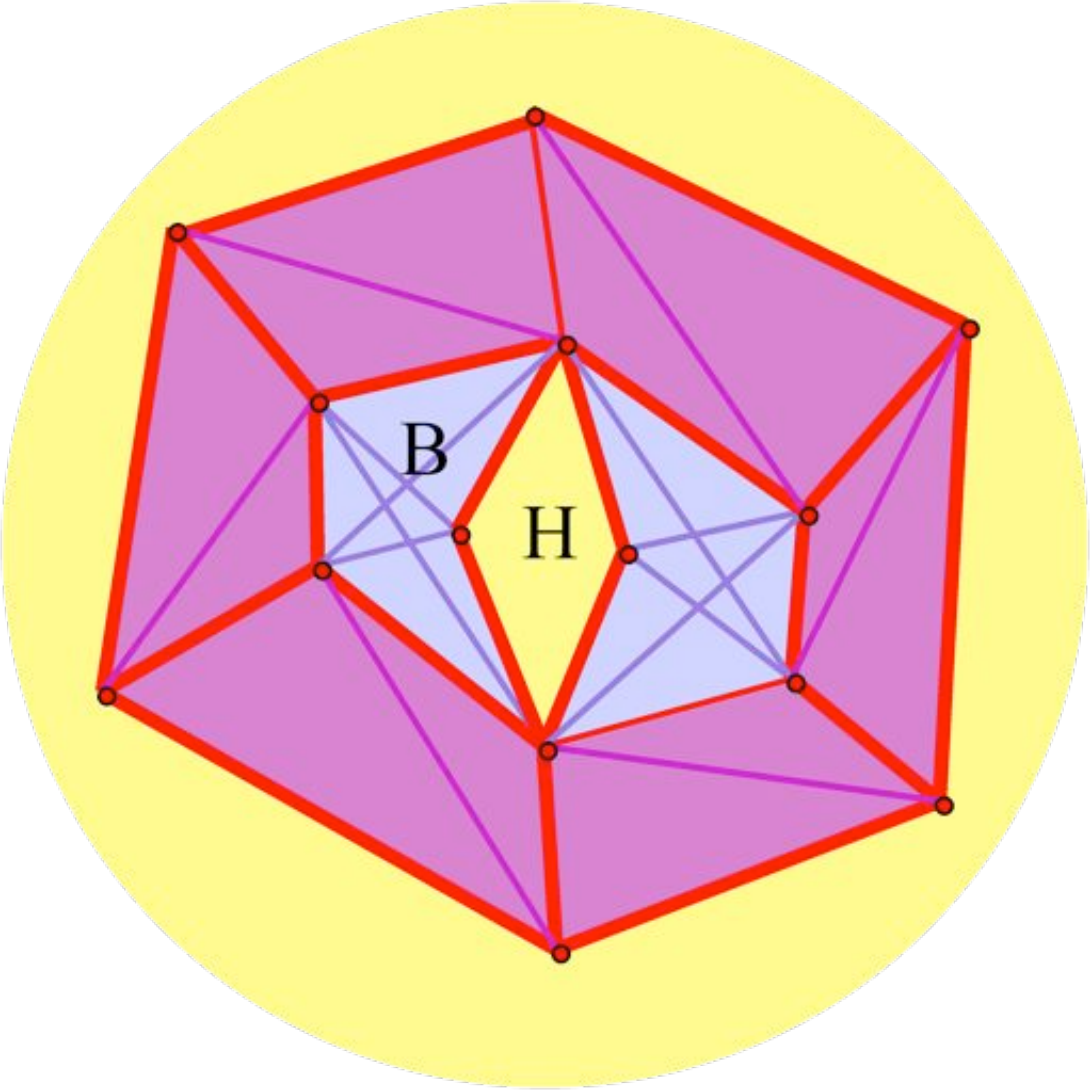, scale=.05}}\quad\quad\quad
  \subfigure[]{ \psfig{file=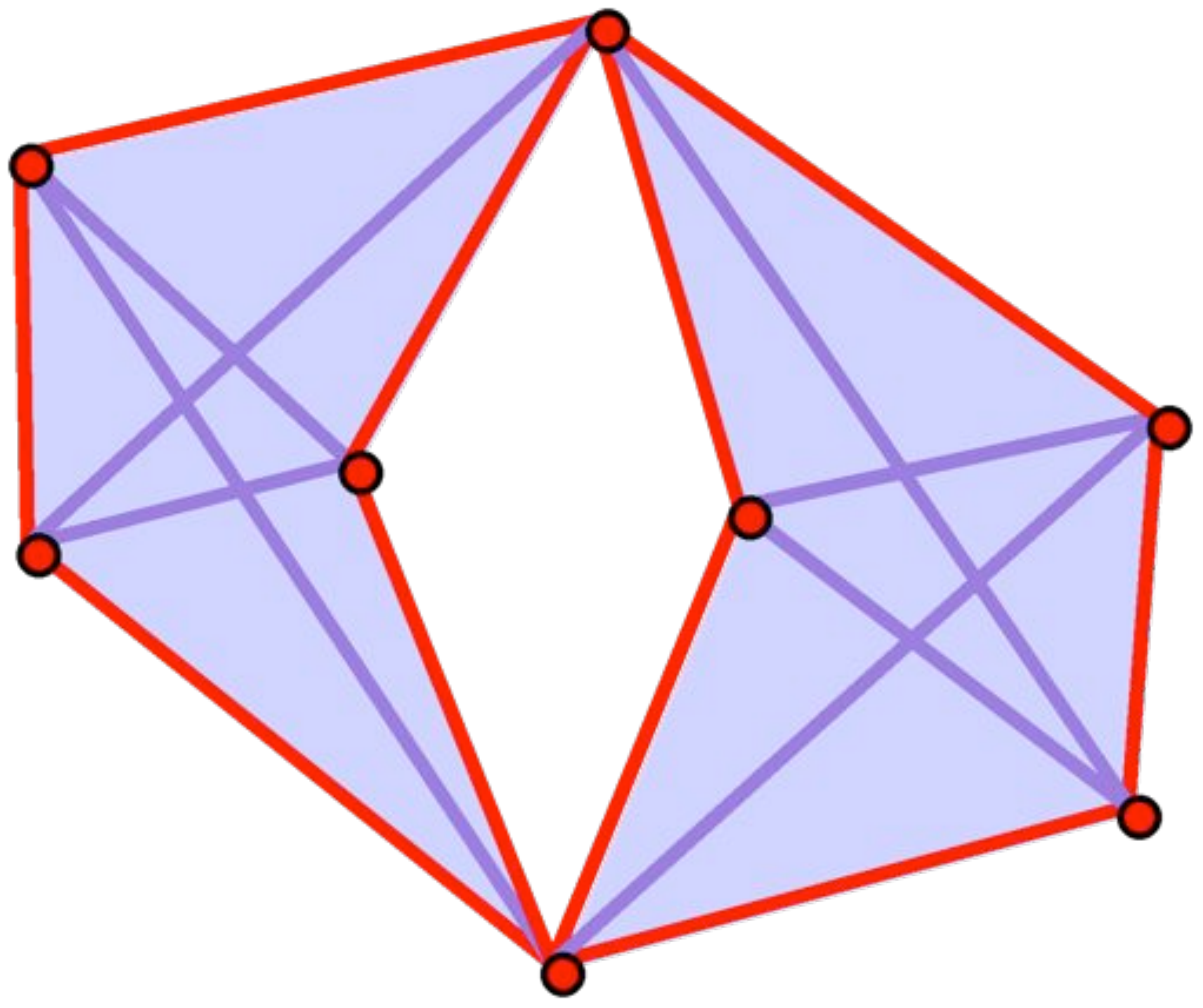, scale=.13}}
 \caption{A block and hole polyhedron which is generically dependent. }\label{fig:implicitbanana}
 \end{center}
 \end{figure}
 
\begin{example} {\rm Consider the block and hole polyhedron in Figure~\ref{fig:implicitbanana}(a) 
In this figure, we have two quadrilateral holes and two quadrilateral blocks.  All the counts above are respected.  
However, the two blocks in Figure~\ref{fig:implicitbanana} (b) are clearly the classical `double banana' and are dependent.  Therefore the larger framework is dependent 
and not isostatic (or rigid).}
\end{example}

  Are there additional necessary conditions which might, together, make sufficient conditions? For example, we need to attend to whether blocks (and holes) are well connected.  That is: whether two blocks which share two vertices also share the edge between these vertices.   
The reader will note that these next codnitions are also `dual' in the use of `hole' and `block'. 

\begin{proposition} If a block and hole polyhedron $\Phat$ is generically isostatic, then: 
\begin{enumerate}
\item two blocks intersect only on $0 \leq k\leq 2$ vertices, and if they intersect on $2$ vertices, they also share the edge between these vertices;  

	two holes intersect only on $0 \leq  k\leq 2$ vertices, and if they intersect on $2$ vertices, they also share the edge between these vertices.  
	
\item  no two vertices of a block are connected by an edge which is not in the block;  

         no two vertices of a hole are connected by an edge which is not in the boundary of the hole.  
\end{enumerate} \label{thm:separation}
\end{proposition}

\subsection{Conjecture on sufficient conditions \label{sec:sufficient}
for isostatic block and hole polyhedra.}
We conjecture that these, together with the counts of Proposition~\ref{thm:cuts} are sufficient. 

\begin{conjecture} A block and hole polyhedron $\Phat$ which satisfies the separation conditions of Proposition~$\ref{thm:separation}$ is generically isostatic for all retriangulations of its triangulated discs if and
 only if,  for every cut cycle $\mathcal{C}$ of size $c$, each component $\Phat_{C}{'}$  satisfies
$$ \sum_{i\in \B'} (b'_{i} -3) - \sum_{j\in \Ho'}( h'_{j}-3 )\leq (c - 3)\rm{.} $$
 \label{cycle-condition}
  \end{conjecture}

Using some simple counting arguments, the reader can verify that this condition is equivalent to the alternative formulation: 

\begin{conjecture} A block and hole polyhedron $\poly$ which satisfies the separation conditions of Proposition~\ref{thm:separation} is generically isostatic for all retriangulations of its triangulated discs  if and
 only if,  overall $|E| = 3|V| -6$ and for every subset $V'$ with at least three vertices, 
$$ |E(V'|) \leq 3|V'| -6.  $$
\label{count-condition}
  \end{conjecture}
  
If verified, the advantage of this alternate conjecture is that this counting condition can be rapidly checked
  using the `pebble game algorithm' \cite{pebble}, and the separation conditions are local around the holes and blocks.

\section{Mathematical Allostery}\label{sec:allostery}

There is a larger theory for applying methods from rigidity theory to many aspects of protein behaviour and function \cite{Counting}. Within our bodies, a ligand (for example a hormone, a neurotransmitter or a drug) binds to a protein, 
causing a change in the rigidity of the protein.  In some proteins, changes in rigidity at one site will in turn effect the shape and  flexibility at a distant site (allostery) and alter the function of the protein \cite{allostery}.  

Motivated by this important phenomenon of transmission through a protein, we investigate the ability of a class of `toy' mathematical examples to model the transmission of rigidity / flexibility across a
structure between two sites - what is called {\it mathematical allostery}.   
The simplest objects in the class are triangulated spheres with 2 holes.  
Clearly these structures are not rigid ($|E| <3|V|-6$).  
We investigate what happens when one of the holes is filled in one edge at a time - up to the point this first hole is replaced with a block which makes one side rigid,  
simulating a ligand binding to a site on the protein.  
Under what conditions is this new structure rigid? or stressed?  How does this change the degree of freedom (dof) of the other end (if at all)?  
Equivalently, if we squeeze the one hole (i.e. mechanically change the shape as binding might) will the other side change shape?  
Our examples will all be generic - in keeping with the approach of the rest of the paper.  \\

Informally, the {\it internal degrees of freedom (idof)} of a framework is the number of distinct non-trivial ways the framework can move.  Alternatively, it is the number of bars that need to be added to make the framework (infinitesimally) rigid.
We note that an $n$-gon hole has $n$ vertices and edges and so has $3n-6 - n = 2(n-3)$ idof as an isolated structure. 
It can have  $2(n-3)$ or fewer idof in the larger polyhedron, as we will see. \\

We can relate the concept of idof to the rigidity matrix equation $M \beta = {\bf 0}$ \cite{infp1, wchapter}. 
The idof of a given framework is the number of independent rows that can be added to the matrix.  
Thus, adding a bar to a framework may or may not reduce the idof by one because it is equivalent to adding a row to the rigidity matrix.   Recall that bar is  independent if it reduces the idof and is called dependent otherwise.  
In the later case, we are adding a {\it redundant} bar and a stress is created within the framework.   
When this is done with a framework with generic coordinates for the points, we are studying generic rigidity, and generic independence.  We speak of the {\it generic idof of a graph $G$} and the generic independence (redundance) of an edge.

\subsection{Simple Examples}\label{sec:allosteryexamples}
We will use the towers of Section~5.1 as our starting model.    A {\it waist} of the block and hole polyhedron, $\poly$, with two holes $H_{1}$ and $ H_{2}$ and no blocks is the smallest cycle 
vertices and edges in $\poly$ whose removal separates the vertices of the two holes of $\poly$.  Note that some vertices and edges of the holes may be in the cycle.  
By Menger's Theorem \cite{dirac}, a triangulated sphere with two holes will have a $k$-waist if, and only if, there are at most 
$k$ vertex disjoint paths from the vertices of one hole to the vertices of the other hole.  
For practical purposes, $k\leq \min \{m,n\}$ since we can consider the vertices of a hole as `separating' the hole from the rest of the polyhedron. 

We define  a {\it cylinder} $\cylinder(m, k, n)$ to be block and hole polyhedron with one $m$-gon hole, $H_1$ and one $n$-gon hole, $H_2$, separated by a $k$-waist in a triangulated cylinder.  
In this language, the proper $n$-towers, with the block replaced by a hole are $\cylinder(n, n, n)$.   
The general class of frameworks, built from a cylinder, with some edges (up to and including a block) inserted in the first hole, are called  {\it tubes}.

In the following examples, we add the bars of a possible block at $H_{1}$ one at a time, and examine the idof of the resulting structure after each bar is added.  Throughout this section we are working with generic configurations, and rigidity means generic rigidity.

Recall that a triangulated sphere (without holes) is generically isostatic.   
In particular, it has $|E| = 3|V|-6$ bars, 0 idof, and no bar of the triangulated sphere is redundant.  
When a triangulated sphere is modified with one or more holes, $(3|V|-6 - |E|)$ determines the number of missing bars. For each bar that is removed from the triangulated sphere to form a hole, a row of the rigidity matrix is removed and the idof is increased by one. 

\begin{example} \label{example1} 
{\rm Consider  $\cylinder(4, 4, 4)$ starting with no bars in $H_{1}$:  
\begin{itemize}
\item each hole is missing one edge from the triangulated sphere;
\item $3|V| - 6 - |E| = 2$, and $\poly(4, 4, 4)$ is not rigid with 2 idof; and
\item there are 2 idof still at each hole.   
\end{itemize}
The third claim above is confirmed by recalling that adding two edges to one of the quadrilaterals gives a rigid block, creating a proper $4$-cylinder, which is isostatic. 

In particular, this also shows that each edge added to $H_{1}$ removes 1 idof (is independent).  
With one bar added to $H_1$, we have:  
\begin{itemize}
\item $(3|V| - 6 - |E|) = 1$, and $\cylinder(4, 4, 4)$  is not rigid with 1 idof; and
\item there is 1 idof still at each hole.   
\end{itemize}
The second claim above is confirmed by recalling that adding one more edge to $H_{1}$ gives a rigid block,  creating a proper $4$-tower, which is isostatic, while adding 1 to the second hole creates a triangulated sphere.

When we insert a second bar into $H_1$,  we effectively replace $H_1$ with a 4-gon block.  
\begin{itemize}
\item $3|V| - 6 - |E| = 0$, and $\poly(4, 4, 4)$  rigid 1; and
\item there is $0$ idof  at each hole.   
\end{itemize}
Since we have assumed $H_{2}$ and this block at $H_{1}$ are 4-connected in a vertex sense, this tower is isostatic, overall. 

$H_1$ becomes rigid when the block is added, and this also removes 2 idof from $H_2$.   
We have {\it transmitted} the two constraints from the site at $H_{1}$ to the site $H_{2}$, altering its behaviour from a distance.  This is an example of mathematical allostery.  }
\end{example}

\begin{example}\label{example2}
{\rm More generally, consider $\cylinder(n, n, n)$, starting with no edges added at $H_{1}$.  We observe that:  
\begin{itemize}
\item each hole is missing $n-3$ edges;
\item $(3|V| - 6 - |E|) = 2(n-3)$, and $\cylinder(n, k, n)$ is not rigid with $2(n-3)$ idof overall;
\item there are $2(n-3)$ idof at each hole (the number of edges needed to make this hole into a block).  
\end{itemize}
Thus each of the $2(n-3)$ bars inserted into $H_1$ of $\poly(n, n, n)$ in order to build up this isostatic block and hole structure is independent.  As in Example \ref{example1}, each of these bars removes $1$ idof from $H_{1}$ and also from $H_{2}$. 
In the end, we have replaced $H_1$ with an $n$-gon block, that is,   we added $2(n-3)$ bars to that hole.  Since the assumption on the waist is equivalent to $n$-connectivity, by Theorem \ref{thm:$n$-tower},   
this is generically isostatic - both rigid and independent. 

We note that for $n>5$ there are bad ways to add these $2(n-3)$ edges to $H_{1}$; ways that do not generate an isostatic block, but create a redundant set of edges 
within this subframework.  We assume that the added edges are part of some isostatic block at $H_{1}$.

With these simple structures, all the $2(n-3)$ added independent constraints at $H_{1}$ {\it transmit} to reduced idof at $H_{2}$ as well. 
}
\end{example}

\begin{example}\label{example3}
{\rm Consider a cylinder $\cylinder(n, n, m)$ with $n<m$.  
\begin{itemize}
\item the two holes are missing $(n+m-6)$ edges, so the cylinder starts with  idof $=(n+m-6)$;
\item there are $2(n-3)$ idof at the first hole (the number of edges needed to make this hole into a block).  
\item there are $(n+m-6)$ idof at the second hole which is less than $(2m-6)$ - the number of edges needed locally to make this polygon into a block.   This hole has all the idof of the structure - but less than a simple $m$-gon would have. 
\end{itemize}
 \begin{figure}[h]
 \begin{center}
   \subfigure[] { \psfig{file=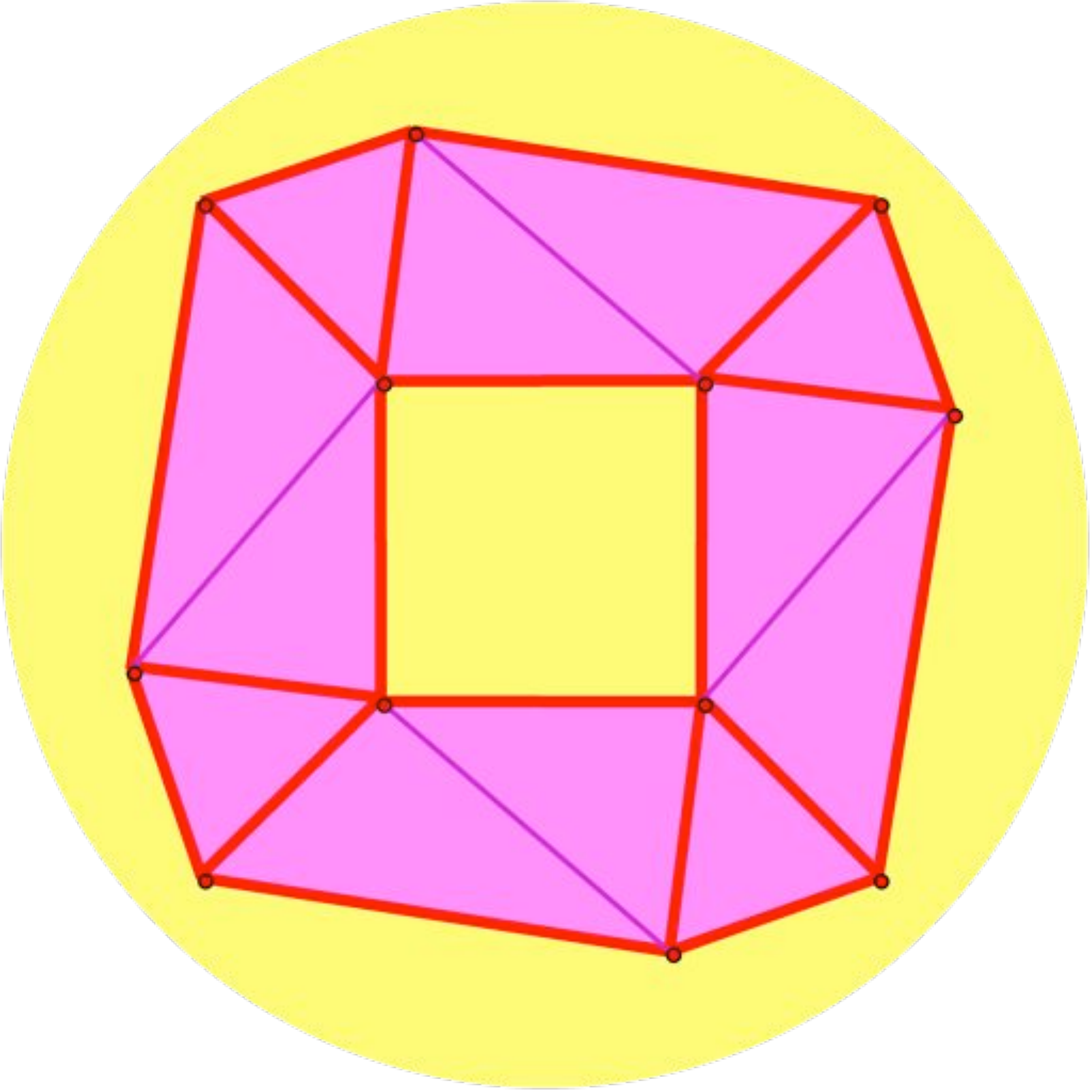, scale=.029}}\quad 
  \subfigure[]{ \psfig{file=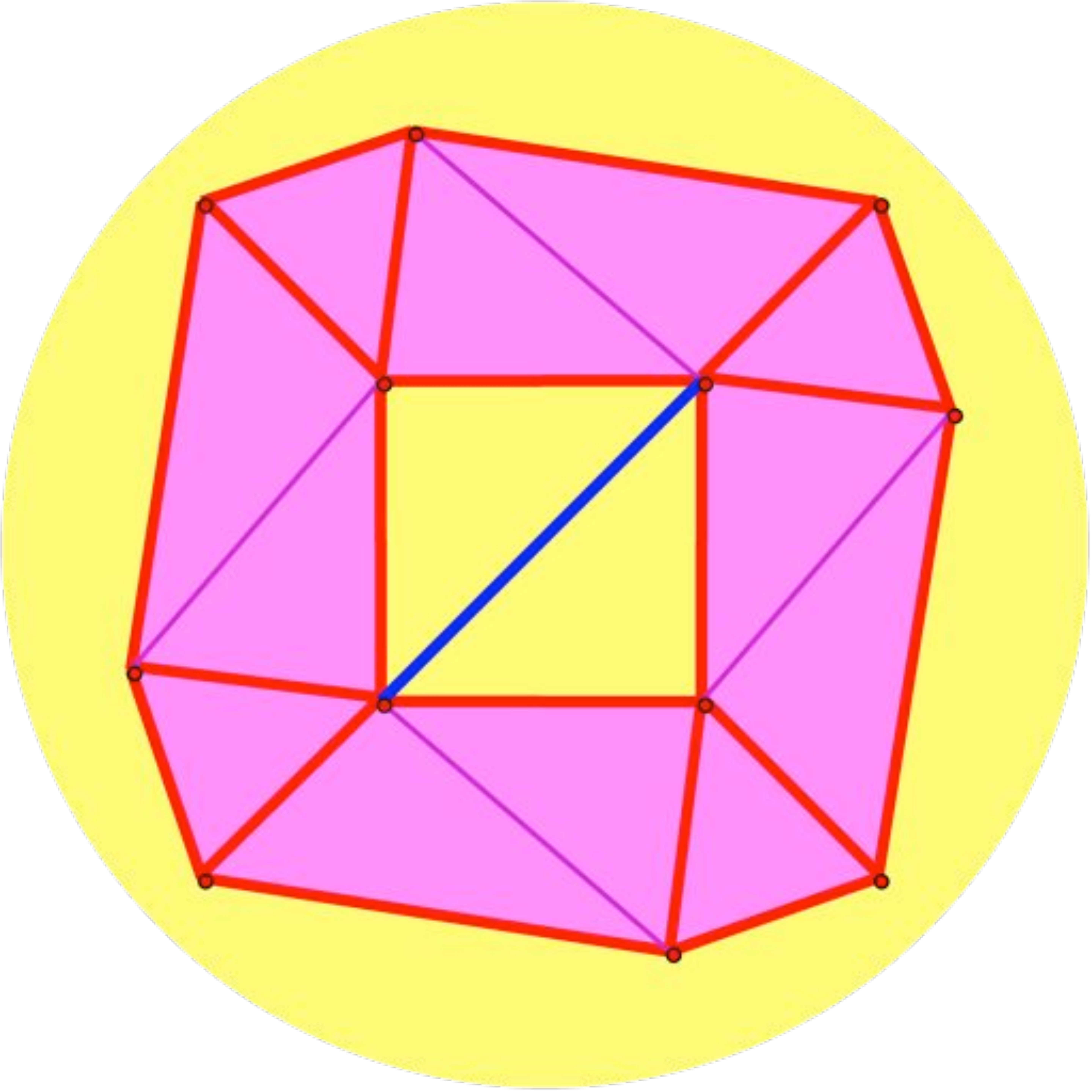, scale=.029}}\quad 
  \subfigure[]{ \psfig{file=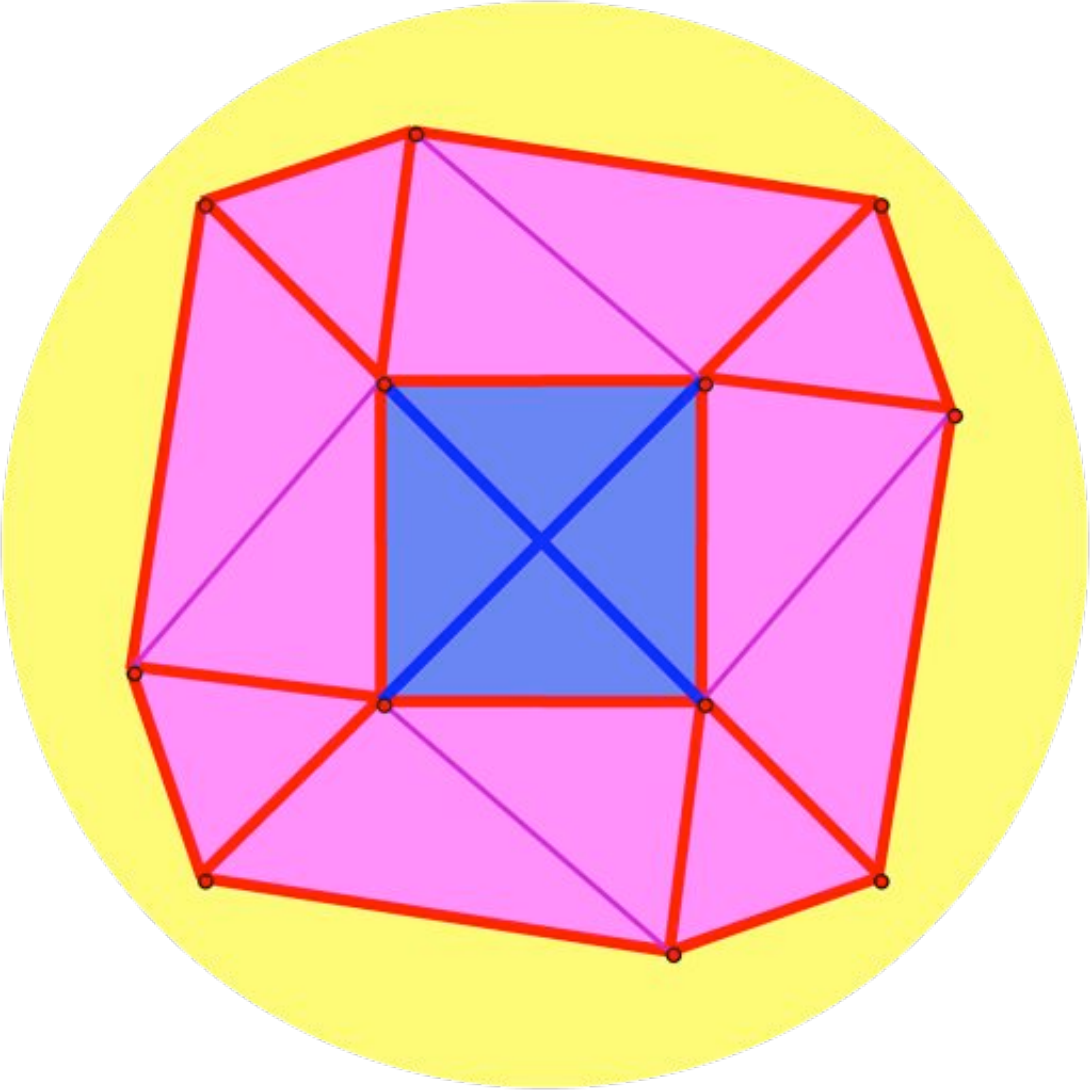, scale=.029}}
 \caption{Tubes on the cylinder $\cylinder(4,4,4)$ adding edges one at a time. }
\label{fig:C448}
 \end{center}
 \end{figure}
 \begin{itemize}
\item we can insert a sequence of $2(n-3)$ independent edges at the first hole hole (the number of edges needed to make this hole into a block);
\item each of these edges removes 1 idof from the overall framework, and therefore from the second hole.  
\item when this first hole $H_1$ is a block, the second hole, and the larger tube, still shows $(m-n)$ idof.  
\end{itemize}
The overall framework is independent, but not generically rigid.
All the $2(n-3)$ added independent constraints at $H_{1}$ {\it transmit} to reduced idof also at $H_{2}$. }
\end{example}

\begin{example}\label{example4}
{\rm Consider a cylinder $\cylinder(n, m, m)$ with $n>m$.  
\begin{itemize}
\item the two holes are missing $n+m-6$ edges, so the cylinder starts with  idof $=n+m-6$;
\item there are $n+m-6$ idof at the first hole (less than the number of edges needed to make this hole into a block).  
\item there is $2(m-3)$ idof at the second hole. 
\end{itemize}
 \begin{figure}[h]
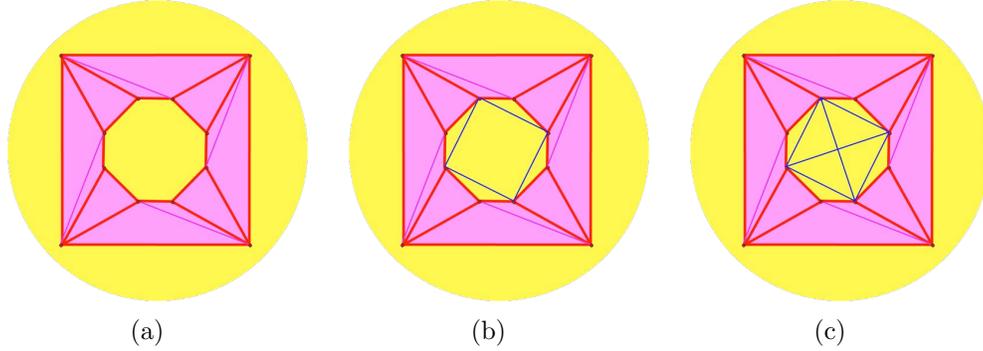

 \begin{center}
   \subfigure[] { \psfig{file=T844_0, scale=.027}}\quad
  \subfigure[]{ \psfig{file=T844_4, scale=.027}}\quad
  \subfigure[]{ \psfig{file=T844_6, scale=.027}}
 \caption{Tubes on the cylinder $\cylinder(8,4,4)$ adding edges one at a time up to rigidity.  We start with $6$ idof (a).  The first $4$ added edges make no change in the idof of the second hole, and the next two make the entire structure rigid.}
\label{fig:C844}
 \end{center}
 \end{figure}
\begin{itemize}
\item When we add the first $n-m$ independent edges to $H_{1}$, there is no change in the idof of $H_{2}$, which remains $2(m-3)$;
\item  Adding further $2(m-3)$ independent edges at $H_{1}$ makes the whole tube generically rigid.  
\item  We do not complete a block at $H_{1}$, unless we want the overall structure to be redundant.
\end{itemize}
The first $n-m$ edges at $H_1$ do not transmit any reduction in idof to $H_2$.  The next  $2(m-3)$ independent edges at $H_{1}$ transmit to reduced idof at $H_{2}$, finishing with idof$=0$. 
Any more edges at $H_{1}$ would be redundant and would not have any impact on the idof. }
\end{example}

\subsection{Examples with narrow waists}
\begin{example}\label{example5}
{\rm Consider a cylinder $\cylinder(n,k,m)$ with $k<n,m$ (Figure~\ref{fig:C8_4_8}).  We split the analysis by splitting the cylinder at some minimum waist $W$. 
 \begin{figure}[h]
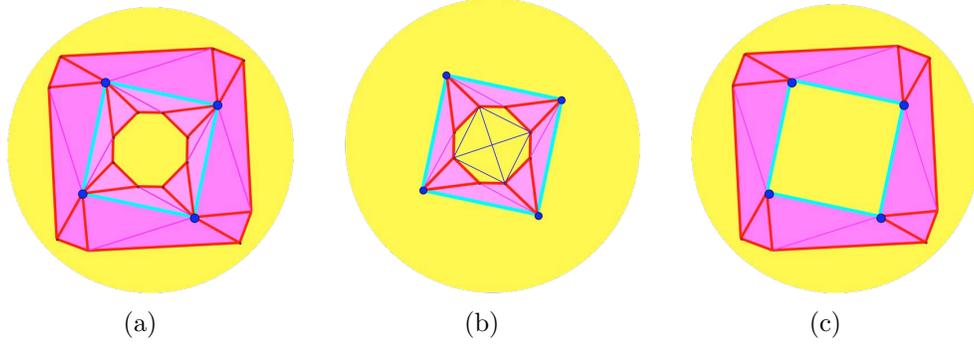

 \begin{center}
   \subfigure[] { \psfig{file=T848_0, scale=.027}} \quad
  \subfigure[]{ \psfig{file=T84-8, scale=.0285}} \quad
  \subfigure[]{ \psfig{file=T8-48, scale=.027}}
 \caption{A tube with a narrow waist can be separated into two pieces at the waist .}
\label{fig:C8_4_8}
 \end{center}
 \end{figure}
 Working with these simpler cylinders $\cylinder(n,k,k)$ and $\cylinder(k,k,m)$ (Figure~\ref{fig:C8_4_8}(b)(c)), we can track the overall transmission through the tube:
\begin{itemize}
\item the two holes are missing $(n+m-6)$ edges, so the cylinder starts with  idof $=(n+m-6)$;
\item the first component (containing $H_{1}$) behaves like Example~\ref{example4};
\item there are $(n+k-6)$ idof at the first hole $H_1$ (less than the number of edges needed to make this hole into a block);  
\item we can add up to $(n-k)$ edges at the first hole with no transmission to the waist (or beyond it);
\item the next $2(k-3)$ edges at the first hole transmit to rigidify the waist.  This is implicitly adding constraints to the waist, which transmits on to the second hole.  
\item The first component is now rigid, even though there is not a full block in the first hole $H_1$.  Any further added edges cause no reduction in idof at $H_{2}$.
\item In the second component, there is an implicit block at the waist and the second hole is still flexible with $(m-k)$ idof (as in Example~\ref{example3})
\end{itemize}
In this case we have transmission only for the added edges between $(n-k)$ and $(n+k-3)$}
\end{example}

We leave it to the reader to extract the summary lemmas on transmission which flow from these examples of tubes. 

In the original motivation from allostery, we mentioned `shape change at a distance'.  If adding an edge at $H_{1}$ transmitted to reducing the idof at $H_{2}$, then making a small change in the length of this bar, in a generic configuration, will produce some change in the distances between some pairs of vertices in $H_{2}$.  Thus transmission of reduced idof is equivalent to transmission of shape change from $H_{1}$ to $H_{1}$.

The larger message is that we can create mathematical models of frameworks  which: 
\begin{enumerate}
\item will connect two sites; 
\item will not transmit changed internal degrees of freedom until a certain number of constraints have been added at one site; 
\item will then faithfully transmit for the next set of added constraints; and 
\item will not completely constrain the other site, no matter what is done at the first site.  
\end{enumerate}
By connecting two `sites' by several tubes, or even interlocking three or more sites (which happens in allostery) we can create moderately complex models of `transmission' within a mathematical model for exploring allosteric transmission.  Whether something analogous occurs in the rigidity and flexibility of proteins is currently an active area of research. 

\section{Cycle to Cylinder Splits }\label{sec:cycle}

The  overall need for  added inductive constructions which preserve generic rigidity in $\R^{3}$ was highlighted in \cite{taywhiteley}.  Implicit in the results and proofs in \S4,5,6 are some  more general principles that can be applied to expand an arbitrary cycle $C_{k}$ in a graph $G$ to a `cylinder' while preserving rigidity related properties in $3$-space.   We summarize this process here. 

\begin{figure}[h]
 \begin{center}
   \subfigure[] { \psfig{file=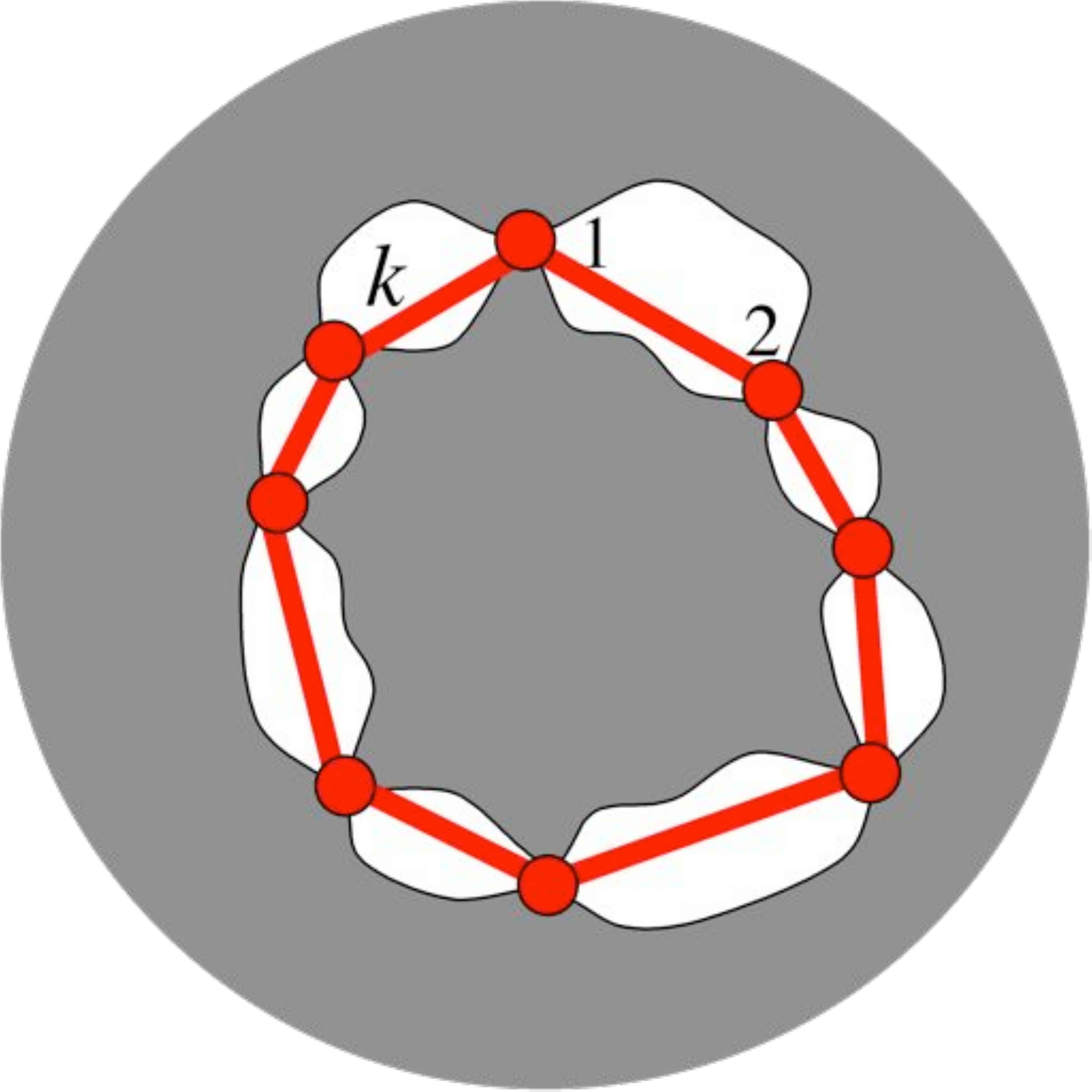, scale=.050}} \quad
    \subfigure[]{ \psfig{file=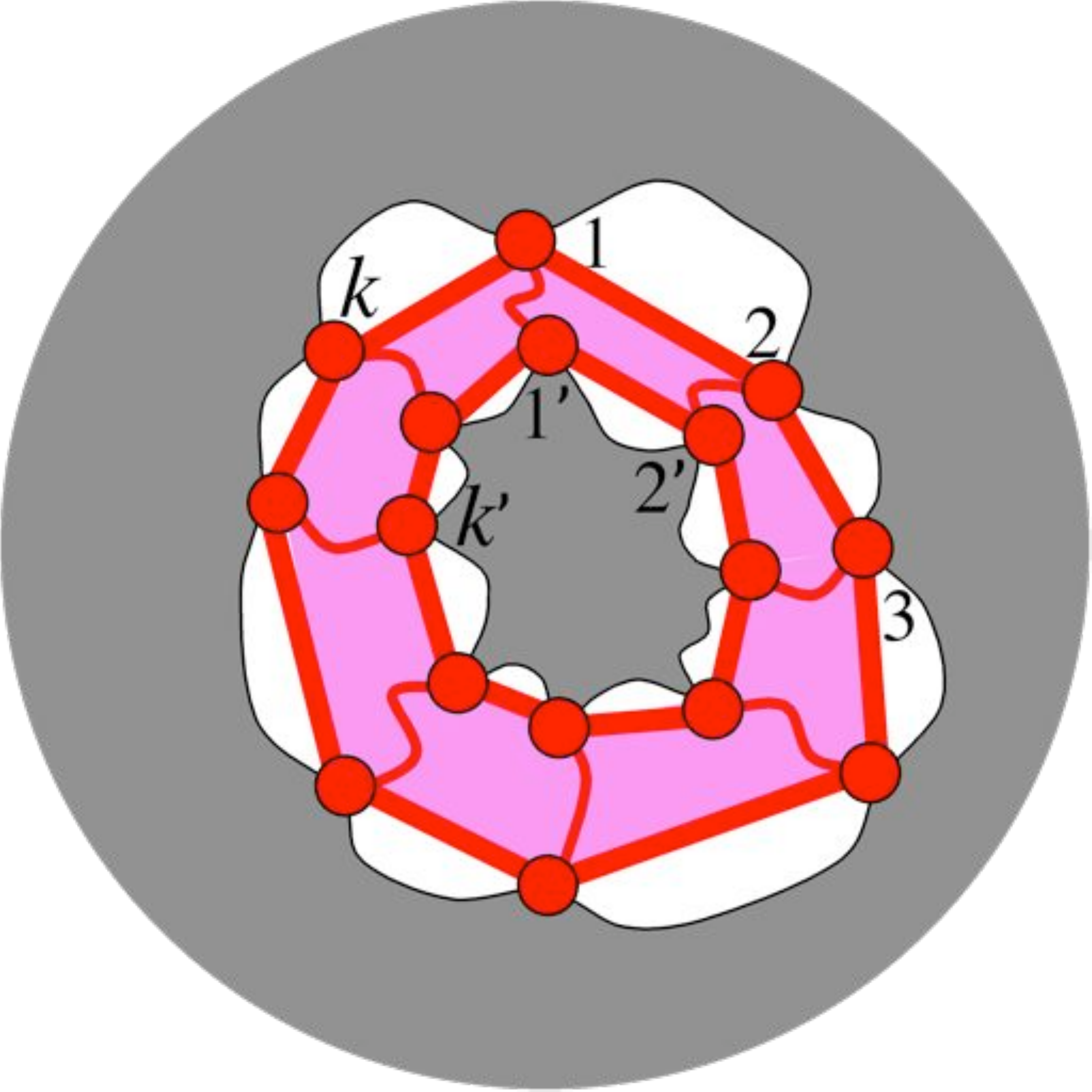, scale=.050}}\quad
       \subfigure[]{ \psfig{file=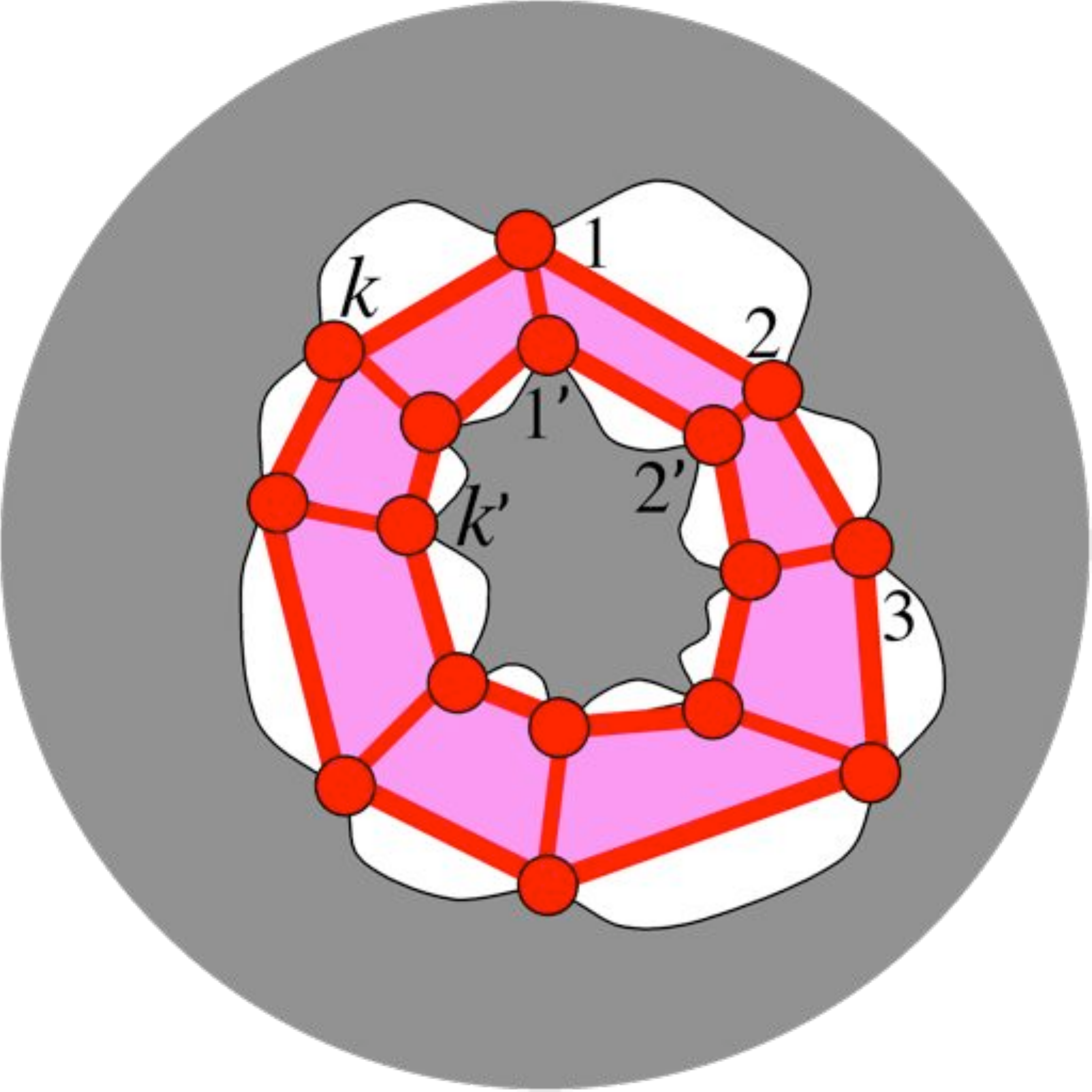, scale=.050}}
 \caption{Given a simple cycle $C_{k}$ in a graph (a), we can split along all  vertices and edges to two cycles, inserting a $k$-connected triangulated cylinder (b).  This preserves generic rigidity in $3$-space. }\label{fig:cyclesplit}
 \end{center}
 \end{figure}

\begin{definition}
{\rm   Given a graph $G = (V, E)$ with a $k$-cycle  $C_{k}=(1,2,\ldots,k)$ and subset of f the other edges at each of these vertices $S=(S_{1},S_{2},\ldots, S_{k})$,  then a {\it cycle split of $C_{k}$ on $(S_{1},S_{2},\ldots, S_{k})$} (Figure~\ref{fig:cyclesplit}(a)) is the extended graph $G * (C_{k},S)$ (Figure~\ref{fig:cyclesplit}(b)) with
\begin{itemize}
\item [(i)] a second cycle of   $C_{k}'=(1',2',\ldots,k')$; 
\item[(ii)] all edges in $S_{i}$ are removed from $i$ and attached to $i'$; and 
\item[(iii)] a $k$-cylinder $\cylinder(k,k,k)$ is inserted between the cycles $C_{k}$ and $C_{k}'$. 
 \end{itemize}
Note that it is permissible that some $i'=i$:  that is  we split only some of the vertices, but otherwise the $i'$ are new vertices. }  
\end{definition}

\begin{theorem} Gven a generically isostatic (generically independent, generically rigid) graph $G$ with a $k$-cycle $C$ with separations $S$, then:
\begin{itemize}
\item[(a)]  any  cycle split  $G * (C_{k},S)$  is obtained from $G$ by a sequence of vertex splits;
\item[(b)]  any  cycle split  $G * (C_{k},S)$ is generically isostatic (generically independent, generically rigid, resp.).
\end{itemize}
\end{theorem} \label{thm:cyclesplit}

\proof  This is implicitly presented in the discussions in Sections 5 and 6. We use the $k$-connectivity to add $k$ vertex-disjoint paths (Figure~\ref{fig:cyclesplit}(b)).   We just note that those arguments {\it reserved} the `faces' of the two cycles in the sense that we do not change any edges around the cycles when simplify the paths to length 1 (Figure~\ref{fig:cyclesplit}(c)).  That is, what other edges there are on the vertices of those cycles do not alter to contractions and the construction of the sequence of vertex splits, because the boundaries of the cycles are `clean' in that analysis. 

Notice that if there were other edges between vertices the cycle $C_{k}$ before the split, they can be left on one side, or the other, or separated to run from one side to the other.  In all cases, because the connected non-adjacent vertices of the cycle, they will connect vertices which do not share a face after the split and the subdivision of the cylinder into $k$ faces.  As such, they cannot make any face boundary unclean and will not impact any of the reductions above. 

For the contraction of the remaining paths of length $1$ between the cycles, we use simply do the contraction and rely on the fact that the original graph $G$ was isostatic, or independent, or rigid as needed. 
\qed

In \S6, we were implicitly using this result for independence when we considered only partial bracing at $H_{1}$. 
\medskip

\begin{figure}[h]
 \begin{center}
   \subfigure[] { \psfig{file=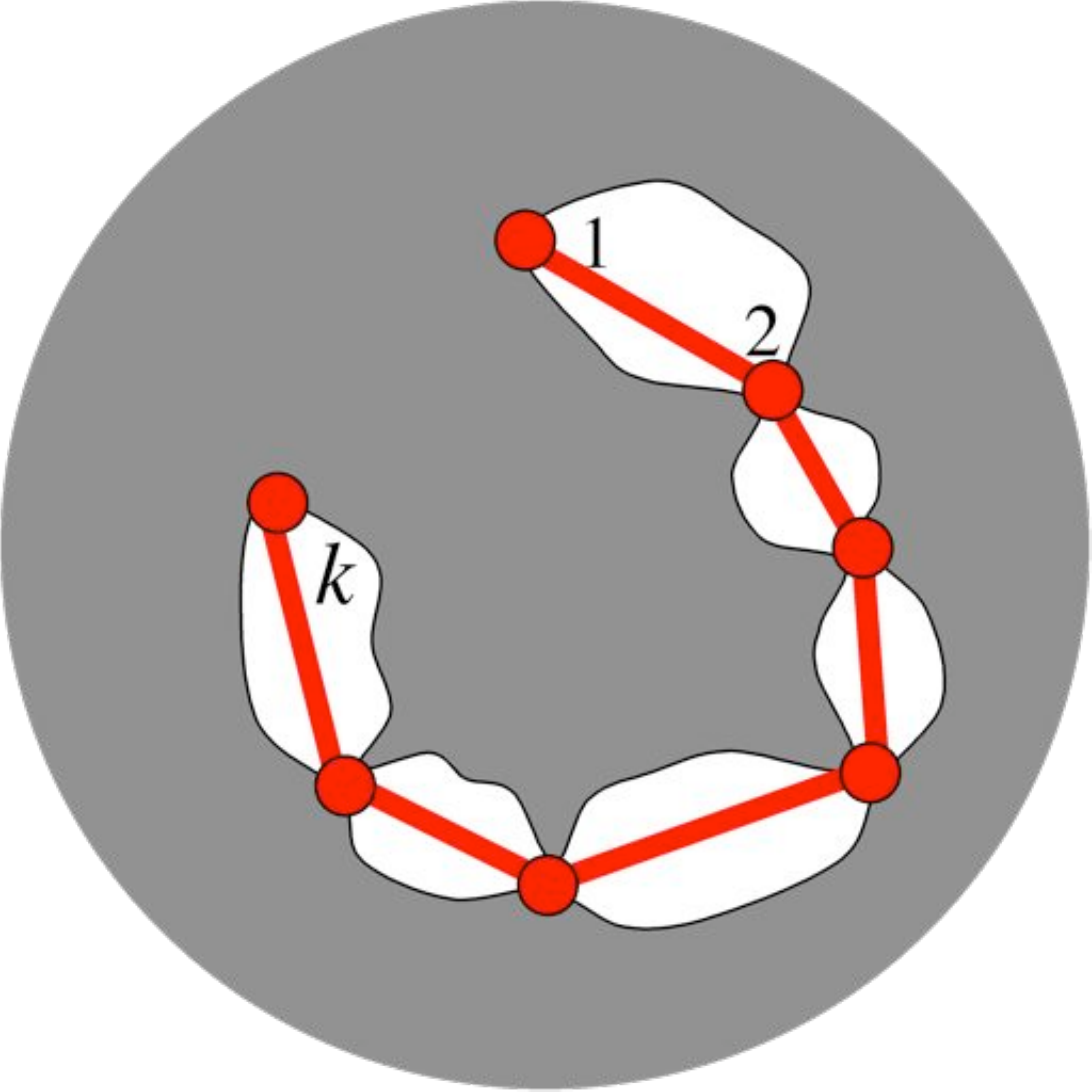, scale=.065}} \quad
    \subfigure[]{ \psfig{file=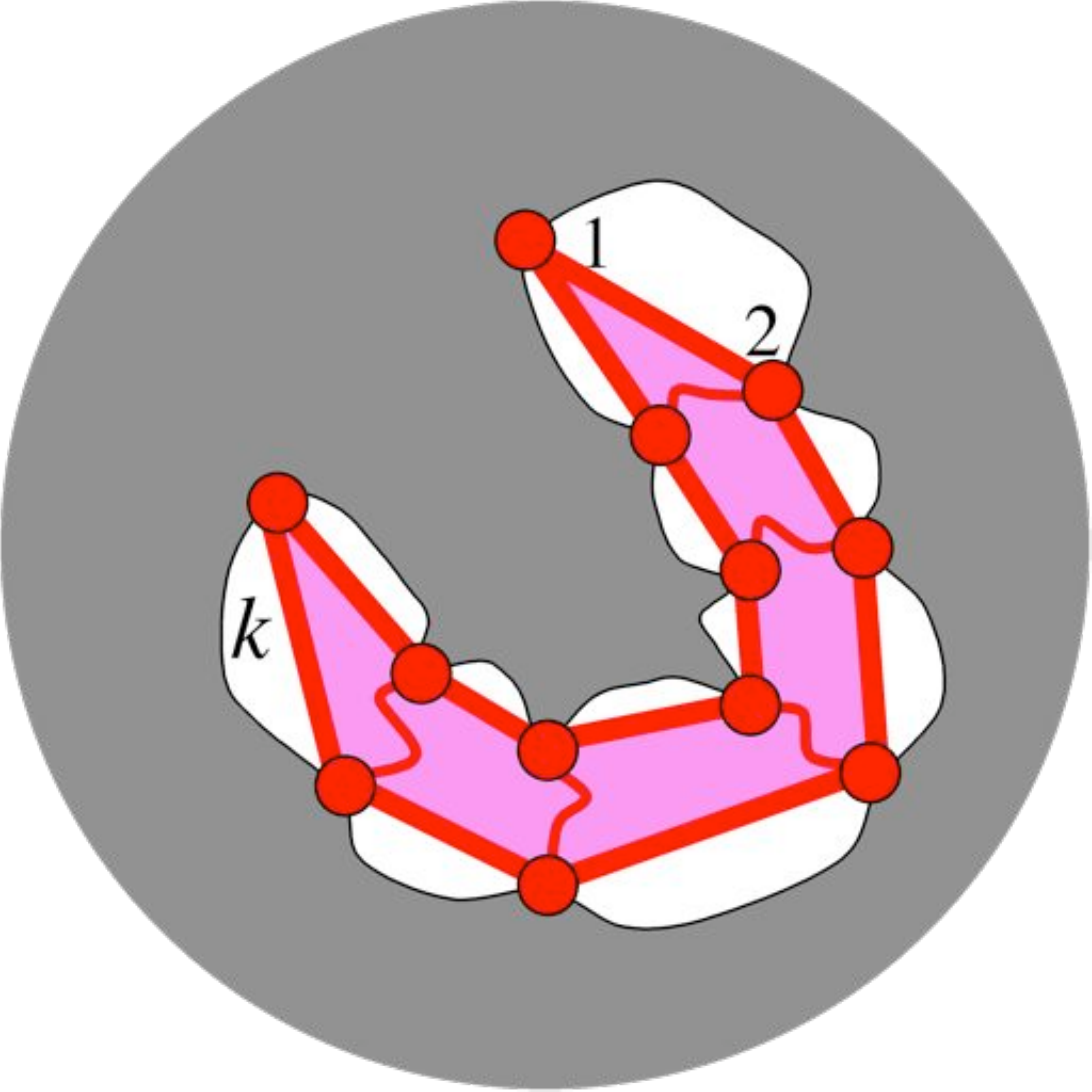, scale=.065}}
 \caption{Given a simple path of length $k$ (a), we can split along the path to create two paths, joined at the end, inserting a $k$-connected fence (b).  This preserves generic rigidity in $3$-space. }\label{fig:pathsplit}
 \end{center}
 \end{figure}

We can extend this reasoning to more general insertions in an arbitrary graph $G$.  For example, consider a {\it path split} in which the cycle is replaced by a path $P = (1,2,\ldots,k)$ of distinct vertices without an edge joining $!,k$.   We now only duplicate vertices $2,\ldots,k-1$ with $2',\ldots,k-1'$ leaving the first and last vertices not split.  We need selections $S=(S_{2},\ldots, S_{k-1})$ to identify which edges move from the original vertices to the duplicate vertices, and we duplicate the entire path (Figure~\ref{fig:pathsplit}(a),(b)).  We insert a triangulated disc in which we can find $k-2$ vertex-disjoint paths between the two paths (also disjoint from the end points) Figure~\ref{fig:pathsplit} (b). 

By an analogous argument, these connecting paths guarantee that this insertion can be generated by a sequence of vertex splits. Therefore the path split preserves generic rigidity and independence.  We observe that the simplest split on a path of length $2$ is the original vertex split. 

\section{Further work}\label{sec:further}
The results and methods in this paper open up some possible extensions and some new questions.  We note a few here. 

\subsection{Spheres with shafts} \label{sec:shaft}
The original conjecture of Kuiper \cite{kuiper,infp2} was that removing one edge from a $4$-connected convex triangulated sphere would leave a finite motion in which the 
dihedral angle at each remaining edge was changing.  This is equivalent to saying that the shaft between the two vertices of the triangles connected to this edge would return 
the framework to infinitesimal rigidity.  This is not quite true for all geometric triangulated convex spheres, but the results in \S 5 verify this is generically true. 
We can go on to investigate more general circumstances when the shaft connects some two non-adjacent vertices.  
With a few details to check, the indication is that these spheres with long shafts are also generically rigid circuits, if the extended graph is $4$-connected.  
We plan to complete this investigation in another paper.  

Jackson and Jordan showed that in 2-dimensions, vertex splitting resulting in both the added vertex and the split-vertex being at least $3$-valent, 
takes a generic rigid circuit (minimal dependent rigid sets) to a $2$-circuit, and a generically globally rigid framework to a generically globally rigid framework.  
We are interested in knowing if vertex splitting preserves generically rigid $3$-circuits, and more generally global rigidity in $3$-dimensions. 

\begin{conjecture}[Connelly, Whiteley \cite{conwhit}] Vertex splitting resulting in both the added vertex and the split-vertex being at least 4-valent preserves the generic rigidity of
$3$-circuits, and global rigidity in $3$-space.  \end{conjecture}  

In many cases (but not all), generically rigid circuits in $3$-space are also globally rigid.  These spheres with shafts provide an example where available methods 
prove that vertex splitting creates new generically rigid $3$-circuits, using the observation that in this class, vertex splitting preserves $4$-connectivity. 
Current unpublished work of Connelly and of Whiteley then shows that these vertex splits preserves global rigidity.

\subsection{Change of metric} \label{sec:metric}
In this paper, the results were presented as applying in Euclidean $3$-space, as were the original results of Cauchy, Dehn, and Alexandrov \cite{cauchy,dehn,alexandrov}.  However, the techniques were clearly based on the combinatorics of the graphs, aided by the underlying topology of the sphere. 

 It is not surprising that the results (as well as the results of Cauchy, Dehn, and Alexandrov) apply to the range of Cayley-Klein geometries such as the $3$-sphere, the hyperbolic $3$-space, or the Minkowski $3$-space \cite{saliola}.  We could verify this by carefully translating the methods, in particular verifying that vertex-splitting preserves generic rigidity in each of the metrics.  Alternatively, the results in other metrics are clear because of the general transfer of generic rigidity results among these metrics,  as presented in  \cite{saliola}.

\subsection{Symmetries on Block and Hole frameworks} \label{sec:symmetry}
A number of recent papers have studied the impact of symmetry on the rigidity of frameworks \cite{cfgsw,BS4,BS6,BS2,BS1,BS3}.  In particular, these have developed a theory of 
frameworks which are as generic as possible within the constraints of the symmetry.  

In many cases, the symmetry has no impact on the rigidity of a framework.  It is now natural to ask whether the results of this paper can be symmetrized. 
That is, if the base polyhedron is symmetric and isostatic, and the proposed final block and hole polyhedron has the same symmetry and is generic within this symmetry, can we find a set of symmetric vertex splittings (vertex splittings which are applied separately in a symmetric fashion) which preserved symmetry generic rigidity? 

This must be done delicately, as there are examples of isostatic base polyhedra with mirror symmetry, and larger block and hole polyhedra with the same symmetries that are
flexible \cite{tarnai}.  The required vertex splits cannot be applied in a symmetric sequence to preserve the symmetry at each stage, and they add vertices or edges which are fixed by at least one symmetry (that is the vertex goes to itself, or the edge pair of vertices goes to the same pair).  We offer one conjecture for cases where we anticipate this will work. 

\begin{conjecture}
Given an abstract block and hole polyhedron with symmetry group on the graph $\Phi(\mathcal{S})$ with no vertices or edges fixed, there is a symmetry preserving set of contractions to a base polyhedron realized with the same  symmetry group $\mathcal{S}$ (no vertices or edges fixed).  If the base polyhedron is isostatic in some symmetric realization, then the expanded block and hole polyhedron is isostatic when realized generically within the symmetries $\mathcal{S}$.
\end{conjecture}

\providecommand{\bysame}{\leavevmode\hbox to3em{\hrulefill}\thinspace}

\end{document}